\documentclass[a4paper]{article}
\usepackage[latin1]{inputenc}
\usepackage[T1]{fontenc}
\usepackage{RR}
\usepackage{hyperref}
\usepackage{amsthm}
\usepackage{amsmath}
\usepackage{dsfont}
\usepackage{amssymb}
\usepackage{graphicx}
\usepackage{tabularx}
\usepackage{pstricks}
 \usepackage{pst-grad} 
\usepackage{algorithm} 
\usepackage{algorithmic}
\newtheorem{theo}{Theorem}

\newtheorem{rem}[theo]{Remark}

\newcommand{\ds}{\displaystyle}

 \newcommand{\V}{\mathbf} 
 \newcommand{\algoref}[1]{Algo. (\ref{#1})} %
\newcommand{\fref}[1]{fig. \ref{fig:#1}} %
\newcommand{\cref}[1]{chapitre \ref{cha:#1}} %
\newcommand{\sref}[1]{sec. \ref{sec:#1}} 
\newcommand{\eq}{\begin{equation}} 
\newcommand{\eeq}{\end{equation}}
\newcommand{\eqs}{\begin{equation*}} 
\newcommand{\eeqs}{\end{equation*}}  	

\newcommand{\Ea}{\begin{eqnarray}} 
\newcommand{\Eas}{\begin{eqnarray*}} 
\newcommand{\Eae}{\end{eqnarray}}
\newcommand{\Ease}{\end{eqnarray*}}
\newcommand{\Ean}{\begin{eqnarray*}} 
\newcommand{\Eane}{\end{eqnarray*}}
\newcommand{\Al}{\begin{align}} 
\newcommand{\E}{\end{align}}
\newcommand{\Bea}{\begin{Beqnarray}} 
\newcommand{\Beae}{\end{Beqnarray}}

\graphicspath{{./}}

\newcommand{\dparx}[1]{\displaystyle \frac{\partial #1}{\partial x}} 
\newcommand{\dparxk}[1]{\frac{\partial #1}{\partial x_k}} 
\newcommand{\deriv}[2]{\frac{d #1}{d #2}} 
 
\RRdate{November 2010}

\RRauthor{
 J\'er\^ome Guterl%
  \thanks[stras]{IRMA, Universit\'e de Strasbourg, 7 rue Ren\'e-Descartes, 67084 Strasbourg Cedex}  
  \and
Jean-Philippe Braeunig%
  \thanks[inria]{INRIA Nancy-Grand Est, 615 rue du Jardin Botanique, 54600 Villers-l\`es-Nancy}%
  \thanksref{stras}
  \thanks{CEA/DIF Bruy\`eres-le-Ch\^atel, 91297 Arpajon Cedex}\\
  \and
Nicolas Crouseilles 
 \thanksref{inria}
  \thanksref{stras}
  \and Virginie Grandgirard
 \thanks[cad]{CEA Cadarache, 13108 St Paul-lez-Durance Cedex}
  \and Guilaume Latu
   \thanksref{cad}
   \thanksref{inria}    
\and  
  Michel Mehrenberger\thanksref{stras}  
   \thanksref{inria}
   \and Eric Sonnendr\"{u}cker\thanksref{stras}  
   \thanksref{inria}
}
\authorhead{Guterl et al}
\RRtitle{Test of some numerical limiters for the conservative PSM scheme for 4D Drift-Kinetic simulations.}
\RRetitle{Test of some numerical limiters for the conservative PSM scheme for 4D Drift-Kinetic simulations.}
\titlehead{Test of some numerical limiters for the PSM scheme}
\RRresume{Ce travail concerne la simulation de plasmas magn\'etis\'es dans le cadre du projet ITER.  Pour cette application, des mod\`eles de type Vlasov-Poisson sont utilis\'es pour simuler la turbulence \`a coeur dans un tokamak, en g\' eom\'etrie toroidale. Ces \'etudes m\`enent \`a r\'esoudre des probl\`emes dans un espace \`a 6 dimensions, 3D en espace 3D en vitesse, qui sont tr\`es lourds \`a simuler en terme de ressources informatiques. Le mod\`ele est r\'eduit \`a un mod\`ele gyrocin\'etique 5D en exploitant les caract\'eristiques de ce plasma, dont le mouvement des particules est fortement influenc\'e par la pr\'esence d'un champ magn\'etique intense. Cependant, il est n\'ecessaire de mettre au point des sch\'emas pr\'ecis et des algorithmes parall\`eles pour mener ces simulations. Ce rapport d\'ecrit une formulation de type Hermite du sch\'ema conservatif PSM qui est tr\`es g\'en\'erique et qui permet d'impl\'ementer diff\'erent sch\'emas semi-Lagrangiens. Nous testons et proposons \'egalement des limiteurs num\'eriques de pente qui doivent accro\^itre la robustesse des simulations en r\'eduisant les oscillations d'origine num\'erique. Dans ce travail, nous l'utilisons pour r\'esoudre le mod\`ele drift-kinetic 4D, qui est le squelette du mod\`ele gyrocin\'etique 5D. Ce mod\`ele 4D est suffisamment pertinent pour la conception d'une m\'ethode num\'erique robuste et pr\'ecise pour le mod\`ele  5D}
\RRabstract{The purpose of this work is simulation of magnetised plasmas in the ITER project framework. In this context, Vlasov-Poisson like models are used to simulate core turbulence in the  tokamak in a toroidal geometry. This leads to heavy simulation because a 6D dimensional problem has to be solved, 3D in space and 3D in velocity. The model is reduced to a 5D gyrokinetic model, taking advantage of the particular motion of particles due to the  presence of a strong magnetic field. However, accurate schemes, parallel algorithms need to be designed to bear these simulations. This paper describes a Hermite formulation of the conservative PSM scheme which is very generic and allows to implement different semi-Lagrangian schemes. We also test and propose numerical limiters which should improve the robustness of the simulations by diminishing spurious oscillations. We only consider here the 4D drift-kinetic model which is the backbone of the 5D gyrokinetic models and relevant to build a robust and accurate numerical method.}
\RRmotcle{simulation num\'erique, sch\'ema conservatif, ITER, turbulence plasma}
\RRkeyword{numerical simulation, conservative scheme, ITER, plasma turbulence}
\RRprojet{CALVI}
\RRdomaine{1} 
\RRtheme{Mod\'elisation, analyse num\'erique}
 \URLorraine 
\RCNancy 

\begin{document}
\RRNo{7467}
\makeRR   

\newpage

\tableofcontents

\newpage
\section{Introduction}

 The ITER device is a tokamak designed to study controlled thermonuclear fusion. Roughly speaking, it is a toroidal vessel containing a magnetized plasma where fusion reactions occur. The plasma is kept out of the vessel walls by a magnetic field which lines have a specific helicoidal geometry. However, turbulence develops in the plasma and leads to thermal transport which decreases the confinement efficiency and thus needs a careful study.  Plasma is constituted of ions and electrons, which motion is induced by the magnetic field. The characteristic mean free path is high, even compared with the vessel size, therefore a kinetic description of particles is required, see {\it Dimits} \cite{dimits}. Then a full 6D Vlasov-Poisson model should be used for both ions and electrons to properly describe the plasma evolution. However, the plasma flow in presence of a strong magnetic field has characteristics that allow some physical assumptions to reduce the model. First, the Larmor radius, i.e. the radius of the cyclotronic motion of particles around magnetic field lines, can be considered as small compared with the tokamak size and the gyration frequency very fast compared to the plasma frequency. Thus this motion can be averaged (gyro-average) becoming the so-called guiding center motion. As a consequence, 6D Vlasov-Poisson model is reduced to a 5D gyrokinetic model by averaging equations in such a way the 6D toroidal coordinate system $(r,\theta,\phi,v_r,v_\theta,v_\phi)$ becomes a 5D coordinate system $(r,\theta,\phi,v_\parallel, \mu)$, with $v_\parallel$ the parallel to the field lines component of the velocity and  $\mu=m~v^2_\perp / 2B$ the adiabatic invariant which depends on the norm of the perpendicular to the field lines components of the velocity $v^2_\perp$, on the magnetic field magnitude $B$ and on the particles mass $m$. Moreover, the magnetic field is assumed to be steady and the mass of electrons $m_e$ is very small compared to the mass of ions $m_i$. Thus  the cyclotron frequency $\omega_{i,e}=q_{i,e} ~B / m_{i,e}$ is assumed to be much higher for electrons than for ions  $\omega_{e} >> \omega_{i} $. Therefore the electrons are assumed to be at Boltzmann equilibrium, i.e. the effect of the electrons cyclotronic motion is neglected. The 5D gyrokinetic model then reduces to a Vlasov like equation for ions guiding center motion: 
\begin{equation}\label{gyro_vlasov}
\begin{array}{ll} 
\dfrac{\partial \bar{f}_{\mu}}{\partial t} + \dfrac{dX}{dt} \cdot \nabla_X \bar{f}_{\mu} + \dfrac{dv_{\parallel}}{dt}\partial_{v_{\parallel}}\bar{f}_{\mu} = 0  
\end{array}
\end{equation}
where $\bar{f}_{\mu} (X,v_\parallel)$ is the ion distribution function with $X=(r,\theta,\phi)$, velocities $dX/dt$ and ${dv_{\parallel}}/{dt}$ define the guiding center trajectories.\\
If $\nabla_{(X,v_{\parallel})} \cdot (dX/dt,{dv_{\parallel}}/{dt})^t =0$, then the model is termed as conservative.\\
This equation for ions is coupled with a quasi-neutrality equation for the electric potential $\Phi(R)$ on real particles position, with $R=X-\rho_L$ (with $\rho_L$ the Larmor radius) :
\begin{equation}\label{gyro_poisson}
\begin{array}{ll} 
-\dfrac{1}{B \omega_{i} } \nabla_\perp \cdot (n_0 \nabla_\perp \Phi) +\dfrac{e}{\kappa T_e}(\Phi-<\Phi>_{\phi})= \int \bar{f}_{\mu} d\mu dv_{\parallel}-n_0 
\end{array}
\end{equation}
where $n_0$ is the equilibrium electronic density, $T_e$ the electronic temperature, $e$ the electronic charge, $\kappa$ the Boltzmann constant for electrons and $ \omega_{i} $ the cyclotronic frequency for ions.\\

These equations are of a simple form, but they have to be solved very efficiently because of the 5D space and the large characteristic time scales considered. However, the adiabatic invariant $\mu$ acts as a parameter, thus it could easily be parallelized. Moreover, we can see that for each $\mu$, we have to solve a  4D advection equation, as accurately as possible but also taking special care on mass and energy conservation, especially in this context of large characteristic time scales. The maximum principle that exists at the continuous level for the Vlasov equation should also be carefully studied at discrete level because there is no physical dissipation process in this model that might dissipate over/undershoots of the scheme. These studies will be achieved first on a relevant reduced model, the 4D drift-kinetic model which corresponds to \eqref{gyro_vlasov} with $\mu=0$ instead a range of $\mu$ values (theoretically $\mathbb{R}^+$).This work follows those of {\it Grandgirard et al} in the GYSELA code, see \cite{gysela} and \cite{grandgirard}.
The geometrical assumptions of this model for ion plasma turbulence are a cylindrical geometry with coordinates $(r,\theta,z,v_{\|})$ and a constant magnetic field $B=B_z ~e_z$,  where $e_z$ is the unit vector in $z$ direction. In this collisionless plasma, the trajectories are governed by the guiding center (GC) trajectories:
\begin{equation}
\begin{array}{lll}
\dfrac{dr}{dt}=v_{GC_r}; ~r \dfrac{d \theta}{dt}=v_{GC_\theta}; ~ \dfrac{d z}{dt}=v_\|; ~\dfrac{d v_\|}{dt}=\dfrac{q_i}{m_i} E_z
\end{array}
\end{equation}
with $v_{GC}=(E \times B)/B^2$ and $E=-\nabla \Phi$ with $\Phi$ the electric potential.\\
The Vlasov equation governing this system, with the ion distribution function $f(r,\theta,z,v_\|,t)$, is the following:
\begin{equation}\label{Vlasov4d}
\begin{array}{lll}
\partial_t f+v_{GC_r} \partial_r f + v_{GC_\theta} \partial_\theta f  + v_\| \partial_z f  +\dfrac{q_i}{m_i} E_z \partial_{v_\|} f =0.
\end{array}
\end{equation}
This equation is coupled with a quasi-neutrality equation for the electric potential $\Phi(r,\theta,z)$ that reads the same as for the 5D gyrokinetic model  \eqref{gyro_poisson} with $\mu=0$.

Let us notice that the 4D velocity field $a=(v_{GC_r},v_{GC_\theta},v_\|, q/m_i ~E_z)^t$ is divergence free:
\begin{equation}\label{diva}
\begin{array}{lll}
\nabla \cdot a = \dfrac{1}{r} \partial_r  (r~v_{GC_r})+  \dfrac{1}{r} \partial_\theta  (v_{GC_\theta})+\partial_z  v_\| + \partial_{v_\|} (q_i/m_i ~E_z)=0
\end{array}
\end{equation}
because of variable independence $\partial_{v_\|}  E_z=\partial_{v_\|} (\partial_z \Phi(r,\theta,z))=0$ and $\partial_z  v_\|=0$. Moreover we have $v_{GC}=(E \times B)/B^2$, with  $E=-\nabla \Phi$ and $B=B_z ~e_z$, thus: 
\begin{equation}\label{def_vGC}
v_{GC_r}=\dfrac{1}{B_z} \left( -\dfrac{1}{r} \partial_\theta \Phi \right) ~~ \mbox{and} ~~ v_{GC_\theta}=\dfrac{1}{B_z} \left( \partial_r \Phi \right)
\end{equation}
and
\begin{equation}\label{div_rtheta}
\nabla_{r \theta} \cdot a = \dfrac{1}{r}  \partial_r  (r~v_{GC_r})+ \dfrac{1}{r}  \partial_\theta  (v_{GC_\theta})=\dfrac{1}{r~B_z}  \left( \partial_r  \left( r~ (-1/r) \partial_\theta \Phi \right) +  \partial_\theta \left( \partial_r \Phi \right) \right)=0.
\end{equation}
Therefore, one can write an equivalent conservative equation to the preceding Vlasov equation (\ref{Vlasov4d}):
\begin{equation}\label{Vlasovc}
\begin{array}{lll}
\partial_t f+\partial_r (v_{GC_r} ~f)+  \partial_\theta (v_{GC_\theta} ~f)  +  \partial_z (v_\| ~f)  + \partial_{v_\|} \left( \dfrac{q_i}{m_i} E_z ~f \right) =0.
\end{array}
\end{equation}

This conservative system will be discretized using a conservative semi - Lagrangian scheme. Following \cite{braeunig} and \cite{PSM1}, we consider two conservative  schemes, which are fourth order in space:
\begin{itemize}
\item \textbf{LAG}: LAGrange polynom method, which uses Lagrangian polynoms to reconstruct the distribution  function.
\item \textbf{PSM}: Parabolic Splines Method, which uses cubic splines to reconstruct the distribution function.
\end{itemize}
These schemes are designed to solve conservative models and they allow a directional splitting of the Vlasov equation \eqref{Vlasovc}. This equation will be solved by using $D$ (dimensions of space) 1D conservative steps, discretized by using the 1D schemes (LAG or PSM). At the continuous level, each 1D step has no maximum principle, it is only the solution after all $D$ directional steps, the solution of the Vlasov equation, that should satisfy a maximum principle \cite{braeunig}. However, high order schemes may create spurious oscillations leading to break this maximum principle. A flux limiting procedure may improve the discrete solution, which may be closer to the maximum principle in the sense of showing less spurious oscillations. The limiter does not ensure a maximum principle in 1D, but should decrease the oscillations amplitude created by the scheme by locally adding numerical diffusion.  That leads us to investigate and compare in details the properties of many limiters, which depend mainly on two issues:
\begin{itemize}
\item How to make the high order schemes degenerate into a  more diffusive scheme?
\item How to detect in the function profile the location where the scheme will produce oscillations? 
\end{itemize}
We first describe a Hermite formalism proposed by \cite{crouseilles} applied to the PSM and LAG schemes leading to a finite volume form equivalent to the original semi-Lagrangian schemes. This formalism is very efficient to introduce limiters in the PSM or LAG schemes. Consequently, we will compare some limiters, focusing in particular on the OScillations limiter (OSL) proposed by \cite{crouseilles}. We also propose and investigate a new limiter (Slope Limited Spline, SLS), based on classical slope limiting methods. We evaluate the performances of each limiter using a benchmark developed in the Gysela code, which runs a 4D drift-kinetic model \cite{grandgirard}.\\   
The outline of this paper is the following : in section \ref{sec:vlasov} will be recalled some important properties concerning the conservative form of the Vlasov equations. Then the Hermite formalism applied to LAG and PSM schemes will be explored. In section \ref{sec:limiter}, some limiters are described and they  are further investigated in section \ref{sec:numerical_results} in the context of a 4D drift-kinetic model. At last we will comment on numerical results.
\section{Numerical schemes for the Vlasov equation}
\label{sec:vlasov}
\subsection{Directional splitting of the advection problem}
In a phase space of dimension $D$, we consider a distribution $f$ which is advected by a velocity field $\V{a}$. The model taken into account satisfies $\nabla \cdot   \V{a} =0$.
\eq
t\in\mathbb{R}^{+},\V{x}\in\mathbb{R}^D, \V{a}(\V{x},t)\in\mathbb{R}^D
\begin{cases}
\partial_t  f+   \nabla_x\cdot (\V{a}f)= 0\\
\nabla \cdot   \V{a} =0\\
f(x,t) \ge 0
\end{cases}
\label{probDD}
\eeq

For instance in cylindrical geometry, $\V{x}=\left(r,\theta,z,v_{\parallel} \right)$ considering the 4D problem dealt by the Gysela code.
In the next sections, we use a directional splitting following \cite{braeunig} by solving the conservative system \eqref{probDD} by D separate 1D problems for each phase space direction which are still under a conservative form. So formally, we will consider the problem \eqref{eq:equ2} for each of the D directions. The generic direction is named $x$.
\eq
t\in\mathbb{R}^{+},x\in\mathbb{R}, a(x,t)\in\mathbb{R}
\begin{cases}
\partial_t  f+   \displaystyle \dparx{\left( a(x,t)f(x,t)\right)}= 0\\
f(x,t) \ge 0
\end{cases}
\label{eq:equ2}
\eeq
In this context, we don't have in general $\forall k \in [\mspace{-2 mu} [  1,D   ]\mspace{-2 mu}] ,\dparxk{a(x,t)}=0$, but only $\nabla \cdot   \V{a} =0$.

\subsection{Distribution function and phase space}
We divide one direction of the phase space, generically  $x$ with a constant step $\Delta x$ to get a regular mesh. The cells are numbered by an integer i from 0 to N and the cell faces by an one-half integer $i\pm1/2$ (see \fref{grid0}). Hence we have N+1 cells and N+2 faces.

 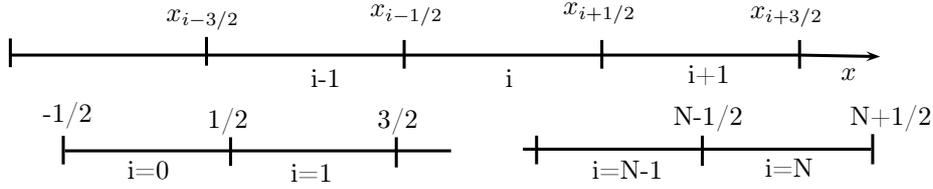
\begin{figure}[!h]
 \centering
 \begin{pspicture}(0,-0.595)(11.44,0.595)

 \rput(7.7521877,0.415){$x_{i+1/2}$}

 \rput(5.2121873,0.395){$x_{i-1/2}$}

 \rput(6.5603123,-0.445){i}

 \rput(9.182656,-0.405){i+1}

 \rput(4.1426563,-0.445){i-1}

 \rput(2.5321875,0.335){$x_{i-3/2}$}

 \rput(10.212188,0.375){$x_{i+3/2}$}
 \psline[linewidth=0.04cm](0.0,-0.135)(10.38,-0.135)
 \psline[linewidth=0.04cm](10.38,0.065)(10.38,-0.315)
 \psline[linewidth=0.04cm](7.78,0.065)(7.78,-0.315)
 \psline[linewidth=0.04cm](5.18,0.065)(5.18,-0.315)
 \psline[linewidth=0.04cm](2.58,0.105)(2.58,-0.275)
 \psline[linewidth=0.04cm](0.0,0.085)(0.0,-0.295)
 \psline[linewidth=0.04cm,arrowsize=0.05291667cm 2.0,arrowlength=1.4,arrowinset=0.4]{->}(10.38,-0.135)(11.44,-0.155)

 \rput(11.026563,-0.385){$x$}
 \end{pspicture}
 \begin{tabular}{c}
 $\:$
  \end{tabular}
 \centering
 \begin{tabular}{cc}
 \begin{pspicture}(0,-0.545)(5.3253126,0.545)
 \psline[linewidth=0.04cm](0.2253125,-0.105)(5.3053126,-0.105)
 \psline[linewidth=0.04cm](0.2053125,0.095)(0.2053125,-0.285)
 \psline[linewidth=0.04cm](2.4053125,0.075)(2.4053125,-0.305)
 \psline[linewidth=0.04cm](4.5853124,0.075)(4.5853124,-0.325)

 \rput(3.4679687,-0.395){i=1}

 \rput(1.3235937,-0.355){i=0}

 \rput(0.26046875,0.365){-1/2}

 \rput(2.3667188,0.285){1/2}

 \rput(4.6004686,0.285){3/2}
 \end{pspicture} & \begin{pspicture}(1.5,-0.545)(6.8253126,0.545)
 \psline[linewidth=0.04cm](2.0,-0.083125)(6.58,-0.083125)
 \psline[linewidth=0.04cm](2.18,0.096875)(2.18,-0.283125)
 \psline[linewidth=0.04cm](4.36,0.096875)(4.36,-0.303125)

 \rput(5.5021877,-0.333125){i=N}

 \rput(3.3726563,-0.353125){i=N-1}

 \rput(4.4348435,0.326875){N-1/2}

 \rput(6.8348436,0.326875){N+1/2}
 \psline[linewidth=0.04cm](6.6,0.116875)(6.6,-0.283125)
 \end{pspicture} 
 \end{tabular}
 \caption{Mesh grid on the $x$ phase space direction (top) and beginning (left bottom) and end (right bottom) of the mesh grid.} 
 \label{fig:grid0}
 \end{figure}

We note the distribution function at the time $t$ and at the position $x$ in the phase space: $\mathbf{g(t,x)}$. We discretize the time space with a constant time step $\Delta t$. Writing $t^n=n\Delta t$, we then note ${g}^{n}(x)=g(t^n,x)$ the value of the distribution function at $t^n$. Using the previous discretization, we define the distribution function at the cells faces as ${g}^{n}_{i+1/2}=g^n(x_{i+1/2})$ (\fref{ghalf0}).  
At last, we define the average of the distribution function on one cell $i$ at $t^n$ by: 
\eqs
\bar{g}^n_i = \frac{1}{\Delta x} \int^{x_{i+1/2}}_{x_{i-1/2}} g^n(x) dx,\  i=0\dots N
\eeqs

 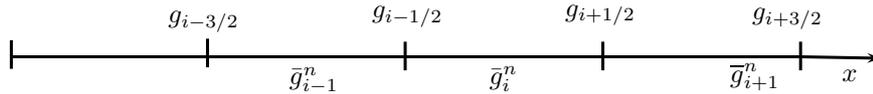
\begin{figure}[!h]\label{fig:ghalf0} 
 \centering
 \begin{pspicture}(0,-0.63390625)(11.821875,0.63390625)

 \rput(7.749219,0.45390624){$g_{i+1/2}$}

 \rput(5.209219,0.43390626){$g_{i-1/2}$}
 
 \rput(9.781406,-0.36609375){$\overline{g}^{n}_{i+1}$}

 \rput(2.5292187,0.37390625){$g_{i-3/2}$}

 \rput(10.209219,0.41390625){$g_{i+3/2}$}
 \psline[linewidth=0.04cm](0.0,-0.09609375)(10.38,-0.09609375)
 \psline[linewidth=0.04cm](10.38,0.10390625)(10.38,-0.27609375)
 \psline[linewidth=0.04cm](7.78,0.10390625)(7.78,-0.27609375)
 \psline[linewidth=0.04cm](5.18,0.10390625)(5.18,-0.27609375)
 \psline[linewidth=0.04cm](2.58,0.14390625)(2.58,-0.23609374)
 \psline[linewidth=0.04cm](0.0,0.12390625)(0.0,-0.25609374)
 \psline[linewidth=0.04cm,arrowsize=0.05291667cm 2.0,arrowlength=1.4,arrowinset=0.4]{->}(10.38,-0.09609375)(11.44,-0.11609375)
 
 \rput(11.026563,-0.34609374){$x$}

 \rput(6.4714065,-0.38609374){$\bar{g}^{n}_i$}

 \rput(3.9814062,-0.40609375){$\bar{g}^{n}_{i-1}$}
 \end{pspicture} 
 \caption{Normative example for the distribution g in the cells and at the nodes}
 \end{figure}

\subsection{Conservative semi-Lagrangian scheme principle}
The mass conservation in a lagrangian volume of the phase space between $t^n$ and $t^{n+1}$ reads as follows:  
\eq
\label{eq:eqvol}
\int^{}_{Vol^{n+1}}g(x,t^{n+1})dx=\int^{}_{Vol^{n}}g(x,t^{n})dx
\eeq
$$\text{with } Vol^{n}=\left\lbrace X(x^{n+1},t^{n}) |X(x^{n+1},t^{n+1}) \in {Vol^{n+1}}\right\rbrace $$
where $X(x,t^{n+1})$ describes the characteristic curve which passes by $x^{n+1}=X(x^{n+1},t^{n+1})$ at $t^{n+1}$.  Thus $X(x^{n+1},t^{n})$ is the point by which the trajectory passes at $t^n$ such as this trajectory also passes by $x^{n+1}$ at $t^{n+1}$. The characteristic curves are obtained by solving the following equation: 
\[
\deriv{X(x,t)}{t}=a(x,t)
\]
with a given initial condition. 
We consider hereafter the trajectories which pass by the cells faces $x_{i+1/2}$ at $t^{n+1}$ which are noted: $x_{i+1/2}=X(x^{n+1}_{i+1/2},t^{n+1})$ and we introduce $x^*_{i+1/2}$ which is the point on each characteristic curve at the time $t^n$: 
\[x^*_{i+1/2} =X(x^{n+1}_{i+1/2},t^{n}).\]
 
The conservation equation \eqref{eq:eqvol}  can thus be written using a 1D discretized form:
\Ea
\Delta x ~\bar{g}^{n+1}_i =  \int^{x_{i+1/2}}_{x_{i-1/2}} g^{n+1}(y) dy = \int^{x^*_{i+1/2}}_{x^*_{i-1/2}} g^{n}(y) dy
\label{eq:eqf}
\Eae
 where \[\begin{cases}
 x^*_{i+1/2}-x^*_{i-1/2}={Vol}^n \\
 x_{i+1/2}-x_{i-1/2}={Vol}^{n+1}=\Delta x 
 \end{cases}
 \]
In the conservative semi-Lagrangian formalism, the right hand side of equation \eqref{eq:eqf}  is numerically computed as follows:
\eq
\Delta x ~ \bar{g}^{n+1}_i =  \int^{x^*_{i+1/2}}_{x^*_{i-1/2}} g^{n}(y) dy=
G(x^*_{i+1/2})-G(x^*_{i+1/2})
\label{eq:semilag}
\eeq
 where $G(x)$ is the cumulative  or primitive function of $g$ defined as:
 \[ {G(x)}=\int^{x}_{x_{-1/2}} g(y)dy. 
\]
 This primitive function can be computed exactly at each cell face of the mesh:
 $$G(x_{i+1/2})=G(x_{-1/2})+ \displaystyle \sum_{k=0}^i  \Delta x ~\bar{g}^{n}_k.$$
 These values at faces  are then interpolated by basis functions to obtain an approximate reconstruction $G_h(x)$ of $G(x)$ for any $x$.\\
 For instance, the PSM scheme uses cubic splines and the LAG scheme uses third order lagrangian polynoma as interpolation functions to obtain the reconstructed function $G_h(x)$.  
 
\subsection{Finite volume form equivalence} 
The equation \eqref{eq:eqf} can be split in three terms:
 \Ea
\Delta x ~\bar{g}^{n+1}_i &=&  \underbrace{\int^{x_{i-1/2}}_{x^*_{i-1/2}} g^n (y) dy}_{\phi_{i-1/2}} +
\underbrace{\int^{x_{i+1/2}}_{x_{i-1/2}} g^n (y) dy}_{\bar{g}^{n}_i \Delta x} + \underbrace{\int^{x^*_{i-1/2}}_{x_{i+1/2}} g^n (y) dy}_{ -\phi_{i+1/2} }
\label{eq:finitevolume}
\Eae
We name $\phi_{i+1/2}$ the following quantity: 
\[
\phi_{i+1/2}=\int^{x_{i+1/2}}_{x^*_{i+1/2}} g^n(y) dy.
\]
We call it 'flux' since it represents the algebraic quantity which is carried through the node $x_{i+1/2}$, by identification with the finite volume formalism.

The equation \eqref{eq:finitevolume} is a finite volume equation (see \fref{flux}) meaning that the new value of the distribution function $\bar{g}^{n+1}_i$ in the cell $i$ is the sum of its value at time $t^n$ and the incoming or outgoing flux $\phi_{i\pm1/2}$.
\begin{equation}\label{fvform}
\bar{g}^{n+1}_i = \bar{g}^{n}_i -\left(\frac{\phi_{i+1/2}-\phi_{i-1/2}}{\Delta x} \right).  
\end{equation}

 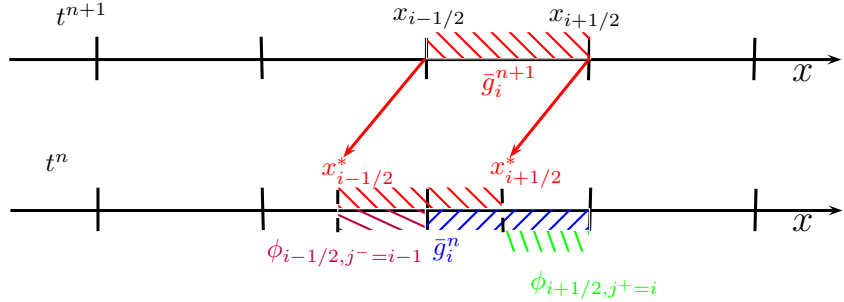
\begin{figure}[!h]
 \centering
  \psset{unit=.72cm}
  \begin{pspicture}(0,-3)(15,3)
 \psline[linewidth=0.04cm,arrowsize=0.05291667cm 2.0,arrowlength=1.4,arrowinset=0.4]{->}(0.0,1.3465625)(15.2,1.3465625)
 \rput(14.5,1.1){\Large $x$}
 \rput(14.5,-1.67){\Large $x$}
 \psline[linewidth=0.04cm](1.6,0.9665625)(1.6,1.7665625)
 \psline[linewidth=0.04cm](4.6,1.7665625)(4.62,0.9865625)
 \psline[linewidth=0.04cm](7.62,1.7465625)(7.62,0.9865625)
 \psline[linewidth=0.04cm](10.6,1.7665625)(10.58,0.9665625)
 \psline[linewidth=0.04cm](13.62,1.7465625)(13.6,0.9665625)
 \psline[linewidth=0.04cm,arrowsize=0.05291667cm 2.0,arrowlength=1.4,arrowinset=0.4]{->}(0.02,-1.4334375)(15.22,-1.4334375)
 \psline[linewidth=0.04cm](1.62,-1.8134375)(1.62,-1.0134375)
 \psline[linewidth=0.04cm](4.62,-1.0134375)(4.64,-1.7934375)
 \psline[linewidth=0.04cm](7.64,-1.0334375)(7.64,-1.7934375)
 \psline[linewidth=0.04cm](10.62,-1.0134375)(10.6,-1.8134375)
 \psline[linewidth=0.04cm](13.64,-1.0334375)(13.62,-1.8134375)
 \psframe[linewidth=0,linecolor=white,fillstyle=vlines,hatchangle=45,hatchcolor=red](7.62,1.36)(10.62,1.86)
 \rput(1.2453125,2.1865625){ $t^{n+1}$}
 \rput(0.8053125,-0.5534375){ $t^{n}$}
 \rput(7.681406,2.0965624){$x_{i-1/2}$}
 \rput(6.371406,-0.7434375){\color{red}$x^*_{i-1/2}$}
 \rput(10.521406,2.0765624){$x_{i+1/2}$}
 \psframe[linewidth=0,linecolor=white,fillstyle=vlines,hatchangle=45,hatchcolor=red](6,-1.4)
 (9,-1)%
 \psframe[linewidth=0,linecolor=white,fillstyle=vlines,hatchangle=65,hatchcolor=purple](6,-1.8)
 (7.6,-1.4)
 \psframe[linewidth=0,linecolor=white,fillstyle=vlines,hatchangle=-45,hatchcolor=blue](7.6,-1.8)
 (10.6,-1.4)
 \psframe[linewidth=0,linecolor=white,fillstyle=vlines,hatchangle=20,hatchcolor=green](9,-2.2)
 (10.6,-1.8)
 \rput(9.4514065,-0.7234375){\color{red}$x^*_{i+1/2}$}
 \psline[linewidth=0.04cm,linestyle=dashed,dash=0.16cm 0.16cm](6.0,-1.0534375)(6.0,-1.8334374)
 \psline[linewidth=0.04cm,linestyle=dashed,dash=0.16cm 0.16cm](9.02,-1.8334374)(9.0,-1.0334375)
 \psline[linewidth=0.04cm,arrowsize=0.05291667cm 2.0,arrowlength=1.4,arrowinset=0.4,linecolor=red]{->}(10.6,1.3665625)(9.14,-0.4534375)
 \psline[linewidth=0.04cm,arrowsize=0.05291667cm 2.0,arrowlength=1.4,arrowinset=0.4,linecolor=red]{->}(7.6,1.3665625)(6.1,-0.45634375)
 \rput(8,-2.1634376){\color{blue}$\bar{g}^{n}_i$}
 \rput(9.131406,0.9365625){\color{red}$\bar{g}^{n+1}_i$}
 \rput(10.691406,-2.8034375){\color{green}$\phi_{i+1/2,j^+=i}$}
 \rput(6.141406,-2.2034376){\color{purple}$\phi_{i-1/2,j^-=i-1}$}
 \end{pspicture}
 \caption{Conservative evolution of a finite volume scheme}
 \label{fig:flux}
 \end{figure}
As for the semi-Lagrangian formalism \eqref{eq:semilag}, the integrals to be computed to obtain the fluxes can be approximated by a reconstruction of the primitive of the distribution function $G_h$:
 \Ea
\phi_{i+1/2}=\int^{x_{i+1/2}}_{x^*_{i+1/2}} g^n(y) dy \approx G_h(x_{i+1/2})-G_h(x^*_{i+1/2}).
\label{eq:flux_prim}
\Eae

\subsection{Time scheme used in the Gysela code}
The time scheme has been modified, because the classical second order leap-frog algorithm used in the Gysela code is not robust enough to deal with the conservative semi-Lagrangian schemes considered here:
\begin{equation*}
\bar{g}^{n+1}_i = \bar{g}^{n-1}_i -\left(\frac{\phi^n_{i+1/2}-\phi^n_{i-1/2}}{\Delta x} \right)  
\end{equation*}
Moreover, this time scheme enforces a constant time step because it involves three different time steps, which is quite restrictive for high iteration numbers simulations. We thus turn the time scheme to a Predictor-Corrector (or Runge Kutta order 2) method, which allows to use a variable time step $\Delta t^n$ computed at time $t^n$. It is computed according to a CFL like condition necessary for the finite volume scheme stability which assesses that the maximum displacement in the domain is less than a fraction of the cells size (using a regular mesh):
\begin{equation} \label{cfl}
\Delta t^n =CFL \displaystyle \min_{d=1,D} \left( \dfrac{\Delta x_d}{\displaystyle \max_{x}  (a^n_d (x))} \right) 
\end{equation} 
with $\Delta x_d$ the space step in direction $d$ and $a^n_d$ the velocity at time $t^n$ in the space direction $d\in[1,D]$. $CFL$ is a coefficient which $0<CFL<1$.

\begin{rem}
 The finite volumes scheme form \eqref{fvform} and the semi-Lagrangian scheme  \eqref{eq:semilag} are strictly equivalent, since the displacement is restricted to $CFL \le 1$. The finite volumes scheme \eqref{fvform} is not defined for displacements bigger than one cell, with $CFL>1$, although the semi-Lagrangian scheme \eqref{eq:semilag}  could be written for any time step. However, the stability of both schemes with $CFL>1$ in a general situation is not demonstrated.  
\end{rem}

\paragraph{Predictor-Corrector Algorithm:}
\begin{itemize}
\item At beginning of the iteration at $t^n$, we compute  $\Delta t^n$ according to \eqref{cfl}.
\item Prediction step : we compute a order 1 in time approximation of the solution at time $t^{n+1/2}$ with half a time step according to values $\bar{g}^{n}_i$ and cell faces fluxes $\phi^n_{d,i+1/2}$ in all direction $d\in [1,D]$ at time $t^n$:
\eqs
\bar{g}^{n+1/2}_i = \bar{g}^{n}_i - \displaystyle \sum_{d=1}^D \left(\frac{\phi^n_{d,i+1/2}-\phi^n_{d,i-1/2}}{\Delta x_d} \right)  
\eeqs
We compute the electric potential at same time $t^{n+1/2}$:
\begin{equation*}
\begin{array}{ll} 
-\dfrac{1}{B \omega_{i} } \nabla_\perp \cdot (n_0 \nabla_\perp \Phi^{n+1/2} ) +\dfrac{e}{\kappa T_e}(\Phi^{n+1/2} -<\Phi^{n+1/2} >_{\phi})= \int  \bar{g}^{n+1/2} d\mu dv_{\parallel}-n_0 
\end{array}
\end{equation*}

\item Correction step : we compute an order 2 in time approximation of the solution at time $t^{n+1}$ according to values $\bar{g}^{n}_i$ at time $t^n$ and cell faces fluxes $\phi^{n+1/2}_{d,i+1/2}$ in all direction $d\in [1,D]$ at time $t^{n+1/2}$:
\eqs
\bar{g}^{n+1}_i = \bar{g}^{n}_i - \displaystyle \sum_{d=1}^D \left(\frac{\phi^{n+1/2}_{d,i+1/2}-\phi^{n+1/2}_{d,i-1/2}}{\Delta x_d} \right)  
\eeqs
We compute the electric potential at same time $t^{n+1}$:
\begin{equation*}
\begin{array}{ll} 
-\dfrac{1}{B \omega_{i} } \nabla_\perp \cdot (n_0 \nabla_\perp \Phi^{n+1} ) +\dfrac{e}{\kappa T_e}(\Phi^{n+1} -<\Phi^{n+1} >_{\phi})= \int  \bar{g}^{n+1} d\mu dv_{\parallel}-n_0 
\end{array}
\end{equation*}

\end{itemize}

\section{The Hermite formalism applied to PSM and LAG} 

Two schemes can be used to compute the new distribution function (reconstructing the flux or the primitive function):
\begin{itemize}
\item The Parabolic Splines Method called PSM, which uses cubic splines functions for the interpolations. 
\item The Lagrangian method (LAG), which uses third order Lagrange polynoms for the interpolations. 
\end{itemize} 

The PSM scheme \cite{braeunig} or the LAG scheme \cite{crouseilles} are conservative semi-Lagrangian schemes that only differs by the interpolation functions used for the reconstruction step. We aim to study some limiters for both schemes and using the Hermite formalism. The Hermite formalism is a generic formulation for the interpolation polynoms. Indeed,  the scheme LAG or PSM in this formalism are set only by the way of computing the distribution function  at the faces $g_{i+1/2}$. Thus we can easily use them simultaneously in a code. 

We rewrite hereafter this generic formulation, based on the Hermite formalism for the conservative schemes PSM and LAG, proposed in \cite{crouseilles}. First, we give the expression of the flux $\phi_{i+1/2}$, then we propose an application to the PSM and LAG schemes.

\subsection{Computation of the flux $\phi_{i+1/2}$ with the Hermite formalism}
\label{sec:hermiteflux}
Following \cite{PSM1} and \cite{crouseilles}, we reconstruct the distribution function $g$ with a second order polynom $P_k$, which interpolates the distribution function in the cell $k$. We first assume that the reconstructed function is continuous at the cell faces. We name the value of $g$ at the cell faces $g_{k+1/2}$. At the end of the section, we will give the general formula for a discontinuous function at the cell faces.
The distribution function is approximated by a second order polynom in cell $k$, which corresponds to interpolate the primitive of the distribution function with a third order polynom:
\[
\forall k, P_k(z)=c^{(k)}+b^{(k)} z+a^{(k)} z^2\text{ with }z\in[0,\Delta x]
\]      
The distribution function $g$ is continuous at the faces. That implies:
\[\begin{cases}
P_k(0)=g_{k-1/2}
\\P_k(\Delta x)=g_{k+1/2}
\end{cases}\]

A third condition on the polynom comes from the mass conservation  in the cell $k$, that reads:
$$ \frac{1}{\Delta x} \int^{x_{k+1/2}}_{x_{k-1/2}} P_k(y) dy=\bar{g}_k$$
Therefore, the polynomial coefficients can be written:
\begin{equation}\label{polcoef}
\begin{array}{l}
a^{(k)} = \left(3g_{k-1/2}+3g_{k+1/2}-6\bar{g}_k\right)/{\Delta x}^2\\
b^{(k)} = \left(-4g_{k-1/2}-2g_{k+1/2}+6\bar{g}_k\right)/{\Delta x}\\
c^{(k)} = g_{k-1/2}
\end{array}
\end{equation}

Finally, $\forall x\in[x_{k-1/2},x_{k+1/2}]$, $P_k(x)$ is an approximation of the distribution function $g$ in the cell $k$. This approximation depends on $g_{k-1/2}$, $g_{k+1/2}$, and $\bar{g}_k$. 
\newline
Using this interpolation, we are able to evaluate the flux $\phi_{k+1/2}$: 
\eq
 \phi_{k+1/2}=\int^{x_{k+1/2}}_{x^*_{k+1/2}} g^n(y) dy. 
\label{eq:fluxtransport}
\eeq

We take the additional assumption that the displacement at the cell face $k+1/2$,  i.e. $\alpha_{k+1/2}={x_{k+1/2}-x^*_{k+1/2}}$ satisfies: \[
|\alpha_{k+1/2}|\leq\Delta x
\] 
It is necessary to know in which cell $j$ ($j=k$ or $j=k+1$ since the displacement $|\alpha_{k+1/2}|\leq\Delta x$) is located the foot $x^*_{k+1/2}$ of the characteristic that passes by $x_{k+1/2}$ at time $t^{n+1}$. 
We can thus give an approximation of the fluxes \eqref{eq:fluxtransport} towards the face $j+1/2$ based on the set of polynoms $P_k$:
\Eas
&&\exists j\:|\: x^*_{k+1/2}\in[x_{j-1/2},x_{j+1/2}]\\
&&\phi_{k+1/2}=\int^{x_{k+1/2}-x_{j-1/2}}_{x^*_{k+1/2}-x_{j-1/2}} g^n(Y) dY \approx  \int^{x_{k+1/2}-x_{j-1/2}}_{x^*_{k+1/2}-x_{j-1/2}} P_j (Y) dY
\Ease
with the change of variable $Y=y-x_{j-1/2} \in [0,\Delta x]$.


We exhibit hereafter the flux $\phi_{k+1/2}$ function of the $g_{k-1/2}$, $g_{k+1/2}$, and $\bar{g}_k$. 
We note $j_{k+1/2}$ the cell where is located the foot of the characteristic passing by the face $k+1/2$: 
\[
x^*_{k+1/2}\in[x_{j_{k+1/2}-1/2},x_{j_{k+1/2}+1/2}]
\]
The index $j_{k+1/2}$ indicates from which cell the flow is coming. For instance, a negative displacement $\alpha_{i+1/2}<0$ on the face $k+1/2$ means that the flow comes from the cell $k+1$, so that  $j_{k+1/2}=k+1$.

So we have: 
\begin{itemize}
\item If  $ x^*_{k+1/2}<x_{k+1/2}$ then 
\begin{itemize}
\item $j_{k+1/2}=k$ 
\item $\phi_{k+1/2}=\int^{x_{k+1/2}-x_{k-1/2}=\Delta x}_{x^*_{k+1/2}-x_{k+1/2}} P_{k}(Y) dY$
\end{itemize}
\item If  $ x^*_{k+1/2}>x_{k+1/2}$ then 
\begin{itemize}
\item $j_{k+1/2}=k+1$
\item $\phi_{k+1/2}=\int^{x_{k+1/2}-x_{k+1/2}=0}_{x^*_{k+1/2}-x_{k+1/2}} P_{k+1}(Y) dY$
\end{itemize}  
\end{itemize}
We introduce \[\delta=\frac{x_{k+1/2}-x_{j_{k+1/2}-1/2}}{\Delta x},\]  thus $\delta=0$ or  $\delta=1$ indicates the upwinding direction.\\
By introducing a normalized displacement 
$$\beta=\frac{x^*_{k+1/2}-x_{k+1/2}}{\Delta x},$$ 
we can write the flux:
\[
 \phi_{k+1/2}(\beta)=\int^{\delta \Delta x}_{(\delta-\beta) \Delta x} P_{j_{k+1/2}}(Y) dY
\] 
That leads to:
\Ea
\frac{ \phi_{k+1/2}}{\Delta x}=a_{{j_{k+1/2}}} \beta +b_{{j_{k+1/2}}} \left(-\beta^2+2\beta\delta \right)+c_{{j_{k+1/2}}} \left(\beta^3+3\delta^2 \beta -3\beta^2\delta\right)
\label{eq:eq5}\Eae
with
\[\begin{cases}     
 a_{j_{k+1/2}}=g_{j_{k+1/2}-1/2}\\    b_{{j_{k+1/2}}}=\frac{-4g_{{j_{k+1/2}}-1/2}+2g_{{j_{k+1/2}}+1/2}+6\bar{g}_{{j_{k+1/2}}}}{2}\\ c_{{j_{k+1/2}}}=\frac{3g_{{j_{k+1/2}}-1/2}+3g_{{j_{k+1/2}}+1/2}-6\bar{g}_{{j_{k+1/2}}}}{3}.
\end{cases}\]

By ordering differently the polynomial expression \eqref{eq:eq5} and replacing the coefficients by their values, we try to get an equivalent formulation to the flux expressions proposed by \cite{crouseilles}.
 
 \Eas
\phi_{k+1/2,{j_{k+1/2}}}(\beta)&=&\Delta x[g_{{j_{k+1/2}}-1/2} \left( \beta(1-\delta)+\beta^2(2-3\delta)+\beta^3 \right)\\&+&g_{{j_{k+1/2}}+1/2}\left( \beta \delta +\beta^2(1-3\delta)+\beta^3\right)\\&+&\bar{g}_{{j_{k+1/2}}}\left(\beta^2(-3+6\delta)+\beta^3 (-2)\right)]
\Ease

with 
\Eas
\beta= \frac{x_{k+1/2}-x^{*}_{k+1/2}}{\Delta x} \in [-1,1]\\
\delta=\left\{
    \begin{array}{ll}
        0 & \mbox{if } x_{k+1/2}<x^{*}_{k+1/2} \\
        1 & \mbox{if } x_{k+1/2}>x^{*}_{k+1/2}.
    \end{array}
\right.
\Ease
It's obvious that $\beta$ depends on k hence formally $\beta=\beta_{k+1/2}$

\paragraph{Positive displacement $\mathbf{\alpha>0}$ i.e. $\mathbf{\delta=1}$\\}
We have $j_{k+1/2}=k$, hence 
\Eas
\frac{\phi_{k+1/2}(\beta)}{\Delta x} &=& g_{{k}-1/2} \left(\beta^2(\beta-1) \right)\nonumber\\
 &+& g_{{k}+1/2}\left( \beta(1-\beta)^2\right)\nonumber\\
 &+& \bar{g}_{{k}}\left(\beta^2(3-2\beta)\right)
\Ease

\paragraph{Negative displacement $\mathbf{\alpha<0}$ i.e. $\mathbf{\delta=0}$\\}

We have $j_{k+1/2}=k+1$, hence 
\Eas
\frac{ \phi_{k+1/2}(\beta)}{\Delta x}&=&g_{k+1/2} \left(\beta(\beta+1)^2 \right)\\
&+&g_{k+3/2}\left( \beta^2(1+\beta)\right)\\
&+&\bar{g}_{k+3/2}\left(\beta^2(-3-2\beta)\right)
\Ease

Thus, we have established a generic formulation for the flux considering displacements smaller than one cell. We have assumed that the distribution function is continuous at the cell face. This generic expression is equivalent to the formulation proposed by \cite{crouseilles}.

\paragraph{Hermite formalism for a discontinuous reconstruction\\}

Performing the same calculations, we can extend the previous formalism to a distribution function which is not continuous at the faces.
We name $g^{+}_{k+1/2}$ and $g^{-}_{k+1/2}$ respectively the left and right values of the distribution function at the cell face $k+1/2$ (see \fref{leftright}).
In general that means:
\[
g^{+}_{k+1/2}\neq g^{-}_{k+1/2} \\
\] 
  \begin{figure}[!h]
 \centering
 \scalebox{1}  
 {
 \begin{pspicture}(0,-0.758125)(8.622812,0.758125)
 \psline[linewidth=0.04cm](2.3209374,-0.1003125)(5.5209374,-0.1203125)
 \psline[linewidth=0.04cm](4.9209375,0.0996875)(4.9209375,-0.3003125)
 \psline[linewidth=0.04cm](3.1209376,0.0996875)(3.1209376,-0.3003125)

 \rput(4.0423436,-0.5303125){$g^{n}_j$}

 \rput(1.2023437,0.4696875){$g^{+}_{j-1/2}$}

 \rput(3.1723437,0.5496875){$g^{-}_{j-1/2}$}

 \rput(5.0823436,0.5296875){$g^{+}_{j+1/2}$}

 \rput(7.2923436,0.5696875){$g^{-}_{j+1/2}$}
 \psline[linewidth=0.04cm,arrowsize=0.05291667cm 2.0,arrowlength=1.4,arrowinset=0.4]{->}(1.1209375,0.2196875)(2.9409375,0.0)
 \psline[linewidth=0.04cm,arrowsize=0.05291667cm 2.0,arrowlength=1.4,arrowinset=0.4]{->}(3.1609375,0.3196875)(3.2809374,0.0196875)
 \psline[linewidth=0.04cm,arrowsize=0.05291667cm 2.0,arrowlength=1.4,arrowinset=0.4]{->}(4.7809377,0.2796875)(4.7609377,-0.0203125)
 \psline[linewidth=0.04cm,arrowsize=0.05291667cm 2.0,arrowlength=1.4,arrowinset=0.4]{->}(7.2009373,0.3596875)(5.1009374,-0.0403125)
 \end{pspicture} 
 }
 \caption{Asymmetric node}
 \label{fig:leftright}
 \end{figure}
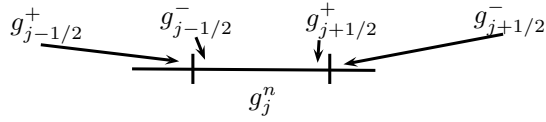  
The conditions on the polynoms at the cell faces are thus changed to:
\[
\begin{cases}
P_k(0)=g^{-}_{k-1/2}\\
P_k(\Delta x) =g^{+}_{k+1/2}
\end{cases}
\]
That leads to the final formulation: 
 \Ea
\phi_{k+1/2,{j_{k+1/2}}}(\beta)&=&\Delta x[g^{-}_{{j_{k+1/2}}-1/2} \left( \beta(1-\delta)+\beta^2(2-3\delta)+\beta^3 \right)\nonumber\\
&+&g^{+}_{{j_{k+1/2}}+1/2}\left( \beta \delta +\beta^2(1-3\delta)+\beta^3\right)\nonumber\\
&+&\bar{g}_{{j_{k+1/2}}}\left(\beta^2(-3+6\delta)+\beta^3 (-2)\right)]
\label{eq:eqflux}
\Eae
where
\Eas
&&\beta=\frac{\alpha_{i+1/2}}{\Delta x}= \frac{x_{k+1/2}-x^{*}_{k+1/2}}{\Delta x} \in [-1,1]\\
&&\delta=\left\{
    \begin{array}{ll}
        0 & \mbox{if } x_{k+1/2}<x^{*}_{k+1/2} \\
        1 & \mbox{if } x_{k+1/2}>x^{*}_{k+1/2}
    \end{array}
\right.
\Ease

\paragraph{Positive displacement $\mathbf{\alpha>0}$ i.e. $\mathbf{\delta=1}$\\}

We have $j_{k+1/2}=k$, hence 
\Eas
\frac{\phi_{k+1/2}(\beta)}{\Delta x} &=& g^{-}_{{k}-1/2} \left(\beta^2(\beta-1) \right)\\ &+& g^{+}_{{k}+1/2}\left( \beta(1-\beta)^2\right)\\ 
&+& \bar{g}_{{k}}\left(\beta^2(3-2\beta)\right)
\Ease

\paragraph{Negative displacement $\mathbf{\alpha<0}$ i.e. $\mathbf{\delta=0}$\\}

We have $j_{k+1/2}=k+1$, hence
\Eas
\frac{ \phi_{k+1/2}(\beta)}{\Delta x}&=&g^{-}_{k+1/2} \left(\beta(\beta+1)^2 \right)\\
&+&g^{+}_{k+3/2}\left( \beta^2(1+\beta)\right)\\&+&\bar{g}_{k+1}\left(\beta^2(-3-2\beta)\right)
\Ease

\subsection{PSM scheme with the Hermite formalism}
\subsubsection{Hermite formalism using splines interpolation}
The type of scheme expressed in the hermite formalism (PSM or LAG) only depends on the manner the face values $g^{+}_{k+1/2}$ and $g^{-}_{k+1/2}$ are computed. 
Here the PSM scheme has the property that the reconstructed function has continuous derivatives of the distribution function $g(x,t)$ at cell faces \cite{PSM1}. Remembering that the distribution function  $g$ is reconstructed by a second order polynom  $g(x)\approx P_i(x)$ in the cell $i$, the continuity of the derivative at face $i+1/2$ can be written:
 \[
\deriv{P_i(Y)}{Y}\mid_{Y=\Delta x}=\deriv{P_{i+1(Y)}}{Y}\mid_{Y=0}
\] 
Using the polynom coefficients expression \eqref{polcoef} we obtain: 
\Eas
\deriv{P_i(Y)}{Y}\mid_{Y=\Delta x}&=&\deriv{P_{i+1}(Y)}{Y}\mid_{Y=0}\\ 
\Leftrightarrow  2 a_{i} \Delta x + b_i &=& b_{i+1}\\
 \Leftrightarrow 
g_{i-1/2}+4g_{i+1/2}+g_{i+3/2}&=&3(\bar{g}_i+\bar{g}_{i+1})
\Ease

This PSM formulation regardless of the boundary conditions is equivalent to the semi-Lagrangian PSM formalism, used in \cite{braeunig} for instance. A rigorous proof of the equivalence between the two formulations is furnished by \cite{crouseilles}, except for the boundary conditions.
    
\subsubsection{Periodic boundary conditions}
\paragraph*{}
In this section, we present the way of imposing the boundary conditions. Extending boundary conditions for the PSM scheme to the Hermite formalism is not simple, especially to get a complete equivalence between the Hermite formalism and the semi-Lagrangian formalism \eqref{eq:semilag} \cite{braeunig}.  
\newline
We named $G$ the primitive function which is defined as:
\eq
 G(x)=\int^{x}_{x_{-1/2}} g(x)dx 
\eeq

We also define the mesh fitted to a periodic domain(\fref{periodic}):
 \begin{figure}[!h]
 {
 \psset{unit=.62cm}
 \centering
 \begin{pspicture}(0,-1.23125)(0.062187,1.23125)
 \psline[linewidth=0.04cm,arrowsize=0.05291667cm 2.0,arrowlength=1.4,arrowinset=0.4]{->}(1.0021875,0.485625)(17.042187,0.445625)
 \psline[linewidth=0.04cm](1.0021875,0.685625)(1.0021875,0.285625)
 \psline[linewidth=0.04cm](3.4021876,0.685625)(3.4021876,0.285625)
 \psline[linewidth=0.04cm](5.7821875,0.665625)(5.7821875,0.265625)
 \psline[linewidth=0.04cm](8.182187,0.665625)(8.182187,0.265625)
 \psline[linewidth=0.04cm](10.602187,0.685625)(10.602187,0.285625)
 \psline[linewidth=0.04cm](13.002188,0.685625)(13.002188,0.285625)
 \psline[linewidth=0.04cm](15.402187,0.685625)(15.402187,0.285625)

 \rput(4.476406,0.100625){ 0}

 \rput(6.8795314,0.080625){ 1}
 \rput(14.477344,0.100625){ N+1=0}
 \rput(11.600938,0.040625){ N}
 \rput(3.4142187,0.980625){ $-1/2=N+1/2$}
 \rput(5.787344,1.020625){ 1/2}

 \rput(8.184218,0.980625){ 3/2}

 \rput(10.5260935,0.960625){ N-1/2}

 \rput(12.956094,0.960625){ $N+1/2=-1/2$}
 \psline[linewidth=0.04cm](14.202188,-0.714375)(4.4021873,-0.694375)
 \psline[linewidth=0.04cm,arrowsize=0.05291667cm 2.0,arrowlength=1.4,arrowinset=0.4]{->}(14.182187,-0.714375)(14.202188,-0.294375)
 \psline[linewidth=0.04cm,arrowsize=0.05291667cm 2.0,arrowlength=1.4,arrowinset=0.4]{->}(4.4221873,-0.694375)(4.4221873,-0.334375)

 \rput(4.4735937,-0.924375){$0$}

 \rput(14.153594,-1.004375){$2\pi$}
 \end{pspicture} 
 }
 \caption{Mesh for a periodic domain}
 \label{fig:periodic}
 \end{figure}
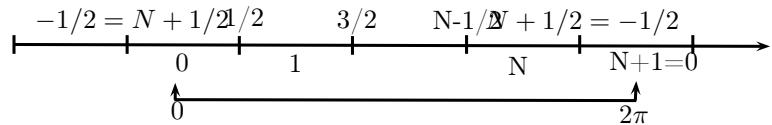
\paragraph{}
The periodic boundary conditions for the semi-Lagrangian scheme are approximated by conditions on the primitive function derivatives:
\eqs \begin{cases}
G'(x_{-1/2})=G'(x_{N+1/2})\\
G''(x_{-1/2})=G''(x_{N+1/2})
\end{cases}
\eeqs
which is equivalent to set continuity of the distribution function and its first derivative:
\eqs \begin{cases}
g(x_{-1/2})=g(x_{N+1/2})\\
g'(x_{-1/2})=g'(x_{N+1/2})
\end{cases}
\eeqs
The periodic boundary conditions for the Hermite formalism are then obtained by setting the same constraint on the polynomial reconstruction:
\Eas
&& \begin{cases}
P_0(0)=P_N(\Delta x)\\
\deriv{P_0(Y)}{Y}|_{Y=0}=\deriv{P_N(Y)}{Y}|_{Y=\Delta x}\\
\end{cases}
\\
\Leftrightarrow
 && \begin{cases}
  g_{-1/2}=g_{N+1/2}\\
-4 g_{-1/2} -2 g_{1/2}+6 \bar{g}_{0}=2 g_{N-1/2} +4 g_{N+1/2}-6 \bar{g}_{N}
\end{cases}
\Ease
The complete linear system to solve to obtain values at nodes is then the following 
 \Eas
 \begin{cases}
g_{i-1/2}+4g_{i+1/2}+g_{i+3/2}=3(\bar{g}_i+\bar{g}_{i+1}),i=0\cdots N-1\\
g_{-1/2}=g_{N+1/2}\\
-4 g_{-1/2} -2 g_{1/2}+6 \bar{g}_{0}=2 g_{N-1/2} +4 g_{N+1/2}-6 \bar{g}_{N}\\
\end{cases}
\\
\Rightarrow
 \begin{cases}
 g_{i-1/2}+4g_{i+1/2}+g_{i+3/2}=3(\bar{g}_i+\bar{g}_{i+1}),i=0\cdots N-2\\
 g_{N-3/2}+4 g_{N-1/2} + g_{-1/2}= 3(\bar{g}_{N-1}+\bar{g}_{N}) \\
  g_{N-1/2}+4 g_{-1/2} + g_{1/2}= 3(\bar{g}_{0}+\bar{g}_{1})\\
\end{cases}\Ease
The matricial system $[A]X=B$ of dimension $N+1$ is finally:
\[
\underbrace{\begin{pmatrix}
4    & 1      & 0      & 0      & \cdots & 1      \\
1     & 4      & 1      & 0      &   \cdots     & \vdots \\
\vdots     & \ddots & \ddots & \ddots & \ddots & \vdots \\ 
\vdots&    \cdots    & 0      & 1      & 4      & 1      \\
1     & \cdots      & 0      &0       & 1      & 4
\end{pmatrix}}_{[A]}
\underbrace{\begin{pmatrix}
g_{-1/2}\\
g_{1/2}\\ 
\vdots\\
\vdots\\
g_{N-3/2}\\
g_{N-1/2}\\
\end{pmatrix}}_{X}
= 
\underbrace{\begin{pmatrix}
3(\bar{g}_{N}+\bar{g}_{0})\\
3(\bar{g}_{0}+\bar{g}_{1})\\ 
\vdots\\
\vdots\\
3(\bar{g}_{N-2}+\bar{g}_{N-1})\\
3(\bar{g}_{N-1}+\bar{g}_{N})\\
\end{pmatrix}}_{B}
\]

\subsubsection{Non-Periodic boundary conditions: natural conditions}

The natural boundary condition for the PSM scheme is approximated by  vanishing the second derivative of the primitive function, what corresponds to annulate the first derivative of the distribution function $g$, at the boundaries:
\[ \begin{cases}
G''(x_{-1/2})=0\\
G''(x_{N+1/2})=0
\end{cases}
\]
The corresponding natural boundary conditions for the Hermite formalism are:
\[ \begin{cases}
\deriv{P_0(Y)}{Y}|_{Y=0}=0\\
\deriv{P_N(Y)}{Y}|_{Y=\Delta x}=0\\
  \end{cases}\Leftrightarrow 
 \begin{cases}
 -4 g_{-1/2} -2 g_{1/2}+6 \bar{g}_{0}=0\\
 2 g_{N-1/2} +4 g_{N+1/2}-6 \bar{g}_{N}=0
\end{cases}
\]
 The system to solve to obtain the values of the function at the nodes is then: 
 \[
 \begin{cases}
 g_{i-1/2}+4g_{i+1/2}+g_{i+3/2}=3(\bar{g}_i+\bar{g}_{i+1}),i=0\cdots N-1\\
 4 g_{-1/2} +2 g_{1/2}=6 \bar{g}_{0}\\
 2 g_{N-1/2} +4 g_{N+1/2}=6 \bar{g}_{N}
\end{cases}
\]
Considering $g_{-1/2}=\bar{g}_0$ and $g_{N+1/2}=\bar{g}_N$, we solve from the following linear system of dimension $N+2$: 
\[
\underbrace{\begin{pmatrix}
4     & 2      & 0      & 0      & \cdots & 0      \\
1     & 4      & 1      & 0      &   \cdots     & \vdots \\
\vdots     & \ddots & \ddots & \ddots & \ddots & \vdots \\ 
\vdots&    \cdots    & 0      & 1      & 4      & 1      \\
0     & \cdots      & 0      &0       & 2      & 4
\end{pmatrix}}_{[A]}
\underbrace{\begin{pmatrix}
g_{-1/2}\\
g_{1/2}\\ 
\vdots\\
\vdots\\
g_{N-1/2}\\
g_{N+1/2}\\
\end{pmatrix}}_{[X]}
= 
\underbrace{\begin{pmatrix}
6(\bar{g}_{0})\\
3(\bar{g}_{0}+\bar{g}_{1})\\ 
\vdots\\
\vdots\\
3(\bar{g}_{N-1}+\bar{g}_{N})\\
6 (\bar{g}_{N})\\
\end{pmatrix}}_{[B]}
\]

\subsection{LAG scheme with the Hermite formalism} 
\subsubsection{Hermite formalism using Lagrange interpolation}
We outline hereafter the Hermite formulation for the LAG scheme. In order to set the distribution function values at the faces $g^{+}_{k+1/2}$ and $g^{-}_{k+1/2}$, we  start from the Lagrange interpolation of the flux and then, we make some calculations to show the equivalence with the Hermite formulation.
\newline
We build the primitive or cumulative function ${G}$ at the characteristic feet:
\Eas
g^{n+1}_i=G(x^*_{i+1/2})-G(x^*_{i-1/2})\\
G(x^*_{i+1/2})=\int^{x^*_{i+1/2}}_{x_{-1/2}} g^{n}(x)dx.
\Ease

The value of $G(x^*_{i+1/2})$ is obtained by an interpolation with third-order Lagrange polynoms of values at cell faces $G(x_{j+1/2})$:
\[G(x)=\sum^{i+2}_{j=i-1} G_{j-1/2}L_j(x), x\in[x_{i-1/2},x_{i+1/2}]\]
where $L_j$ are the Lagrange polynoms defined as:
\[L_j(x)=\prod_{k=i-1,k\neq j}^{i+2}\frac{x-x_{k-1/2}}{x_{j-1/2}-x_{k-1/2}}\]
and $G_{j-1/2}$ is the value of $G$ at the cell face $j-1/2$.\newline
Some calculations (given in annex \sref{annexeA}) lead to the expressions \eqref{eq:G1} and \eqref{eq:G2} for $G(x^{*}_{i+1/2})$:
\begin{itemize}
\item if $\beta_{i+1/2}>0$,
\Ea
G(x^{*}_{i+1/2})&=&\beta^3(1/6 G_{i-3/2}-1/2 G_{i-1/2}+1/2G_{i+1/2}-1/6 G_{i+3/2}) \nonumber
\\ &+&\beta^2 (1/2 G_{i-1/2}-G_{i+1/2}+1/2 G_{i+3/2}) \nonumber
\\&+&\beta(-1/6G_{i-3/2}+G_{i-1/2}-1/2 G_{i+1/2}-1/3G_{i+3/2})\nonumber
\\&+&G_{i+1/2}
\label{eq:G1}
\Eae
\item if $\beta_{i+1/2}<0$,
\Ea
G(x^{*}_{i+1/2})&=&\beta^3(1/6 G_{i-1/2}-1/2 G_{i+1/2}+1/2G_{i+3/2}-1/6 G_{i+5/2}) \nonumber
\\ &+&\beta^2 (1/2 G_{i-1/2}-G_{i+1/2}+1/2 G_{i+3/2} ) \nonumber
\\&+&\beta(1/3 G_{i-1/2}+1/2G_{i+1/2}-G_{i+3/2}+1/6G_{i+5/2})
\label{eq:G2}
\Eae
\end{itemize}
with $G_{k+1/2}=G(x_{k+1/2})$ and the function $\beta_{i+1/2}$:
\[
\Delta x ~\beta_{i+1/2} =x_{i+1/2}-x 
\]
We omit the subscript for $\beta$ when the context makes the index obvious. \\
The values of the primitive function at nodes $G_{k+1/2}$ are related to the cell centred values of the distribution function $g^{n}_k$ by:
 \[g^{n}_k\Delta x=G_{k+1/2}-G_{k-1/2}.\]
Thus, by replacing the values $G_{k+1/2}$ function of $g^{n}_{k}$ in such a way only the primitive function value   $G_{i+1/2}$ at cell face $i+1/2$ remains at the right hand side of relations \eqref{eq:G1} and \eqref{eq:G2}, we obtain:
\begin{itemize}
\item if $\beta_{i+1/2}>0$,
\Eas
G(x^{*}_{i+1/2})/\Delta x &=&\beta^3(-1/6g^{n}_{i-1}+1/3g^{n}_{i}-1/6g^{n}_{i+1} )
\\ &+&\beta^2 (-1/2g^{n}_{i}+1/2g^{n}_{i+1})
\\&+&\beta(1/6g^{n}_{i-1}-5/6g^{n}_{i}-1/3g^{n}_{i+1})
\\&+&G_{i+1/2}/\Delta x
\Ease
\item if $\beta_{i+1/2}<0$
\Eas
G(x^{*}_{i+1/2})/\Delta x&=&\beta^3(-1/6g^{n}_{i} +2/6g^{n}_{i+1}-1/6g^{n}_{i+2})
\\ &+&\beta^2 (-1/2g^{n}_{i}+1/2g^{n}_{i+1} )
\\&+&\beta(-1/3g^{n}_{i}-5/6g^{n}_{i+1}+1/6g^{n}_{i+2})
\\&+& G_{i+1/2}/\Delta x
\Ease
\end{itemize}

Let us recall the relation \eqref{eq:flux_prim} between the flux at cell face $i+1/2$ and the primitive function:
 \Ean
\phi_{i+1/2}= G(x_{i+1/2})-G(x^*_{i+1/2}).
\Eane
Therefore, we have an expression of the flux $\phi_{i+1/2}(\beta)$ function of cell centred values $(\bar{g}_k)_k$:
\begin{itemize}
\item if $\beta_{i+1/2}>0$,
\Eas
  \frac{\phi_{i+1/2}(\beta)}{\Delta x}&=&\beta^3( 1/6 g^{n}_{i-1}-1/3 g^{n}_{i}+1/6 g^{n}_{i+1} )
\\ &+&\beta^2 (1/2g^{n}_{i}-1/2g^{n}_{i+1})
\\&+&\beta(-1/6g^{n}_{i-1}+5/6g^{n}_{i}+1/3g^{n}_{i+1})
\Ease
\item if $\beta_{i+1/2}<0$
\Eas
  \frac{\phi_{i+1/2}(\beta)}{\Delta x}&=&\beta^3(1/6g^{n}_{i} -2/6g^{n}_{i+1}+1/6g^{n}_{i+2})
\\ &+&\beta^2 ( 1/2g^{n}_{i}-1/2g^{n}_{i+1} )
\\&+&\beta(1/3g^{n}_{i}+5/6g^{n}_{i+1}-1/6g^{n}_{i+2})
\Ease
\end{itemize}
Let us remember the generic Hermite formulation for the flux \eqref{eq:eqflux} function of values at cell faces $g^\pm_{k+1/2}$: 
\begin{itemize}
\item if $\beta_{i+1/2}>0$ then
\Eas
 \frac{\phi_{i+1/2}(\beta)}{\Delta x}
&=&\beta^3(g^{+}_{{i}-1/2}+g^{-}_{{i}+1/2}-2 g^n_{{i}})
\\&+& \beta^2(-g^{+}_{{i}-1/2}-2g^{-}_{{i}+1/2}+3 g^n_{{i}})
\\&+& \beta(g^{-}_{i+1/2})
\Ease
\item if $\beta_{i+1/2}<0$ then
\Eas
\frac{\phi_{i+1/2}(\beta)}{\Delta x} &=& \beta^3(g^{-}_{i+3/2}+g^{+}_{i+1/2}-2{g}^{n}_{i+1})
\\ &+& \beta^2(2g^{+}_{i+1/2}+g^{-}_{i+3/2}-3g^{n}_{i+1})
\\ &+& \beta (g^{+}_{i+1/2})
\Ease
\end{itemize}
By identifying in the two last expressions of the flux $\phi_{i+1/2}(\beta)$ the coeficients of each degree of the polynom in variable $\beta$, we obtain a necessary expression for values $g^{\pm}_{{k}-1/2}$ at cell faces: 
\eq
\begin{cases}
g^{+}_{{k}-1/2}=1/3g^{n}_{k-1}+5/6 g^{n}_{k}-1/6g^{n}_{k+1}  \\ 
g^{-}_{{k}+1/2}=-1/6g^{n}_{k-1}+5/6 g^{n}_{k}+1/3g^{n}_{k+1}.
\label{eq:finaleq}
\end{cases}
\eeq

As a conclusion, we have explicit values of $g^{+}_{{k}-1/2}$ and $g^{-}_{{k}+1/2}$ function of cell centred values of the distribution function at time $t^n$ with \eqref{eq:finaleq}, so we can directly  compute the generic Hermite formulation for the flux \eqref{eq:eqflux}  for the LAG scheme.

\subsubsection{Boundary conditions}
The boundary conditions for the LAG scheme are more simple than the PSM ones. Indeed, the coefficient are explicitly depending on the value of the distribution function $g^{n}_{k}$. So we can explicitly set values  at boundaries.
\paragraph{Non-periodic boundary conditions: natural conditions\\}

\[
\begin{cases}
g^{+}_{N+3/2}=1/3g^{n}_{N+1}+5/6 g^{n}_{N+1}-1/6g^{n}_{N+1}  \\ 
g^{-}_{{N}+3/2}=-1/6g^{n}_{N}+5/6 g^{n}_{N+1}+1/3g^{n}_{N+1}\\
g^{+}_{-1/2}=1/3g^{n}_{0}+5/6 g^{n}_{0}-1/6g^{n}_{1}  \\ 
g^{-}_{-1/2}=-1/6g^{n}_{0}+5/6 g^{n}_{0}+1/3g^{n}_{0}\\
\end{cases}
\]
\paragraph{Periodic boundary conditions} 
\[
\begin{cases}
g^{+}_{-1/2}=1/3g^{n}_{N-1}+5/6 g^{n}_{0}-1/6g^{n}_{1}  \\ 
g^{-}_{{}-1/2}=-1/6g^{n}_{N-2}+5/6 g^{n}_{N-1}+1/3g^{n}_{0}\\
g^{+}_{N-1/2}=1/3g^{n}_{N-1}+5/6 g^{n}_{0}-1/6g^{n}_{1}  \\ 
g^{-}_{N-1/2}=-1/6g^{n}_{N-2}+5/6 g^{n}_{N-1}+1/3g^{n}_{0}\\
\end{cases}
\]

\newpage

\section{Limiters}
\label{sec:limiter}
High-order numerical schemes may create spurious oscillations when stiff profiles occur in the distribution function. As usual in the numerical framework of finite volume schemes, we employ a flux limiter to lower these oscillations. The principle of the limiter is to introduce numerical diffusion when stiff profiles are detected, by modifying the flux at the cell faces. We will test limiters found in the literature, the so called Oscillation Limiter (OSL) described in \cite{crouseilles},  and a the one proposed in this report that we call Slope Limited Splines (SLS). These limiters are mainly provided with the PSM scheme, but some might be used for any finite volume scheme. 
\subsection{ENTropic flux limiter (ENT)}
The principle of this limiter proposed in \cite{braeunig} is to make degenerate the fourth order PSM scheme to a second order centred flux to reduce the anti-diffusive behaviour of the scheme when it occurs. The position where anti-diffusion occurs is detected looking at the second order equivalent equation solved by the scheme, obtained by a Taylor expansion. This equivalent equation shows a diffusion term at second order, which the sign should be positive, leading to numerical diffusion and thus stability. When this sign is negative, the scheme is anti-diffusive and it is replaced by a centred scheme. The corresponding algorithm is the following:  
\paragraph{Algorithm}
\begin{itemize}
\item Computation of the PSM flux $\phi^{PSM}_{i+1/2}$ and the centred flux $\phi^{CEN}_{i+1/2}=\alpha_{i+1/2}\frac{\bar{g}^{n}_{i}+\bar{g}^{n}_{i+1}}{2}$. 
\item if $\left(\phi^{CEN}_{i+1/2}-\phi^{PSM}_{i+1/2}\right)\left(\bar{g}^{n}_{i+1}-\bar{g}^{n}_{i}\right)>0$ then the flux $\phi^{PSM}_{i+1/2}$ is supposed to be diffusive.
\item if $\left(\phi^{CEN}_{i+1/2}-\phi^{PSM}_{i+1/2}\right)\left(\bar{g}^{n}_{i+1}-\bar{g}^{n}_{i}\right)<0$ then the flux $\phi^{PSM}_{i+1/2}$ is supposed to be anti-diffusive, thus it is switched with the centred flux:
\[\phi^{PSM}_{i+1/2}=\phi^{CEN}_{i+1/2}
\]
\end{itemize}
\subsection{UMEDA's limiter (UMEDA)}
The 4D advection equation is split in four 1D advection equations. Although the maximum principle is satisfied by the 4D equation, this principle is not fulfilled for each 1D equation. Therefore the extrema of the distribution function (minimum and maximum) are unknown in 1D. The PFC limiter \cite{filbet} employ these extrema, which are not known in the context of a directional splitting. Therefore, following {\it Umeda} \cite{umeda}, we modify the extrema definition to get a limiter working in 4D. The Umeda's limiter have been written with Lagrange polynoms.
\paragraph{Algorithm}
\begin{itemize}
\item We evaluate ${g}_{min1},{g}_{min2},{g}_{max1},{g}_{max2}$ 
\[
\begin{cases}
{g}_{min1}=\max\left[\max\left(\bar{g}^{n}_{i-1},\bar{g}^{n}_{i}\right),\min\left(2\bar{g}^{n}_{i-1}-\bar{g}^{n}_{i-2},2\bar{g}^{n}_{i}-\bar{g}^{n}_{i+1}\right)\right]\\
{g}_{min2}=\max\left[\max\left(\bar{g}^{n}_{i+1},\bar{g}^{n}_{i}\right),\min\left(2\bar{g}^{n}_{i+1}-\bar{g}^{n}_{i+2},2\bar{g}^{n}_{i}-\bar{g}^{n}_{i-1}\right)\right]\\
{g}_{max1}=\min\left[\min\left(\bar{g}^{n}_{i-1},\bar{g}^{n}_{i}\right),\max\left(2\bar{g}^{n}_{i-1}-\bar{g}^{n}_{i-2},2\bar{g}^{n}_{i}-\bar{g}^{n}_{i+1}\right)\right] \\
{g}_{max2}=\min\left[\min\left(\bar{g}^{n}_{i+1},\bar{g}^{n}_{i}\right),\max\left(2\bar{g}^{n}_{i+1}-\bar{g}^{n}_{i+2},2\bar{g}^{n}_{i}-\bar{g}^{n}_{i-1}\right)\right]
\end{cases}
\]
\item We set ${g}_{min},{g}_{max}$ 

\[
\begin{cases}
\bar{g}_{min}=\max\left[0,\min\left({g}_{min1},{g}_{min2}\right)\right]
\\ \bar{g}_{max}=\max\left[{g}_{max1},{g}_{max2}\right]
\end{cases}
\]
\item We define $L^{+}_i$ and $L^{-}_i$
\[
L^{+}_i=\begin{cases}
\text{if } \bar{g}_{i+1}-\bar{g}_{i}\geq 0, \min(2(\bar{g}_i-\bar{g}_{min}),\bar{g}_{i+1}-\bar{g}_i) \\ 
\text{if } \bar{g}_{i+1}-\bar{g}_{i}< 0, \max(2(\bar{g}_i-\bar{g}_{max}),\bar{g}_{i+1}-\bar{g}_i) 
\end{cases}
\]
\[
L^{-}_i=\begin{cases}
\text{if } \bar{g}_{i}-\bar{g}_{i-1}\geq 0, \min(2(\bar{g}_{max}-\bar{g}_i),\bar{g}_{i}-\bar{g}_{i-1}) \\
\text{if } \bar{g}_{i}-\bar{g}_{i-1}< 0, \max(2(\bar{g}_{min}-\bar{g}_i),\bar{g}_{i}-\bar{g}_{i-1}) 
\end{cases}
\]
\item We finally redefine the flux: \[G_{i+1/2}(x)=\bar{g}_i+x(1-x)(2-x)(L^{+}_i/6)+x(1-x)(1+x)(L^{-}_i/6)\]
\end{itemize}

 We obtain the LAG reconstruction without limiter by setting:
  \[\begin{cases}
L^{+}_k= \bar{g}_{k+1}-\bar{g}_{k}\\
L^{-}_k= \bar{g}_{k}-\bar{g}_{k-1}\\
\end{cases}
\]

\subsection{OScillation Limiter (OSL)}
The OScillation Limiter (OSL) proposed in \cite{crouseilles} is really using the Hermite formalism. It compares the values at call faces $g^{\pm}_{k-1/2}$ obtained with the LAG and the PSM schemes to the value computed with a linear reconstruction. This later consists in computing an average at the face of the left and right cells centred values. If the values computed with PSM and LAG are not simultaneously upper or lower than the average value, then we take the average value. If not, a mixed scheme between PSM, LAG and the average reconstruction is performed. The OSL limiter includes a parameter C>1 determining the proportion in the average of PSM and LAG fluxes in the limited flux.
     \paragraph{Algorithm}
     \begin{itemize}
     \item Computation of $g^{\pm}_{k-1/2}$ values for both PSM et LAG schemes.
\item Average value at the $k-1/2$ and $k+1/2$ nodes: $g^{ave}_{k-1/2}=(\bar{g}^{n}_k+\bar{g}^{n}_{k-1})/2$ and  $g^{ave}_{k+1/2}=(\bar{g}^{n}_k+\bar{g}^{n}_{k+1})/2$.
 \item The following formula perform the choice for the limited flux according to the regularity of the cell faces values: 
 \begin{itemize}
\item if ($g^{+}_{k-1/2,LAG}-g^{ave}_{k-1/2})(g_{k-1/2,PSM}-g^{ave}_{k-1/2}) >0$ then\\
$g^{+}_{k-1/2}=g^{ave}_{k-1/2}$
else
\begin{equation}
\begin{array}{l}
g^{+}_{k-1/2}=g^{ave}_{k-1/2} \\
+sign(g_{k-1/2,PSM}-g^{ave}_{k-1/2}) \\
\min(C|g^{+}_{k-1/2,LAG}-g^{ave}_{k-1/2}|,|g_{k-1/2,PSM}-g^{ave}_{k-1/2}|)
\end{array}
\end{equation}

 \item if ($g^{-}_{k+1/2,LAG}-g^{ave}_{k+1/2})(g_{k+1/2,PSM}-g^{ave}_{k+1/2}) >0$ then\\
$g^{-}_{k+1/2}=g^{ave}_{k+1/2}$
else
\begin{equation}
\begin{array}{l}\label{osl}
g^{-}_{k+1/2}=g^{ave}_{k+1/2} \\
+sign(g_{k+1/2,PSM}-g^{ave}_{k+1/2}) \\
\min(C|g^{-}_{k+1/2,LAG}-g^{ave}_{k+1/2}|,|g_{k+1/2,PSM}-g^{ave}_{k+1/2}|)
\end{array}
\end{equation}
\end{itemize}
\end{itemize}

\subsection{Slope Limited Splines (SLS)}
\begin{figure}
 \centering
 \begin{tabular}{c}
   \includegraphics[scale=0.35]{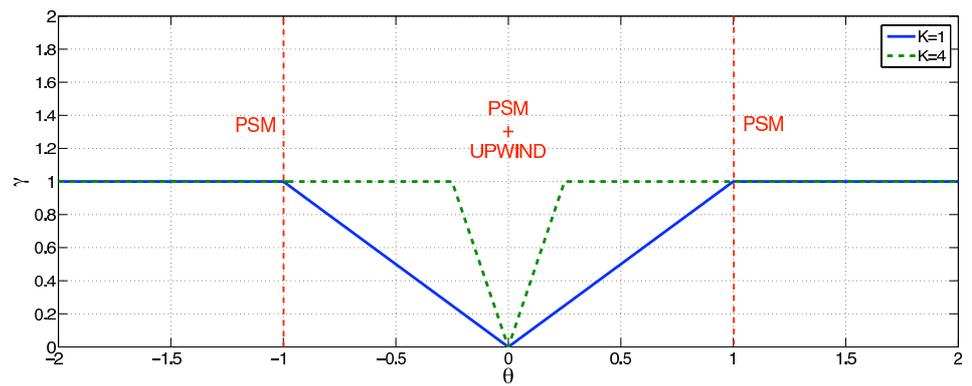} \\
 \end{tabular}
 \caption{$\gamma$ function for the SLS limiter}
 \label{fig:gamma_theta}
 \end{figure}

The limiting procedure basically aims to cut the oscillations generated by strong gradients in the distribution function profile, where high order schemes will produce overshoots and spurious oscillations. We propose here to measure these gradients and to add diffusion where strong gradients are detected. The diffusion is added by mixing the high order scheme with a first order upwind flux. The more the gradient is steep, the more we raise the proportion of upwind flux in the average with the high order scheme. The evaluation of the gradient is given by a function $\theta$ and we estimate the diffusion needed with a function $\gamma(\theta) \in [0,1]$ based on the minmod like limiter function (\fref{gamma_theta}):
\[
\phi^{new}_{i+1/2}=\gamma(\theta_{i+1/2}) ~\phi^{PSM}_{i+1/2}+(1-\gamma(\theta_{i+1/2})) ~\phi^{upwind}_{i+1/2}
\]
 where 
\[\begin{cases}
\phi^{upwind}_{i+1/2}= \alpha_{i+1/2}\frac{\bar{g}^{n}_{i}+\bar{g}^{n}_{i+1}}{2}- sign(\alpha_{i+1/2}) \frac{\bar{g}^{n}_{i+1}-\bar{g}^{n}_{i}}{2}\\
 \alpha_{i+1/2}= \Delta t a_{i+1/2}
\end{cases}
\] 
We define $\theta_{i+1/2}$ as the classical slope ratio of the distribution which depends on the direction of the displacement (\fref{thetapente}):
\[ \theta_{i+1/2}= \left\{
    \begin{array}{ll}
        \frac{\color{blue}\bar{g}^{n}_{i}-\bar{g}^{n}_{i-1}}{\color{red}\bar{g}^{n}_{i+1}-\bar{g}^{n}_{i}} & \text{ if } \alpha_{i+1/2}>0  \\
      \frac{\color{green}\bar{g}^{n}_{i+2}-\bar{g}^{n}_{i+1}}{\color{red}\bar{g}^{n}_{i+1}-\bar{g}^{n}_{i}} & \text{ if } \alpha_{i+1/2}<0   \end{array}
\right.
\]
However, the classical limiters as minmod $\gamma_{i+1/2}=max(0,min(\theta_{i+1/2},1))$, set $\gamma$ to 0 when $\theta<0$. That means that the scheme turns to order 1 when an extrema exists, i.e. the slope ratio $\theta<0$. These extrema are thus quickly diffused and that leads to loose the benefits of a high order method.
For SLS, the choice is to let the high-order scheme deal with the extrema and only add diffusion when high gradients occurs, i.e.  the slope ratio $\theta \approx 0$. We also introduce a constant K in relation to control the maximum slope allowed without adding diffusion. The SLS limiter function is thus to set  $\gamma=1$ for any values of $\theta$, except close to $|\theta|=0$ where strong gradients occurs, see figure \ref{fig:gamma_theta}:  
\begin{equation}\label{sls}
\gamma_{i+1/2}=\max(0,\min(K|\theta_{i+1/2}|,1))
\end{equation}

\begin{figure}
\begin{center}
\scalebox{1.4} 
{
\begin{pspicture}(0,-1.858125)(10.202812,1.358125)
\definecolor{color68}{rgb}{0.6,0.6,0.6}
\definecolor{color130}{rgb}{0.4,0.4,0.4}
\definecolor{color463}{rgb}{0.0,0.0,0.8}
\psline[linewidth=0.04cm](0.8609375,-0.5003125)(8.4609375,-0.5003125)
\psline[linewidth=0.04cm](1.2609375,-0.3003125)(1.2609375,-0.7203125)
\psline[linewidth=0.04cm](2.0609374,-0.3003125)(2.0609374,-0.7203125)
\psline[linewidth=0.04cm](2.8409376,-0.3003125)(2.8409376,-0.7203125)
\psline[linewidth=0.04cm](3.6409376,-0.3003125)(3.6409376,-0.7203125)
\psline[linewidth=0.04cm](4.4609375,-0.3003125)(4.4609375,-0.7203125)
\psline[linewidth=0.04cm](5.2609377,-0.3003125)(5.2609377,-0.7203125)
\psline[linewidth=0.04cm](6.0609374,-0.3003125)(6.0609374,-0.7203125)
\psline[linewidth=0.04cm](6.8609376,-0.3003125)(6.8609376,-0.7203125)
\psline[linewidth=0.04cm](7.6409373,-0.2803125)(7.6409373,-0.7003125)
\psline[linewidth=0.04cm,linecolor=color68](1.2809376,0.0796875)(2.0409374,0.0796875)
\psline[linewidth=0.04cm,linecolor=color68](2.0809374,-0.1403125)(2.8409376,-0.1403125)
\psline[linewidth=0.04cm,linecolor=color68](2.8609376,-0.8003125)(3.6209376,-0.8003125)
\psline[linewidth=0.04cm,linecolor=color68](3.6809375,-0.2803125)(4.4409375,-0.2803125)
\psline[linewidth=0.04cm,linecolor=color68](6.8809376,0.4796875)(7.6409373,0.4796875)
\psline[linewidth=0.04cm,linecolor=color68](6.1009374,0.1396875)(6.8609376,0.1396875)
\psline[linewidth=0.04cm,linecolor=color68](5.2809377,-0.3603125)(6.0409374,-0.3603125)
\psline[linewidth=0.04cm,linecolor=color68](7.6809373,0.8396875)(8.440937,0.8396875)
\psdots[dotsize=0.15,linecolor=color130](8.040937,0.8396875)
\psdots[dotsize=0.15,linecolor=color130](7.2609377,0.4796875)
\psdots[dotsize=0.15,linecolor=color130](6.4809375,0.1396875)
\psdots[dotsize=0.15,linecolor=color130](5.6809373,-0.3603125)
\psdots[dotsize=0.15,linecolor=color130](4.0809374,-0.2803125)
\psdots[dotsize=0.15,linecolor=color130](2.4809375,-0.1203125)
\psdots[dotsize=0.15,linecolor=color130](1.6409374,0.0996875)
\psdots[dotsize=0.15,linecolor=color130](3.2609375,-0.8203125)
\psframe[linewidth=0.04,linecolor=white,dimen=outer](5.2009373,-0.1203125)(4.5209374,-0.7003125)
\psframe[linewidth=0.04,linecolor=white,dimen=outer](5.1009374,-0.2603125)(4.6409373,-0.6603125)
\psframe[linewidth=0.04,linecolor=white,dimen=outer](5.0009375,-0.2403125)(4.6609373,-0.6403125)
\psframe[linewidth=0.04,linecolor=white,dimen=outer](4.9009376,-0.3203125)(4.7609377,-0.6003125)
\usefont{T1}{ptm}{m}{n}
\rput(8.462344,1.1696875){\color{green}$\bar{g}^{n}_{i+2}$}
\usefont{T1}{ptm}{m}{n}
\rput(4.382344,0.196875){\color{green}$\bar{g}^{n}_{i+2}$}
\usefont{T1}{ptm}{m}{n}
\rput(7.422344,0.1096875){\color{red}$\bar{g}^{n}_{i+1}$}
\usefont{T1}{ptm}{m}{n}
\rput(3.2023437,-1.1303124){\color{red}$\bar{g}^{n}_{i+1}$}
\usefont{T1}{ptm}{m}{n}
\rput(2.8123438,0.2096875){\color{red}$\bar{g}^{n}_{i}$}
\usefont{T1}{ptm}{m}{n}
\rput(6.412344,0.5096875){\color{red}$\bar{g}^{n}_{i}$}
\usefont{T1}{ptm}{m}{n}
\rput(5.682344,-0.9103125){\color{blue}$\bar{g}^{n}_{i-1}$}
\usefont{T1}{ptm}{m}{n}
\rput(1.6223438,-0.2103125){\color{blue}$\bar{g}^{n}_{i-1}$}
\rput(3.2023437,-1.903124){left}
\rput(6.2023437,-1.903124){right}
\psline[linewidth=0.051999997cm,linecolor=green,arrowsize=0.05291667cm 2.0,arrowlength=1.4,arrowinset=0.4]{<->}(3.6409376,-0.8403125)(3.6409376,-0.26103125)
\psline[linewidth=0.051999997cm,linecolor=green,arrowsize=0.05291667cm 2.0,arrowlength=1.4,arrowinset=0.4]{<->}(7.6809373,0.4396875)(7.6609373,0.8796875)
\psline[linewidth=0.051999997cm,linecolor=red,arrowsize=0.05291667cm 2.0,arrowlength=1.4,arrowinset=0.4]{<->}(6.8809376,0.0796875)(6.8609376,0.5196875)
\psline[linewidth=0.051999997cm,linecolor=red,arrowsize=0.05291667cm 2.0,arrowlength=1.4,arrowinset=0.4]{<->}(2.8609376,-0.8003125)(2.8609376,-0.1203125)
\psline[linewidth=0.051999997cm,linecolor=color463,arrowsize=0.05291667cm 2.0,arrowlength=1.4,arrowinset=0.4]{<->}(2.0609374,-0.2403125)(2.0609374,0.1796875)
\psline[linewidth=0.051999997cm,linecolor=color463,arrowsize=0.05291667cm 2.0,arrowlength=1.4,arrowinset=0.4]{<->}(6.0809374,0.1996875)(6.0809374,-0.4203125)
\end{pspicture} 
}
\caption{Computation of $\theta_{i+1/2}$}
\label{fig:thetapente}

\end{center}\end{figure}
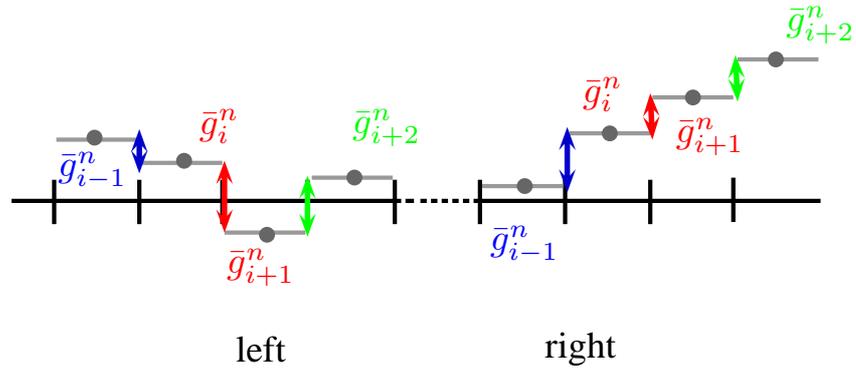

\section{Numerical results}
\label{sec:numerical_results}
We consider both 1D and 4D test cases. The 1D advection test case is performed with a constant velocity to give a qualitative view on the effects of each limiter. We are particularly interested in advection of profiles with strong gradients,  as discontinuous initial distribution functions as a step. However, we have to keep in mind that this situation of a discontinuous shape advected with a constant velocity field does not occur in Vlasov 4D drift-kinetic simulations for many reasons: the advections are not constant, the long time constant advection configuration does not appear, and discontinuous functions (steps) do not exist in the Vlasov model. However, it is relevant to investigate advection of discontinuous functions  since Vlasov models leads to stretched structures in the flow that makes appear strong gradients, which could be assimilated to discontinuities.\\

The 4D drift-kinetic test case is performed with the 5D Gysela code following the benchmark presented in \cite{grandgirard}. This benchmark is to simulate instabilities growing in the plasma leading to shapes like filaments and vortex. We compare the results to evaluate the limiters effects on the development of turbulent structures. We will focus on two directions, $(r,\theta)$, among the four $(r,\theta, \phi, v_{\parallel})$ since strong gradients of the distribution function appear mainly in $(r,\theta)$ planes (further details in \cite{braeunig}). We give algorithms used in Gysela to solve the advection equation in \sref{annexeB}.                

\subsection{Tools}
We propose hereafter some quantitative tools to compare the efficiency of each limiter. The $L^2$ norm and the entropy may be used to investigate the diffusivity of each limiter, since each should theoretically be kept constant. The total energy conservation is related to the treatment of small structures provided by the instabilities. In addition, we also introduce the total variation (TV) norm to estimate the numerical oscillations created by the scheme. 
\paragraph{ Total Variation\\}
The total variational norm (TV) is used to estimate the rate of oscillations produced by a scheme, compare to another.
\begin{itemize}
\item In \textbf{1D}:  \\
$TV(g(t))=\displaystyle  \int^{x_{max}}_{x_{min}} |\deriv{g(t,x)}{x}|dx$
  \[
TV(g^n)=\frac{1}{\Delta x}\sum^{N_x}_{i=0} |\bar{g}^n_{i+1}-\bar{g}^n_{i}|
\]
\item In \textbf{4D}, we restrict the diagnostic to an $(r,\theta)$ plane:  \\
$TV(g(t))=\displaystyle \int^{r_{max}}_{r_{min}}\int^{2 \pi}_0 |\nabla_{r,\theta} ~{g(t,r,\theta)}|_{L^2} dr d\theta $
\[
TV(g^n)=\sum^{N_r}_{i=0}\sum^{N_\theta}_{j=0} \sqrt{ \left(\frac{\bar{g}^n(i+1,j)-\bar{g}^n(i,j)}{\Delta r }\right)^2+\left(\frac{\bar{g}^n(i,j+1)-\bar{g}^n(i,j)}{\Delta \theta }\right)^2}
\]
\end{itemize}
\paragraph{ $L^2$ norm and entropy\\}
The $L^2$ norm is defined as: 
\[ 
L^2(f(t))=\int (f(t,X,v_{\parallel}))^2 dXdv_{\parallel},\:dX=rdr d\theta d\phi.
\] 
The entropy measures the information created or destroyed by a phenomena. We can thus expect that a raising entropy estimates the information lost by numerical diffusion.
\[ S(f(t))=-\int f(t,X,v_{\parallel})log(f(t,X,v_{\parallel})) dXdv_{\parallel},\: dX=rdrd\theta d\phi \] 
 \paragraph{Total Energy\\}
The conservation of the total energy $\epsilon_{tot}$ writes as the sum of the kinetic energy $\epsilon_{kin}$ and the potential energy $\epsilon_{pot}$ defined as follow: 
\[
\epsilon_{tot}=\epsilon_{kin}+\epsilon_{pot}=\int \frac{1}{2}m_i v^2_{\parallel} (f-f_{eq})dVdv_{\parallel}+\frac{1}{2}\int e\phi(n_i-n_0)dX.
\]
For further details, see \cite{grandgirard}.

\paragraph{Quality factor\\}
We propose here to define a quality factor to establish a quantitative evaluation of the limiters efficiency. Since the $L^2$ norm should be conserved by the Vlasov equation, it furnishes a good estimation of the diffusive behaviour of a scheme, because it decreases this norm. On the other hand, the TV norm estimates the numerical oscillations, or also the diffusion in a way, but more 'locally'. The quality of a scheme may be defined as its ability to limit the oscillations (small TV) with the less numerical diffusion (high $L^2$ norm). Therefore we introduce the quality criterion $Q$ defined as: 
\begin{equation}\label{quality}
Q(f(t))=\frac{L^2(f(t))}{TV(f(t))}
\end{equation}

\subsection{Test of the limiters on constant 1D advections}
We test the limiters in 1D on a step function \footnote{the distribution function f is set to 1 on a part of the domain else to 0}. The step exhibits the problems for which the limiters are required. This 1D benchmark gives a first overview on the limiters capabilities.\\
The benchmark consists in solving a constant advection on a periodic domain divided in 80 cells. The displacement is set to 0.2 cell per iteration. The indicated times are the number of iterations.

\begin{figure}
 \centering
 \begin{tabular}{c}
   \includegraphics[width=\textwidth]{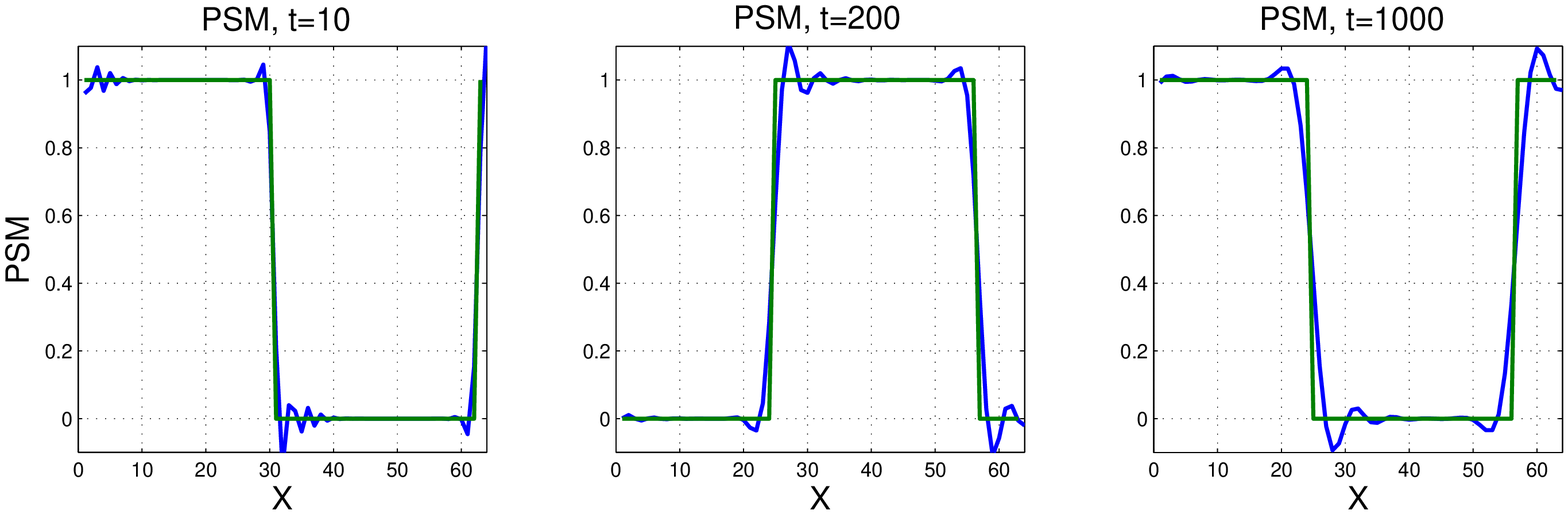} \\
   \includegraphics[width=\textwidth]{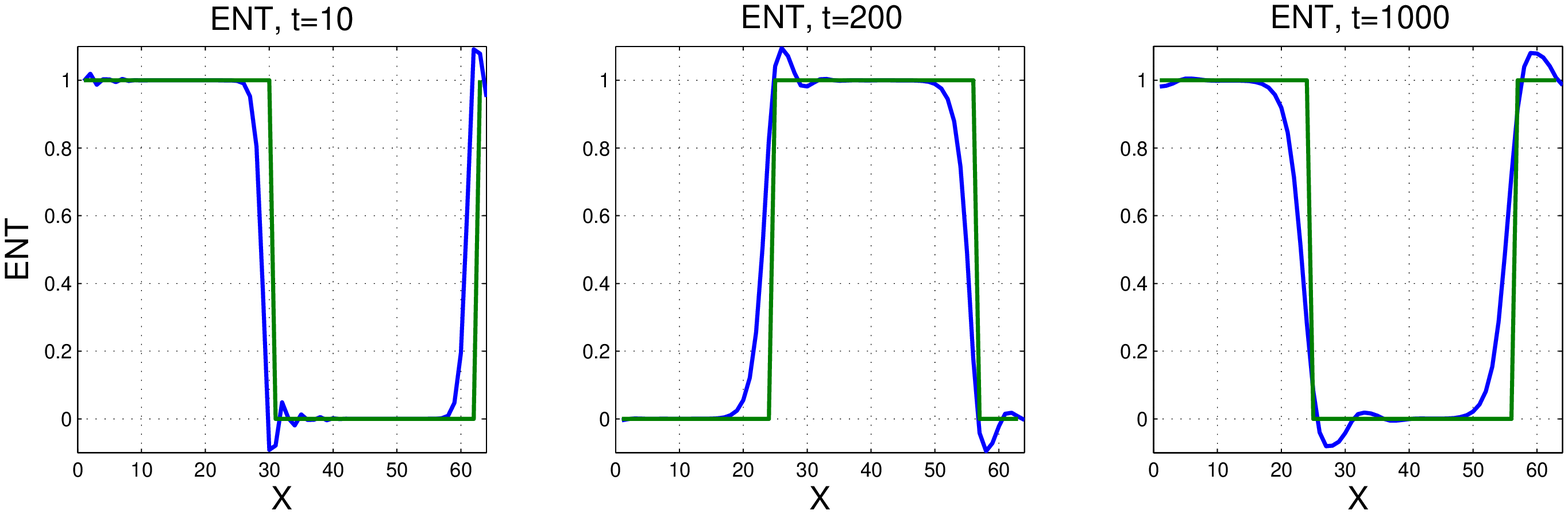} \\
 \end{tabular}
 \caption{Constant advection on a step with the PSM  scheme and the PSM scheme with entropic limiter ENT. The domain is meshed on 80 cells with periodic boundary conditions and the displacement is set to 0.2 cell per iteration.}
 \label{fig:psm-ent}
 \end{figure}

\newpage
\paragraph{PSM-ENT}
We observe on figure \ref{fig:psm-ent} that the ENT limiter cuts off the spurious oscillations of the PSM scheme  at the left side of the discontinuity. However, at the right side, oscillations are weakly damped.
 
\begin{figure}
 \centering
 \begin{tabular}{c}
    \includegraphics[width=\textwidth]{PSM-1.eps} \\
  \includegraphics[width=\textwidth]{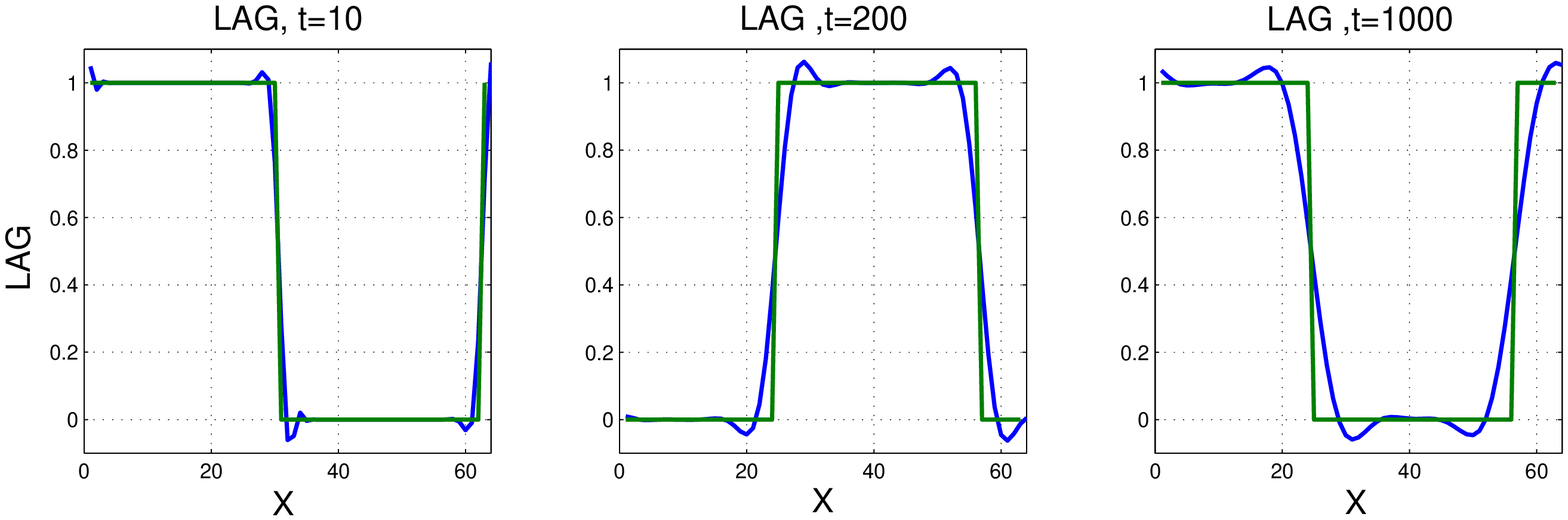} \\
   \includegraphics[width=\textwidth]{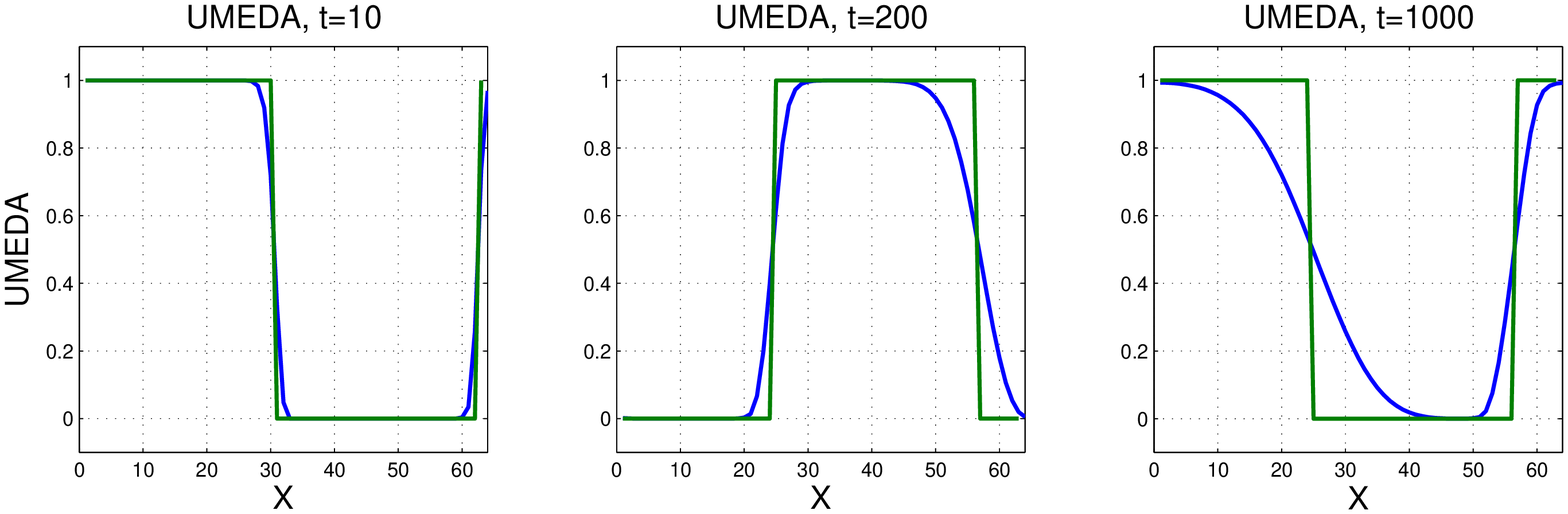} \\
 \end{tabular}
 \caption{Constant advection on a step with the LAG  scheme and the LAG scheme with the Umeda limiter. The domain is meshed on 80 cells with periodic boundary conditions and the displacement is set to 0.2 cell per iteration.}
 \label{fig:lag-umeda}
 \end{figure}

\newpage
\paragraph{LAG-UMEDA}
We observe on figure \ref{fig:lag-umeda} that the LAG scheme is much more diffusive than the PSM scheme without really better satisfying a maximum principle. The under/overshoots of the LAG scheme are cut off by the Umeda limiter. Although the LAG scheme with Umeda limiter respects a maximum principle, the step signal is quickly diffused and it seems not appropriate for long time simulations.
 
\begin{figure}
 \centering
 \begin{tabular}{c}
  \includegraphics[width=\textwidth]{PSM-1.eps} \\
  \includegraphics[width=\textwidth]{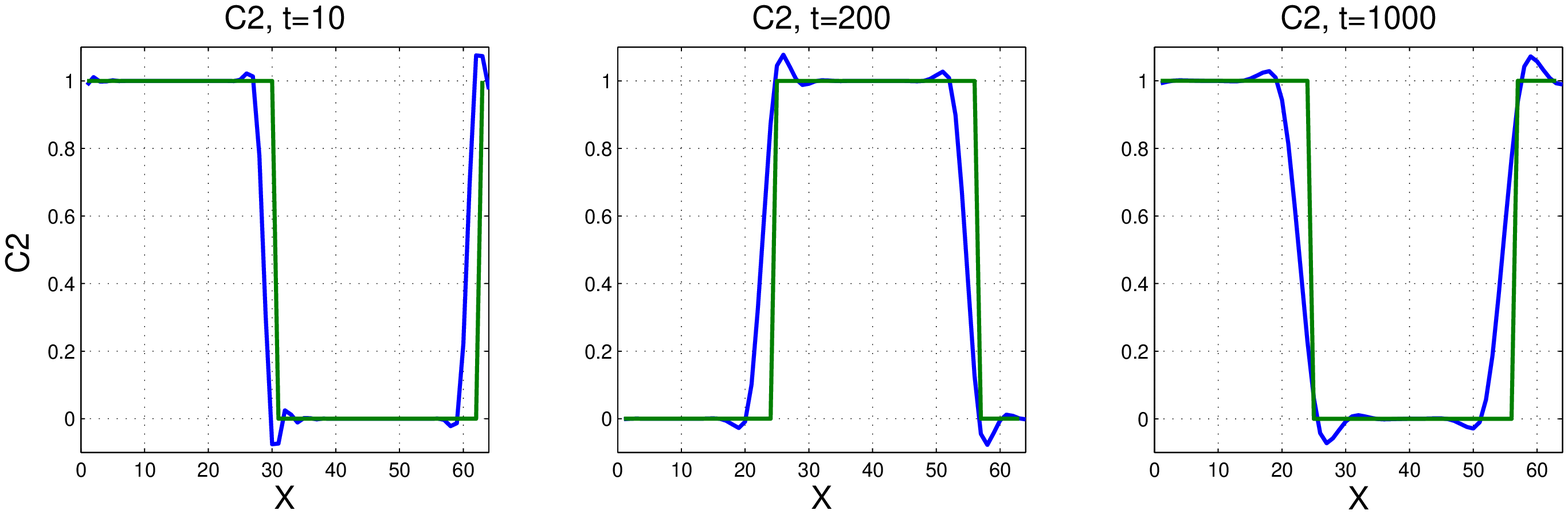} \\
  \includegraphics[width=\textwidth]{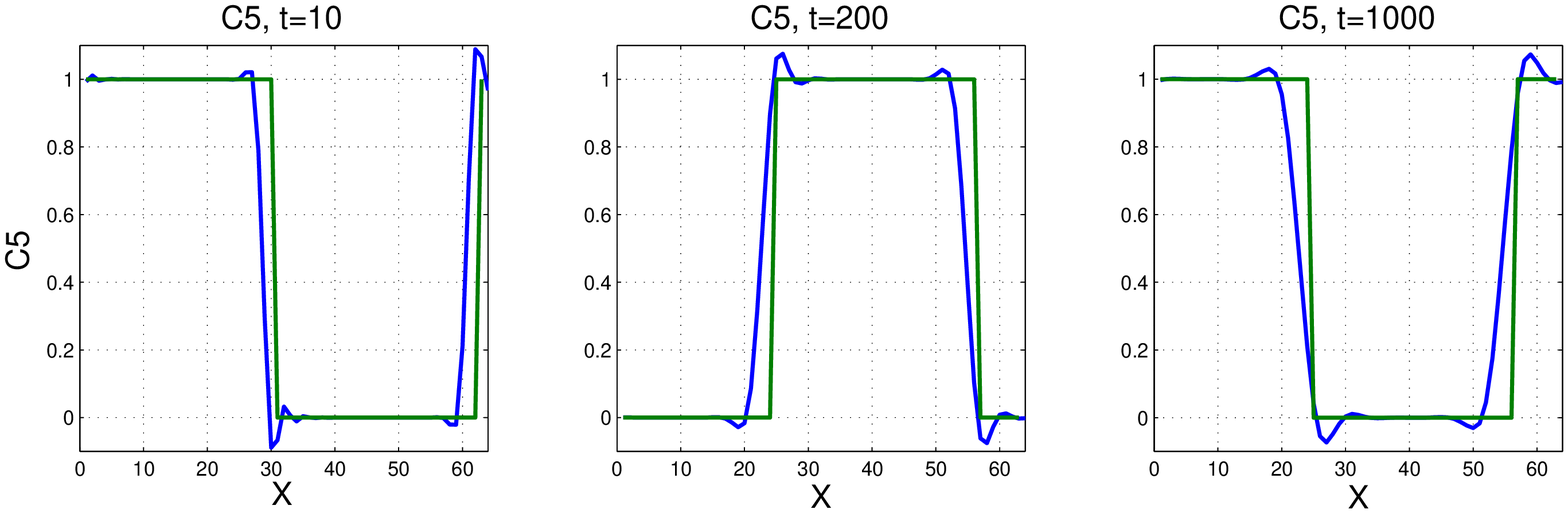} \\
  \includegraphics[width=\textwidth]{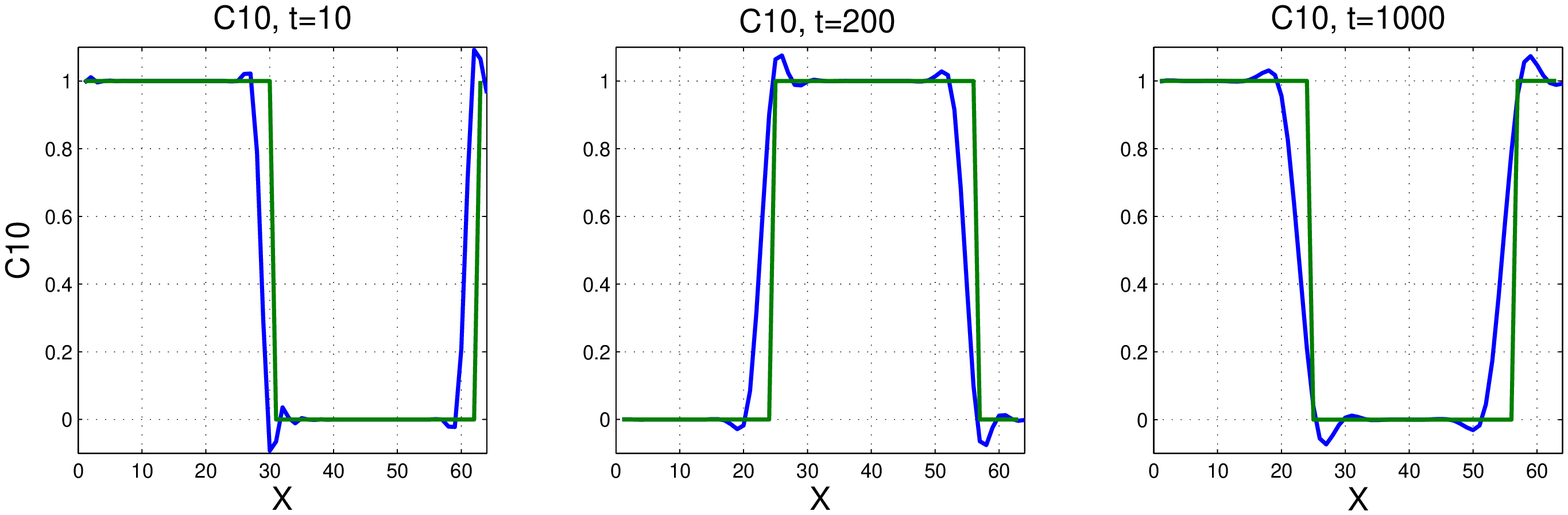} \\
 \end{tabular}
 \caption{Constant advection on a step with the PSM  scheme and the PSM scheme with the OSL limiter with different values of parameter $C$. The domain is meshed on 80 cells with periodic boundary conditions and the displacement is set to 0.2 cell per iteration.}
 \label{fig:psm-osl}
 \end{figure}

 \newpage
\paragraph{PSM-OSL}
We observe on figure \ref{fig:psm-osl} that the OSL limiter reduces the oscillations without introducing much diffusion. However, there is still under/overshoots at both sides of the discontinuity and a slight offset of the solution occurs compare to the exact solution, but it seems not increase with time. The constant $C$ of the OSL limiter has almost no influence on the results.

\begin{figure}
 \centering
 \begin{tabular}{c}
  \includegraphics[width=\textwidth]{PSM-1.eps} \\
  \includegraphics[width=\textwidth]{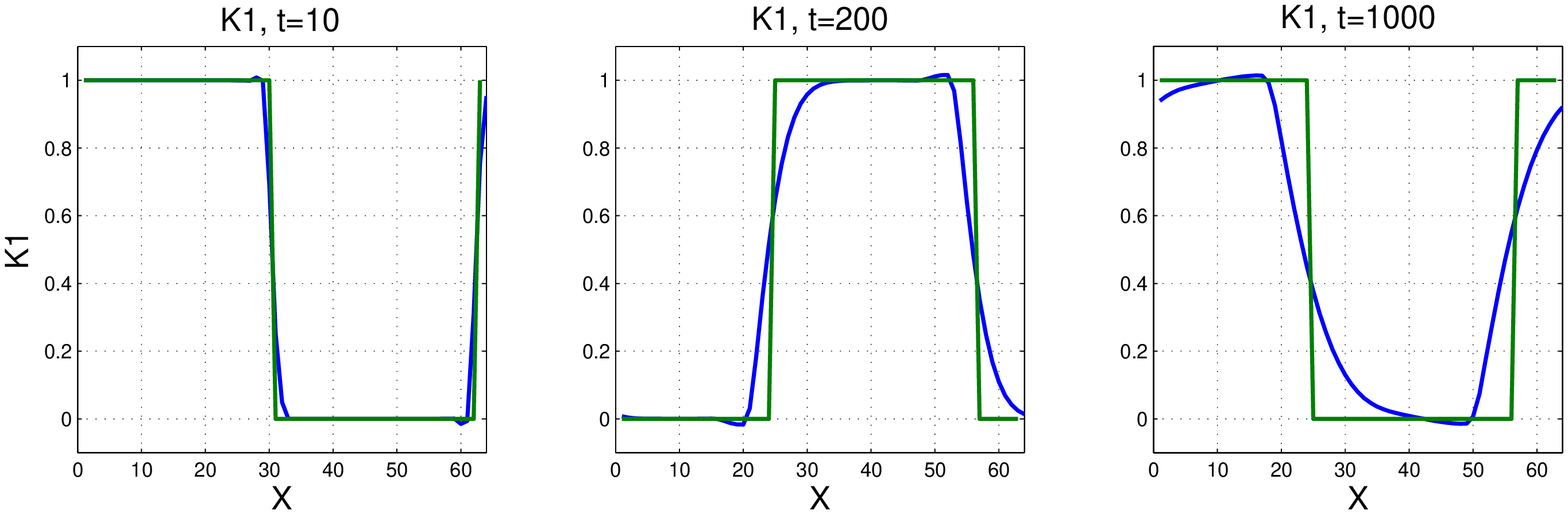}     \\
  \includegraphics[width=\textwidth]{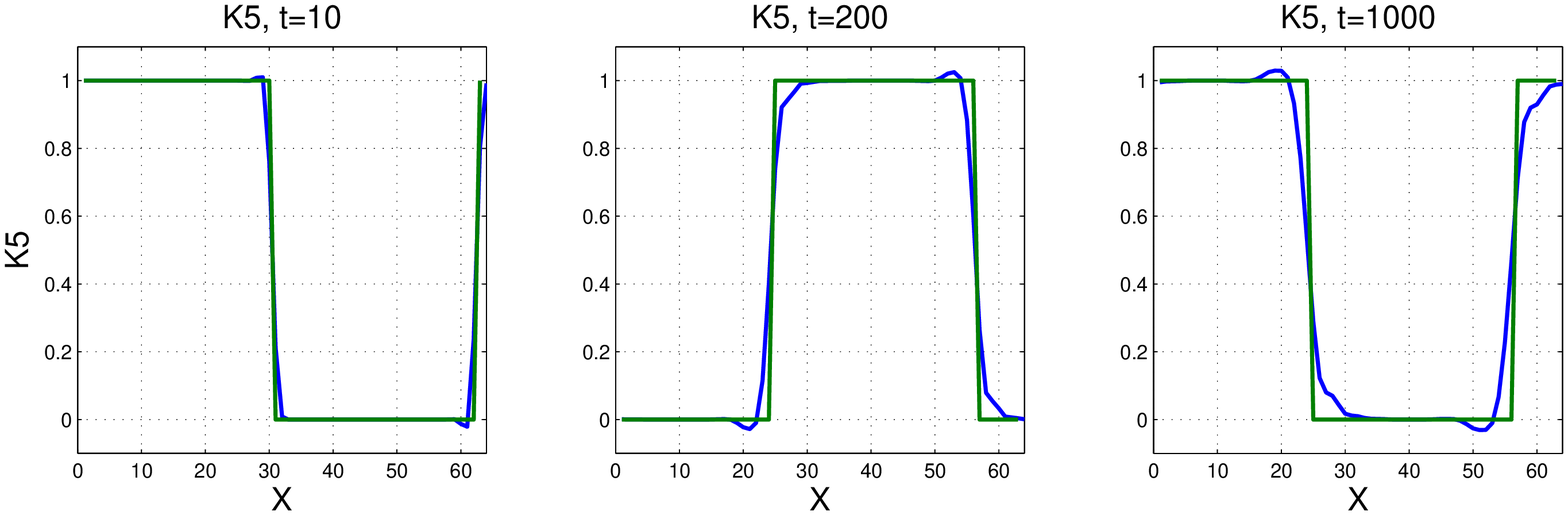}     \\
  \includegraphics[width=\textwidth]{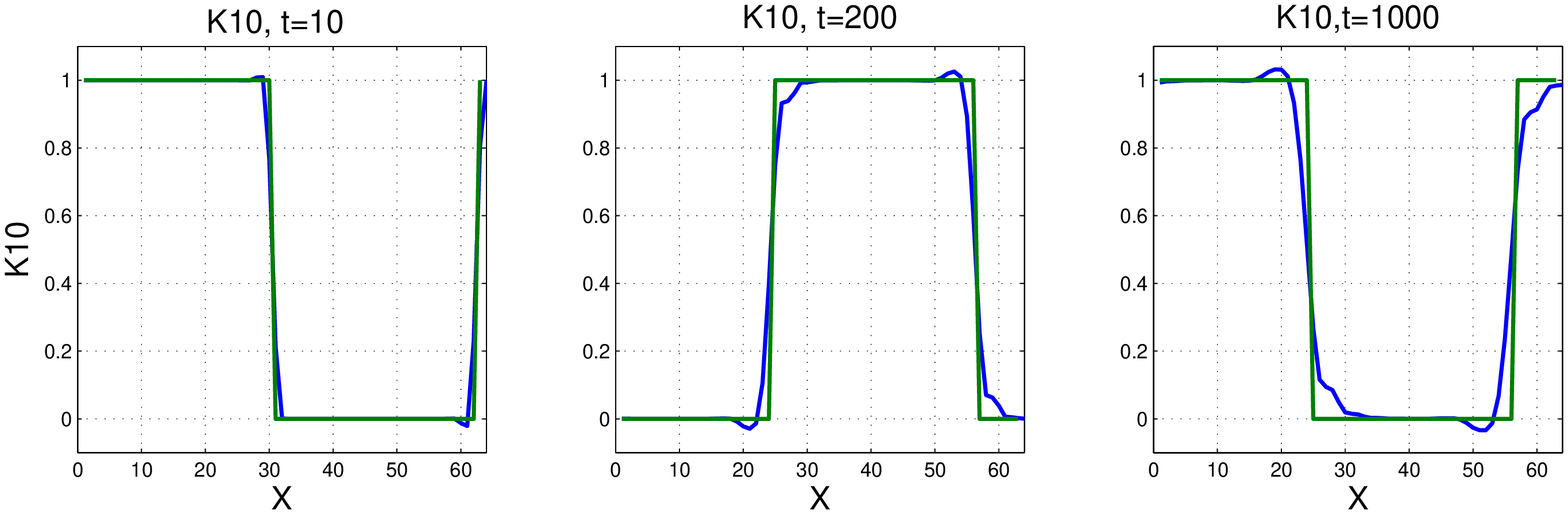}   \\
 \end{tabular}
 \caption{Constant advection on a step with the PSM  scheme and the PSM scheme with the SLS limiter with different values of parameter $K$. The domain is meshed on 80 cells with periodic boundary conditions and the displacement is set to 0.2 cell per iteration.}
  \label{fig:psm-sls}
 \end{figure}

\newpage
\paragraph{PSM-SLS}
We observe on figure \ref{fig:psm-sls} that the SLS limiter cuts off oscillations at the right side of the discontinuity and keeps the overshoot at left side. The results of SLS limiter with $K=1$ is too much diffusive (equivalent to a minmod limiter). With $K=5$, the result is almost the same as with  $K=10$ and the accuracy at discontinuities (slope of the solution at discontinuities) is almost the same as with the PSM scheme.    
  
\newpage
\paragraph{Comparison of all limiters for PSM}
It is interesting to observe the action of each different limiter for PSM on the oscillations. The ENT limiter cuts off the oscillations at left side of the discontinuities, the SLS $K=5$ limiter cuts off the oscillations  at right side of the discontinuities and the OSL limiter reduces the oscillations at both sides. 
 \begin{figure}
 \centering
 \begin{tabular}{c}
   \includegraphics[width=\textwidth]{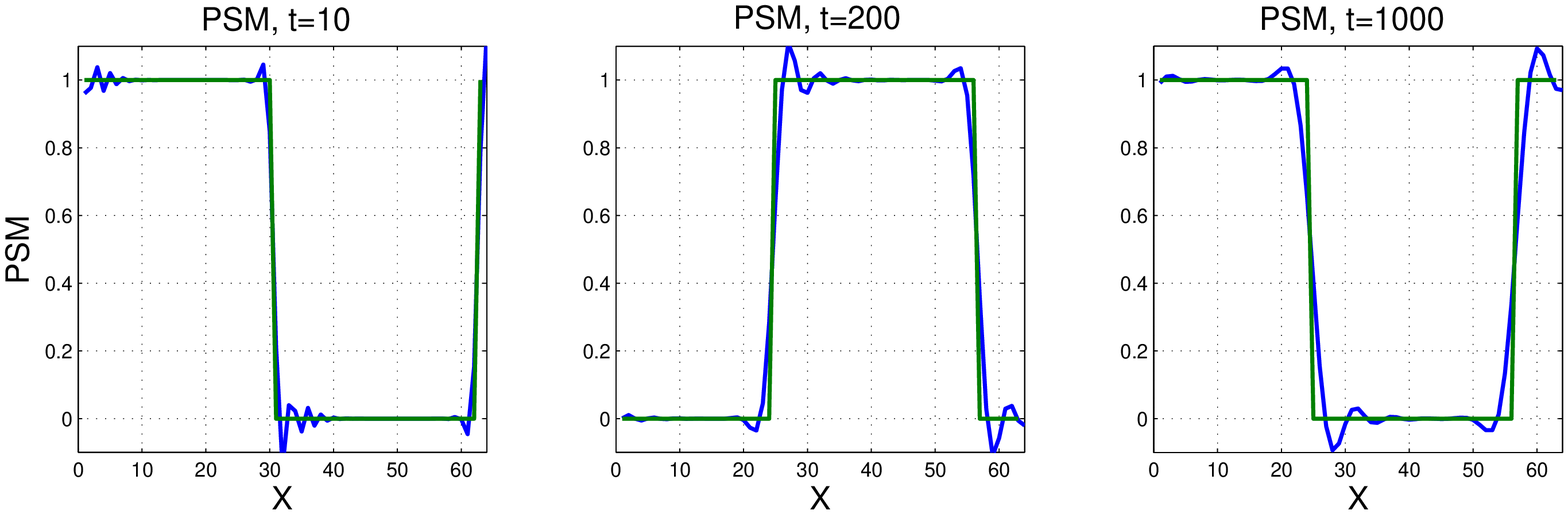} \\
   \includegraphics[width=\textwidth]{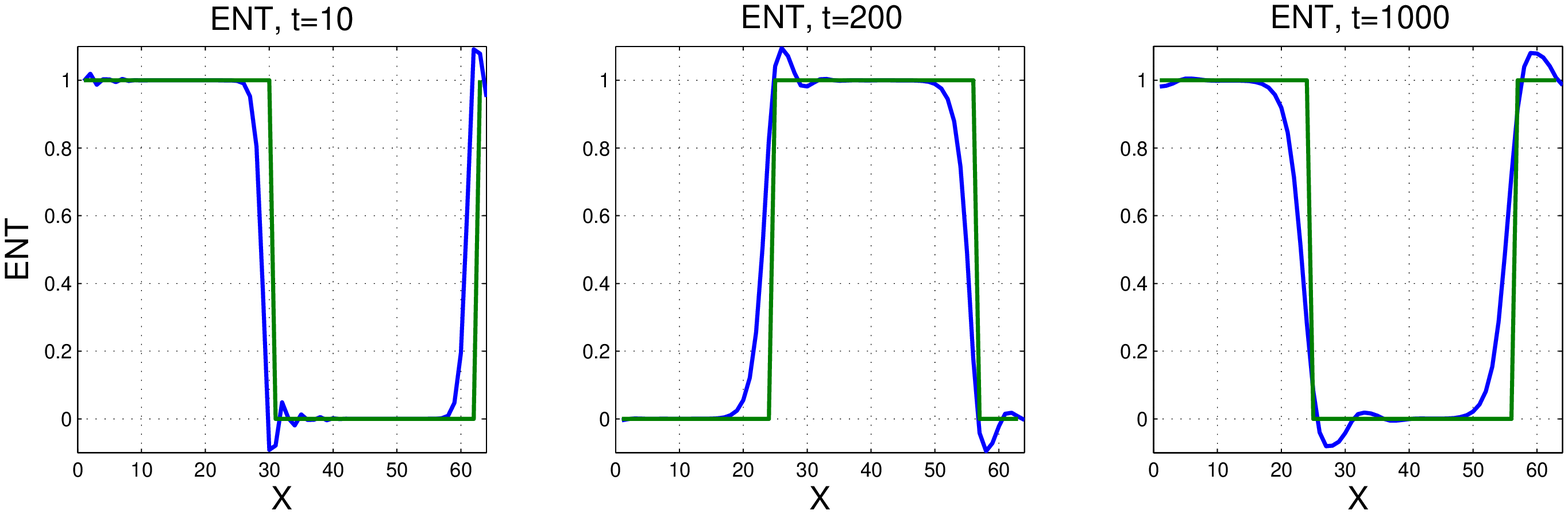} \\
   \includegraphics[width=\textwidth]{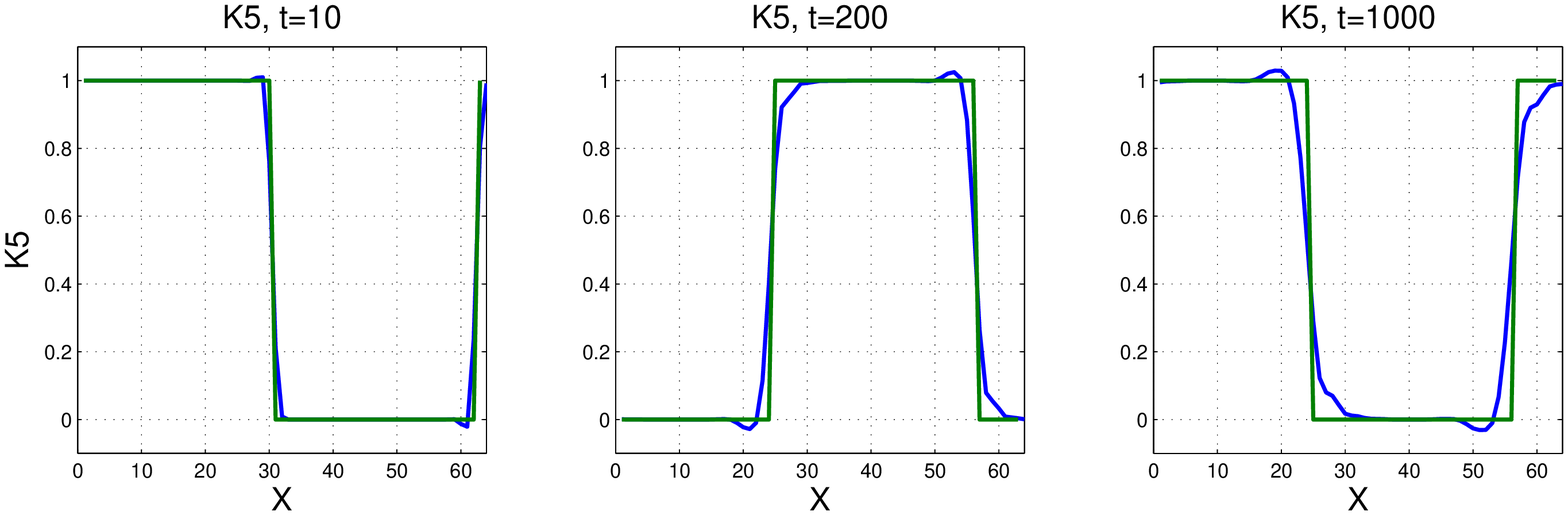} \\
   \includegraphics[width=\textwidth]{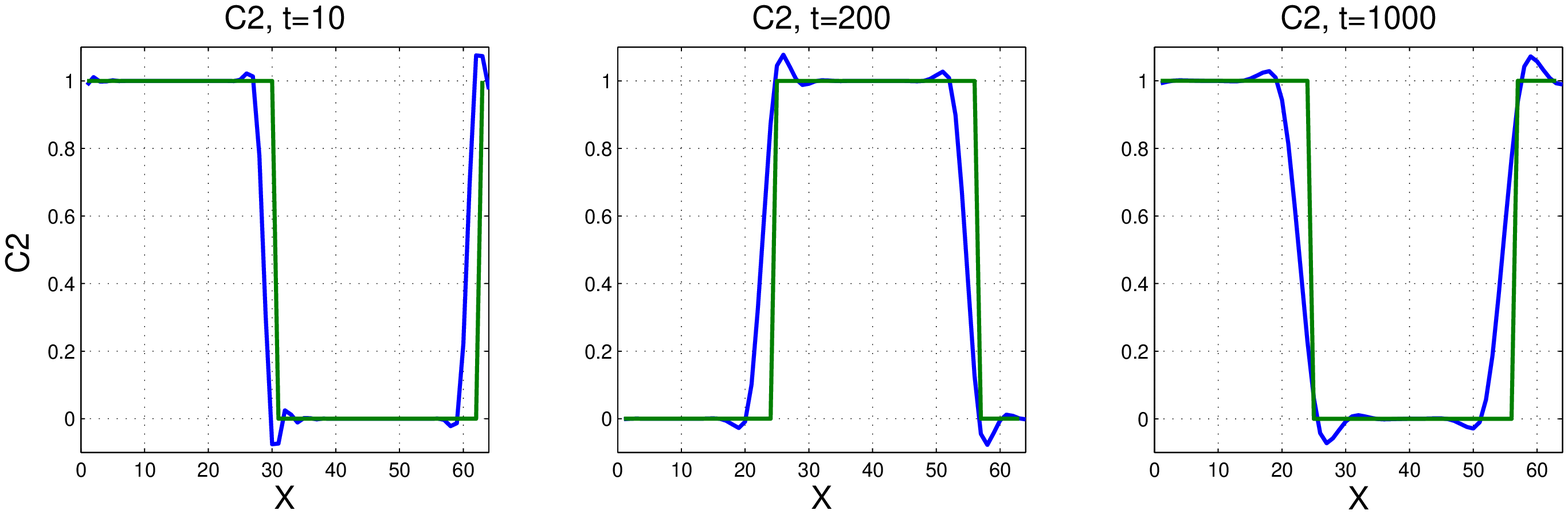} \\
 \end{tabular}
 \caption{Constant advection on a step with the PSM  scheme, the PSM scheme with the ENT limiter, the SLS limiter with parameter $K=5$ and the OSL limiter with parameter $C=2$. The domain is meshed on 80 cells with periodic boundary conditions and the displacement is set to 0.2 cell per iteration.}
 \label{fig:limpsmall}
 \end{figure}

 \newpage
\paragraph{Quality factor}
 \begin{figure}
 \centering
  \includegraphics[width=\textwidth]{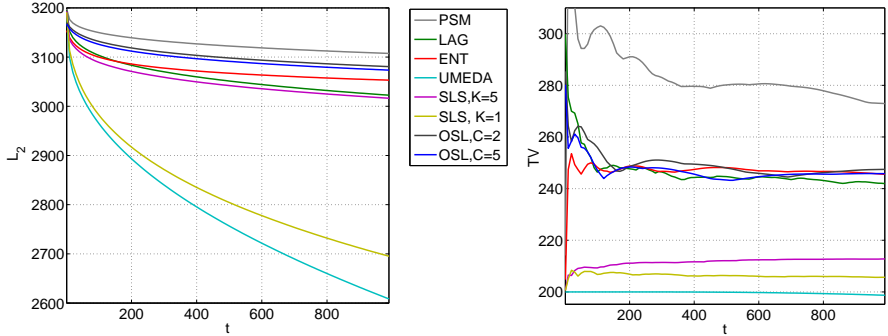}
 \caption{$L^2$ norm (left) and total variation (right) as a function of the time for the different methods applied to a step with 1D constant advection on a periodic domain. }
 \label{fig:}
 \end{figure}
 \begin{figure}
 \centering
  \includegraphics[scale=0.3]{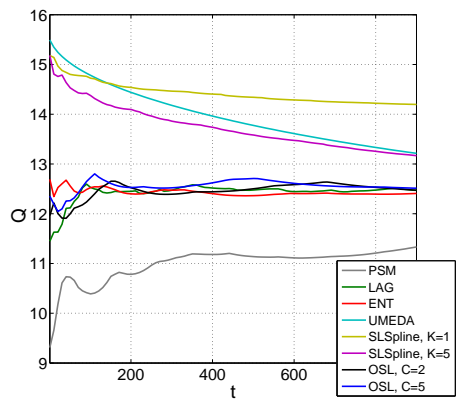}
 \caption{Quality factor for all the methods as a function of the time for the different methods applied to a step with 1D constant advection on a periodic domain.}
 \label{fig:qfactor}
 \end{figure}
The quality factor(\fref{qfactor}), as defined previously in \eqref{quality} as the ratio of $L^2$ norm over $TV$ norm, shows three groups of curves:
\begin{itemize}
\item  On the lower part, the PSM scheme curve shows that even conserving well the $L^2$ norm, the $TV$ norm is much higher than the other schemes because of spurious oscillations. The quality factor is then poor for PSM.
\item The group of four curves at the middle (LAG, PSM-ENT, PSM OSM with C=2 and C=5) are conserving the $L^2$ norm with the same order (best PSM-OSL C=2, worst LAG) and with almost the same level of oscillations according to the $TV$ norm. The quality factor is thus equivalent for this group of schemes.
\item On the upper part, LAG-UMEDA and PSM-SLS K=1 (minmod limiter) are very diffusive according to the $L^2$ norm, so they kill all the  oscillations and the $TV$ norm is very weak. The quality factor is high because these diffusive schemes have too low TV norms and cannot properly be compared with the others. However, the PSM-SLS K=5 scheme conserves the $L^2$ norm in a comparable way with the other schemes, but with a very low $TV$ norm. The quality factor is thus good for this scheme and can be compared with to the other schemes.
\end{itemize}

\paragraph{}
As a conclusion for this comparison of all schemes on a linear advection of a step function test case, we could say that the LAG-UMEDA and PSM-SLS K=1 (minmod limiter) schemes are too diffusive to be compared with the other schemes considering the quantitative quality factor \eqref{quality}. However, this way of quantitative comparison seems to fit well for this test case with the qualitative comparison or the "feeling" we might have by observing the plotted results in Figure \ref{fig:limpsmall}. 

\newpage

\subsection{4D Drift-kinetic  model}
In this section, we evaluate the limiter capacities on a 4D drift-kinetic model simulation of instabilities with the Gysela code \cite{grandgirard}. The benchmark providing instabilities is described in \cite{braeunig}. We simulate this test case with two mesh resolutions: a low resolution to provide results about all the limiters and a high resolution to get refined results with the best schemes to test also robustness. 
However, let us first do a quick review of the 4D drift-kinetic model and the time scheme used for high resolution simulations.\\

We recall hereafter the model described in \cite{gysela}. The geometrical assumptions of this model for ion plasma turbulence are a cylindrical geometry with coordinates $(r,\theta,z,v_{\|})$ and a constant magnetic field $B=B_z ~e_z$,  where $e_z$ is the unit vector in $z$ direction. In this collisionless plasma, the trajectories are governed by the guiding center (GC) trajectories:
\begin{equation}
\begin{array}{lll}
\dfrac{dr}{dt}=v_{GC_r}; ~r \dfrac{d \theta}{dt}=v_{GC_\theta}; ~ \dfrac{d z}{dt}=v_\|; ~\dfrac{d v_\|}{dt}=\dfrac{q_i}{m_i} E_z
\end{array}
\end{equation}
with $v_{GC}=(E \times B)/B^2$ and $E=-\nabla \Phi$ with $\Phi$ the electric potential.\\
The Vlasov equation governing this system, where the ion distribution function is $f(r,\theta,z,v_\|,t)$, is the following:
\begin{equation}\label{Vlasov2}
\begin{array}{lll}
\partial_t f+v_{GC_r} \partial_r f + v_{GC_\theta} \partial_\theta f  + v_\| \partial_z f  +\dfrac{q_i}{m_i} E_z \partial_{v_\|} f =0.
\end{array}
\end{equation}
This equation is coupled with a quasi-neutrality equation for the electric potential $\Phi(r,\theta,z)$ that reads:
\begin{equation}
\begin{array}{lll}
-\nabla_\perp \Phi \cdot \left( \dfrac{n_0(r)}{B~ \Omega_0} \nabla \Phi \right) + \dfrac{e~n_0(r)}{T_e(r)} (\Phi-<\Phi>_{z})=n_i-n_0
\end{array}
\end{equation}
with $n_i=\ds \int_{v_\|} f(r,\theta,z,v_\|) d v_\|$ and constant in time physical parameters $n_0$, $ \Omega_0$, $T_e$ and $e$.\\
Let us notice that the 4D velocity field $a=(v_{GC_r},v_{GC_\theta},v_\|, q/m_i ~E_z)^t$ is divergence free:
\eqs
\begin{array}{lll}
\nabla \cdot a = \dfrac{1}{r} \partial_r  (r~v_{GC_r})+  \dfrac{1}{r} \partial_\theta  (v_{GC_\theta})+\partial_z  v_\| + \partial_{v_\|} (q/m_i ~E_z)=0
\end{array}
\eeqs
because of variable independence $\partial_{v_\|}  E_z=\partial_{v_\|} (\partial_z \Phi(r,\theta,z))=0$, $\partial_z  v_\|=0$ and \\
$v_{GC_r}=\dfrac{-1}{r~B_z} \partial_\theta \Phi $ and  $v_{GC_\theta}/r=\dfrac{1}{r~B_z}  \partial_r \Phi,$\\
 such that 
 $$ \dfrac{1}{r} \partial_r  (r~v_{GC_r})+  \dfrac{1}{r} \partial_\theta  (v_{GC_\theta})=0.$$
Therefore, one can write an equivalent conservative equation to the preceding Vlasov equation \eqref{Vlasov2}:
\eqs
\begin{array}{lll}
\partial_t f+\partial_r (v_{GC_r} ~f)+  \partial_\theta (v_{GC_\theta} ~f)  +  \partial_z (v_\| ~f)  + \partial_{v_\|} \left( \dfrac{q_i}{m_i} E_z ~f \right) =0.
\end{array}
\eeqs

\subsection{Low resolution simulations}\label{low}

\begin{figure}
 \centering
 \begin{tabular}{c}
  \includegraphics*[scale=0.5]{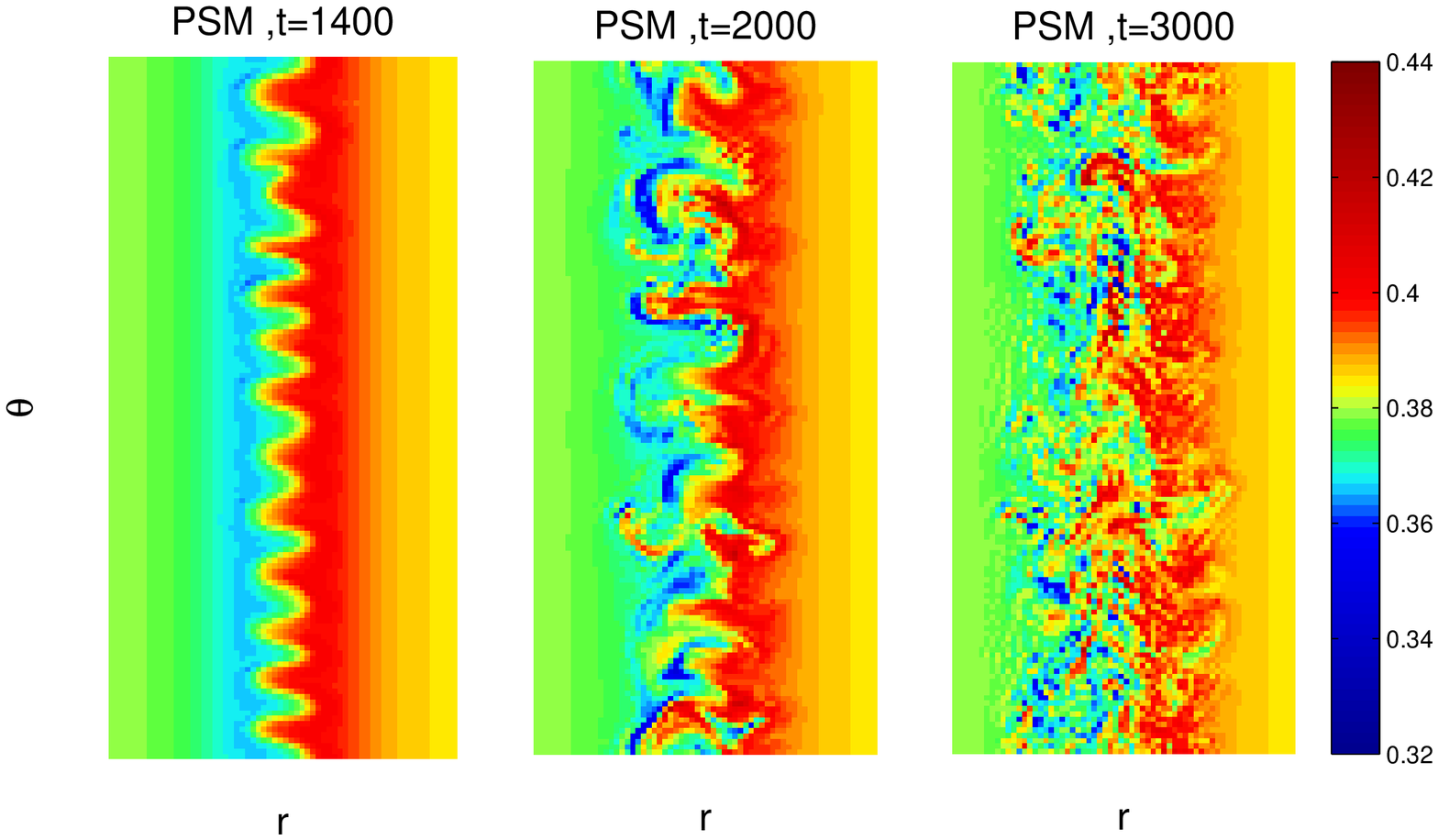} \\
   \includegraphics*[scale=0.5]{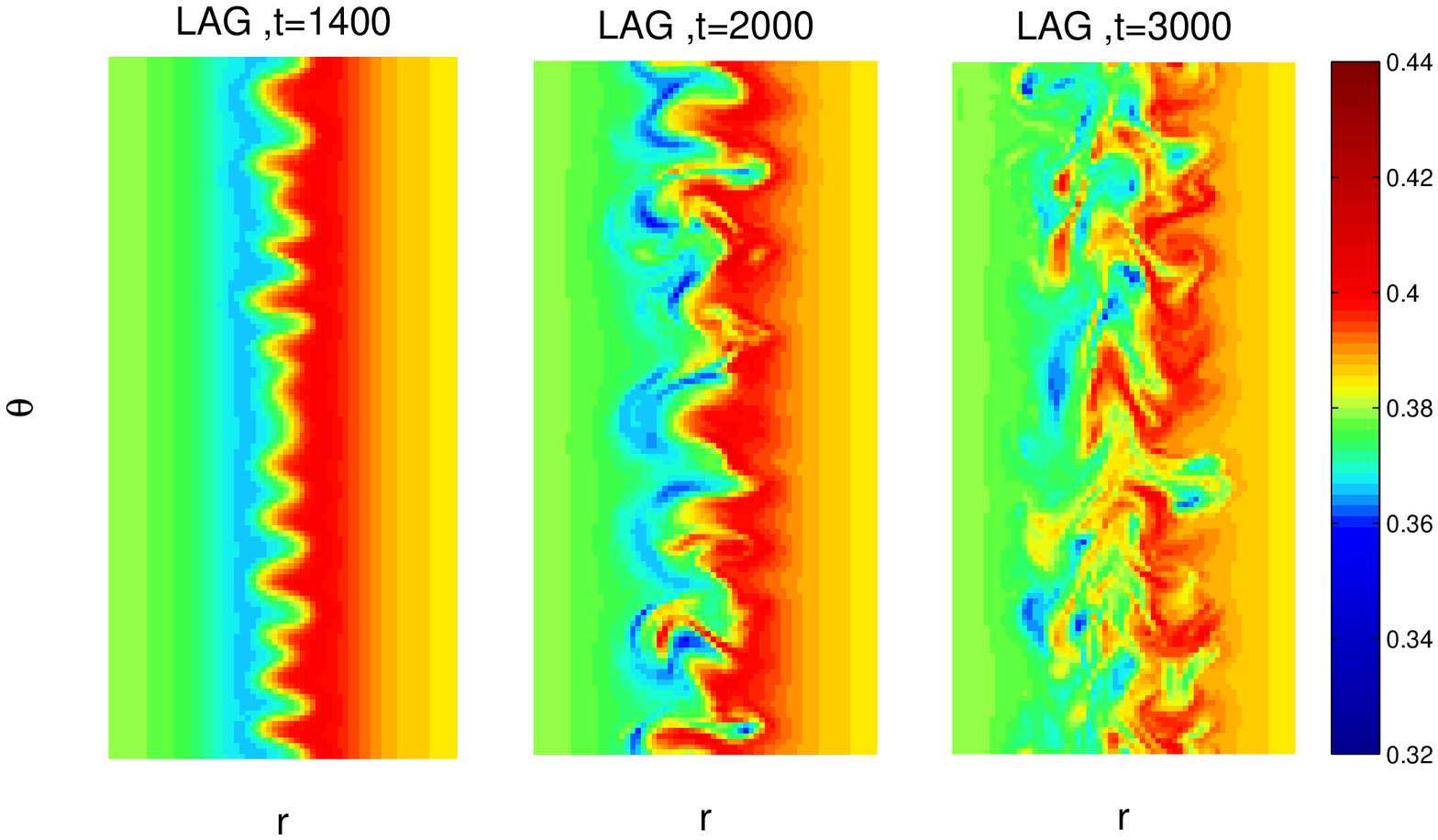} \\
   \includegraphics*[scale=0.5]{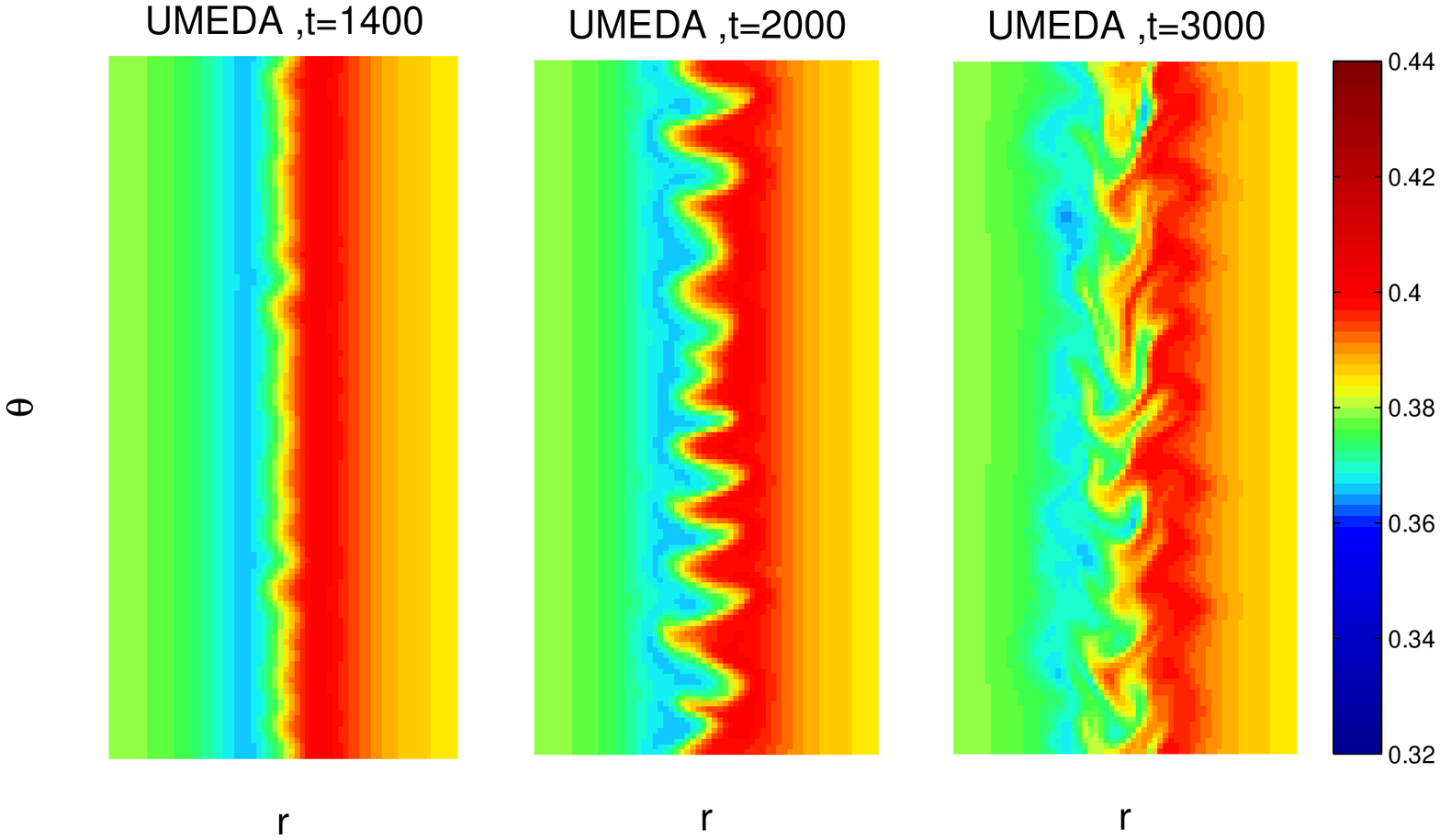} 
 \end{tabular}
 \caption{\label{fig:4Dlag} Solution in a $(r,\theta)$ plane with $N_r \times N_\theta=64 \times 128$ cells (low resolution) of a 4D test case with PSM, LAG, LAG-UMEDA  methods. }
 \end{figure}

\begin{figure}
 \centering
 \begin{tabular}{c}
   \includegraphics*[scale=0.5]{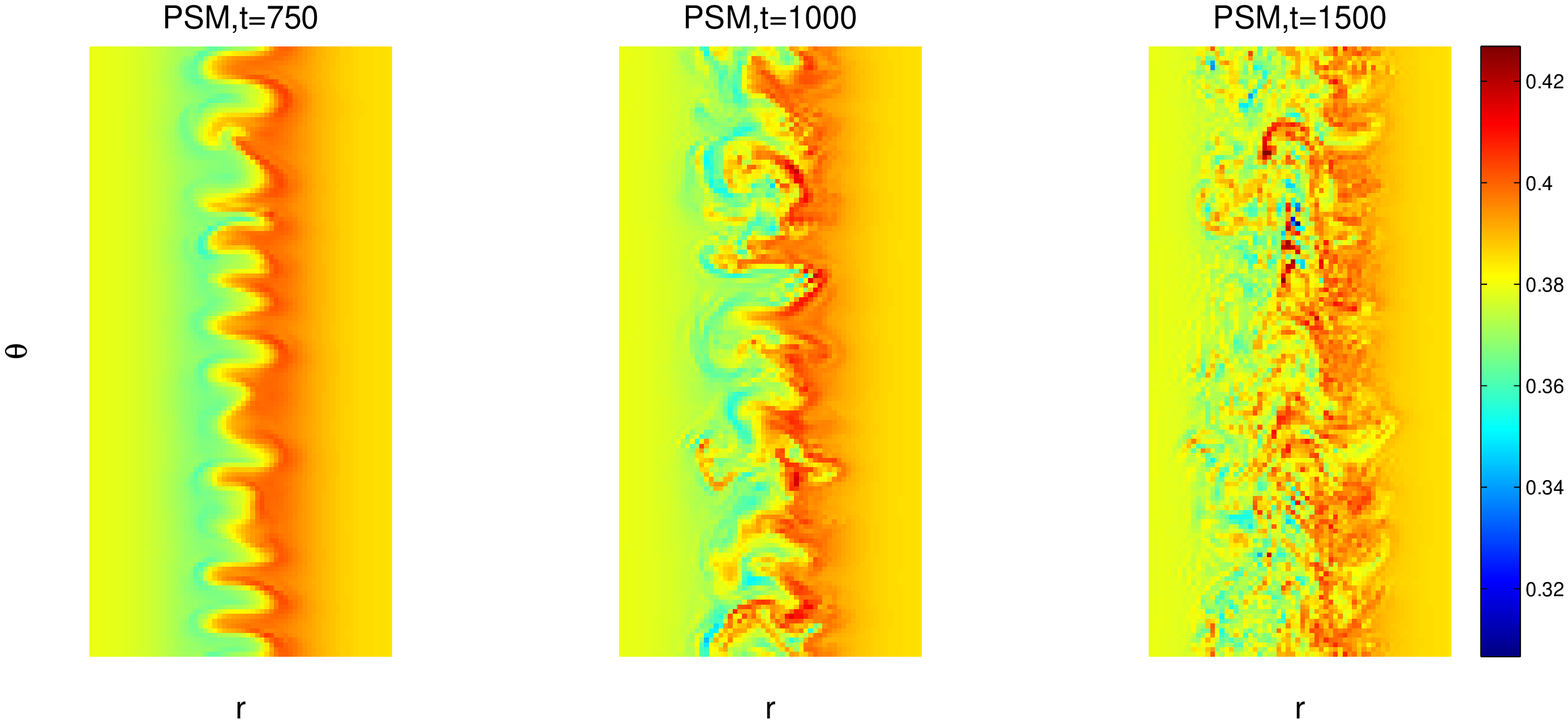} \\
   \includegraphics*[scale=0.5]{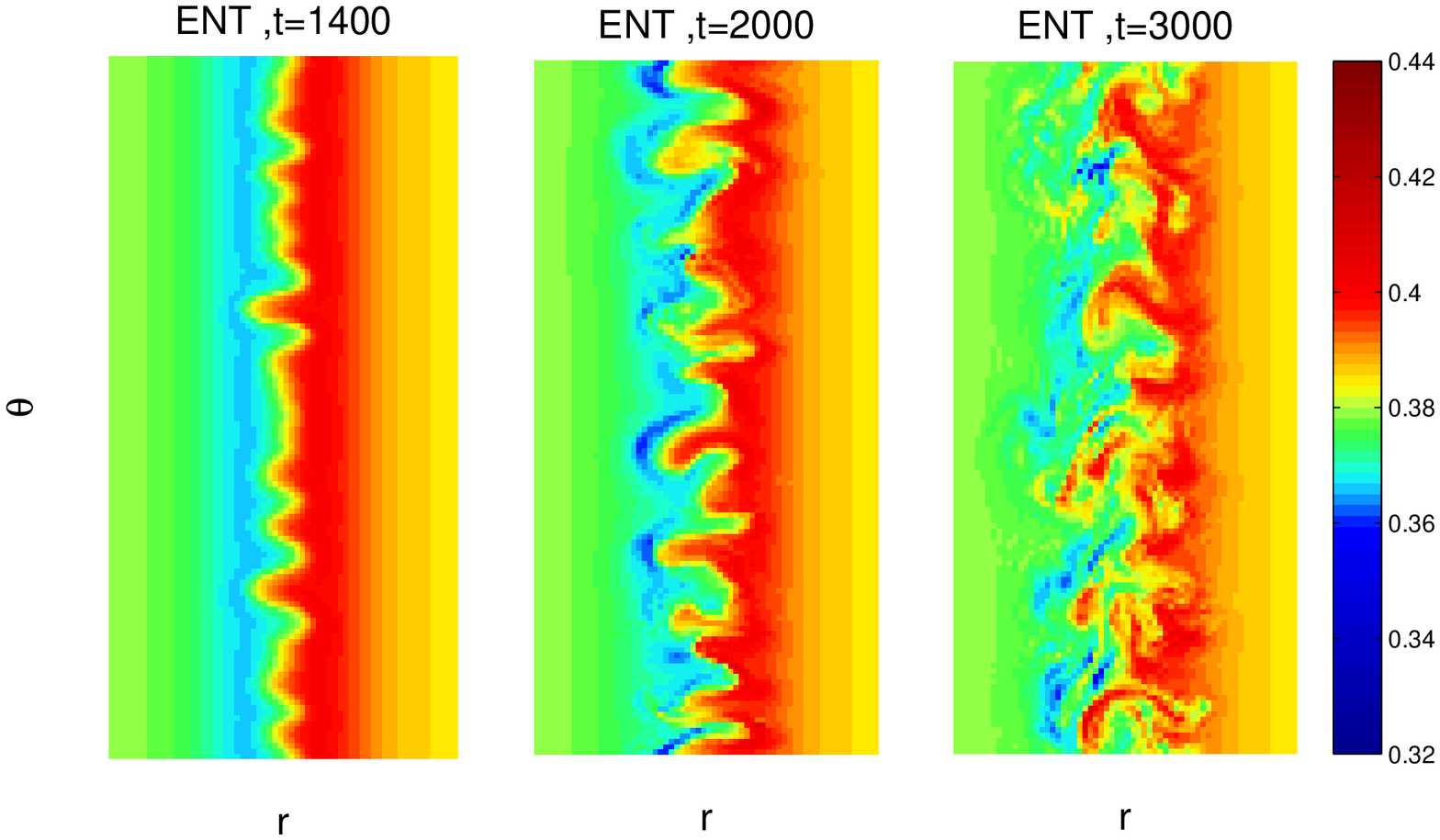} \\
    \includegraphics*[scale=0.5]{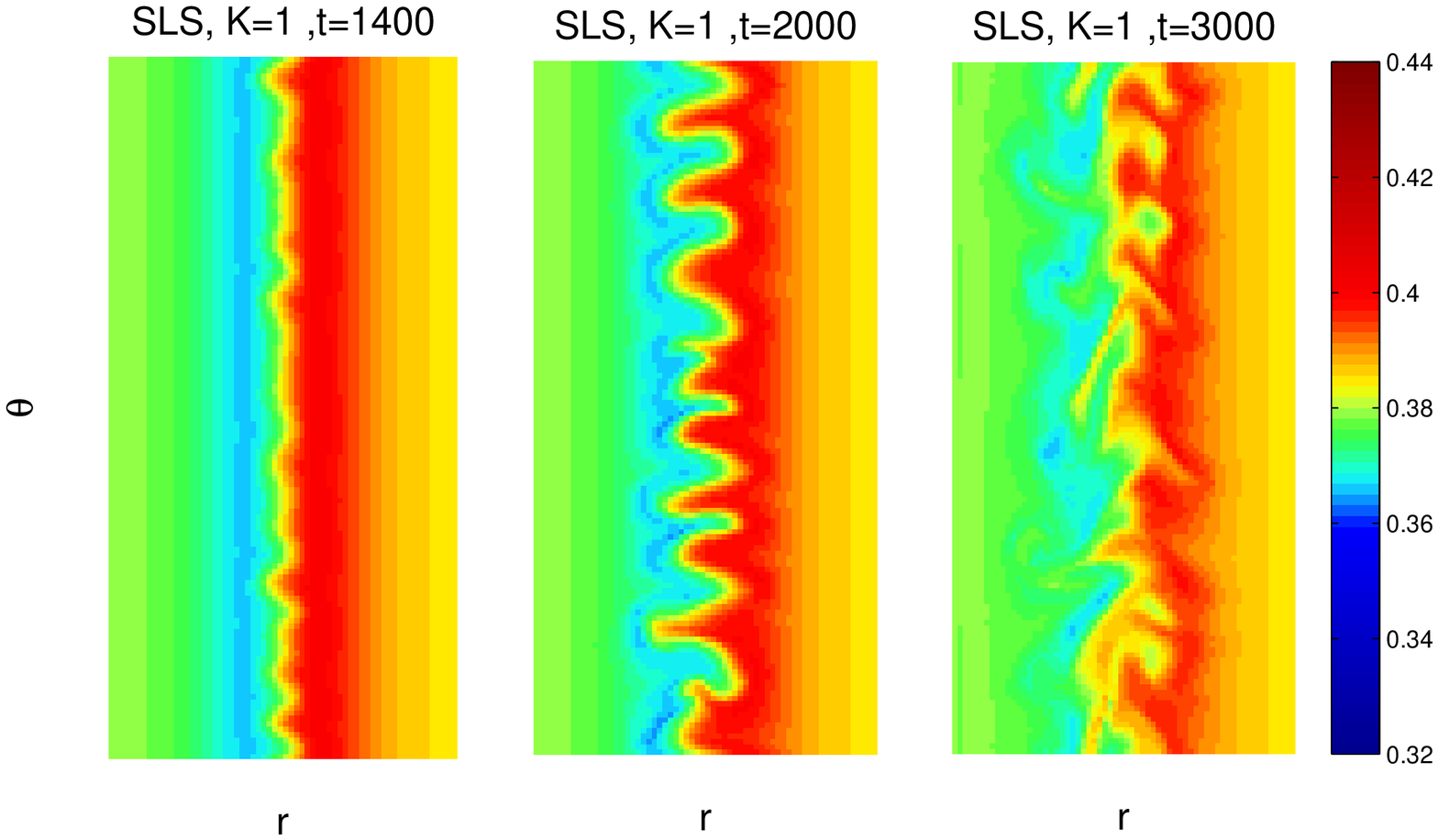} 
 \end{tabular}
 \caption{\label{fig:4Dent} Solution in a $(r,\theta)$ plane with $N_r \times N_\theta=64 \times 128$ cells (low resolution) of a 4D test case with PSM, PSM-ENTand PSM K=1 methods. }
 \end{figure}

\begin{figure}
 \centering
 \begin{tabular}{c}
   \includegraphics*[scale=0.5]{profil/PSM.eps} \\
   \includegraphics*[scale=0.5]{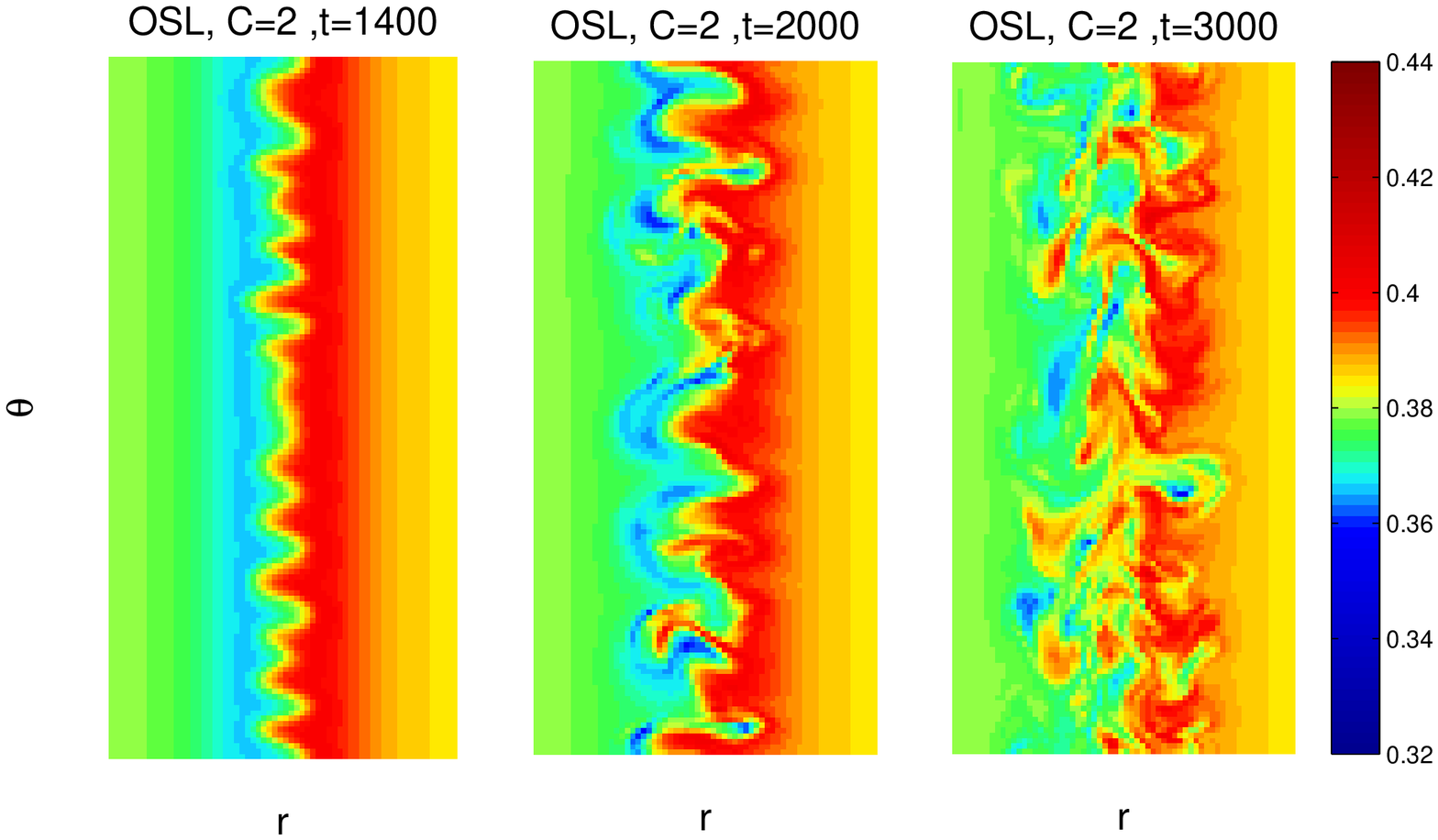} \\
   \includegraphics*[scale=0.5]{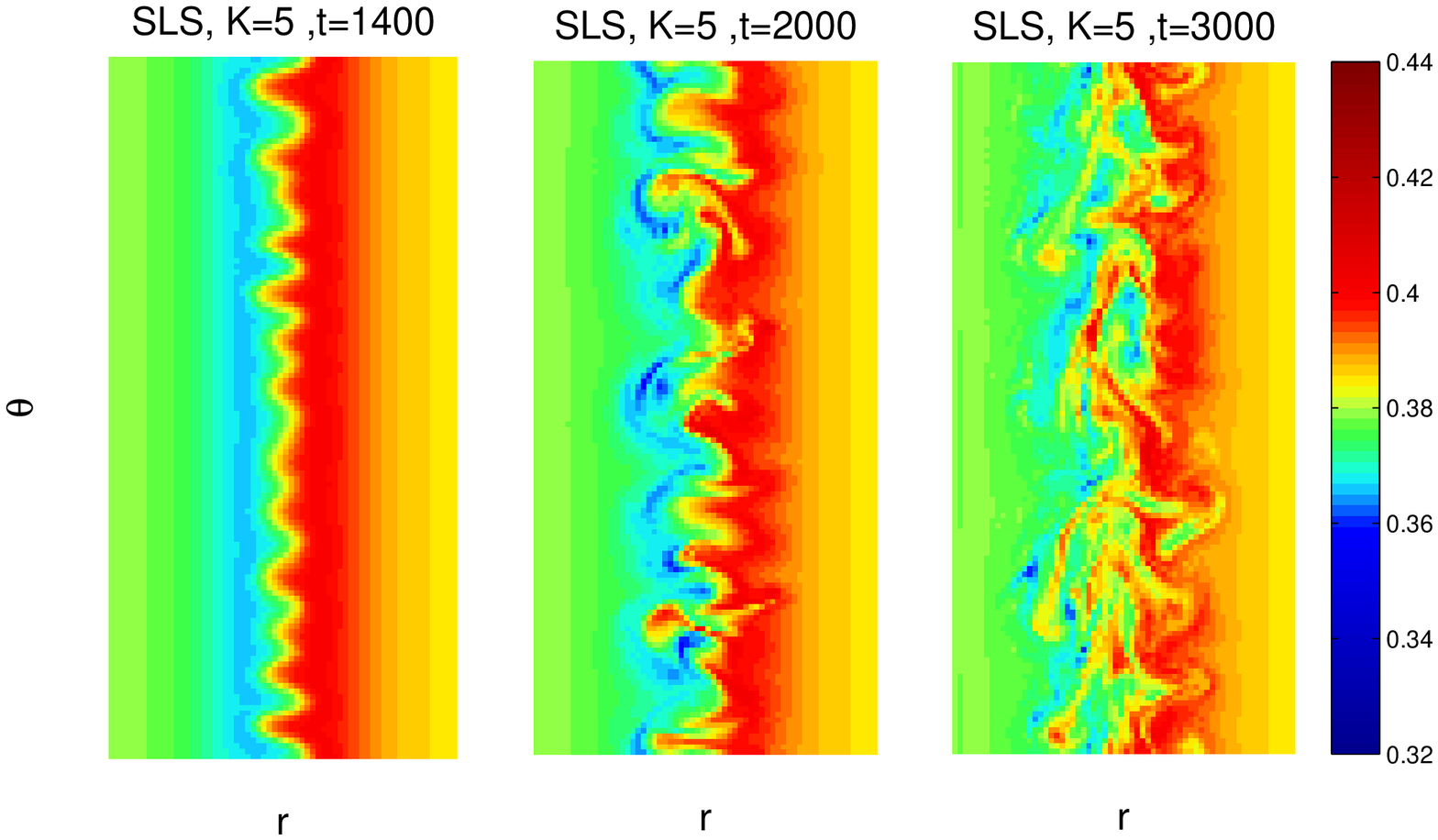} 
 \end{tabular}
 \caption{\label{fig:4Doslsls} Solution in a $(r,\theta)$ plane with $N_r \times N_\theta=64 \times 128$ cells (low resolution) of a 4D test case with PSM, PSM-OSL C=2 and PSM-SLS K=5 methods. }
 \end{figure}

 The low-resolution simulations run on a 64x128x16x16 grid ($N_{r} \times N_{\theta} \times N_{\phi} \times  N_{v_{\parallel}}$). The resolution is too low to consider the total energy as a relevant measure. The entropy, $L^2$ norm and $TV$ norm are used to gauge the limiter's effects. We present in Figures  
\ref{fig:4Dlag} \ref{fig:4Dent}, \ref{fig:4Doslsls} three states of the instabilities: the linear phase (left), the beginning of the non-linear phase (center) and the strong non-linear phase (right).  \\
 
In figure \ref{fig:4Dlag}, we see the difference of behaviour of the PSM scheme based on a spline reconstruction method and the LAG scheme based on Lagrangian polynoms. The LAG scheme diffuses much more the solution even at early times, and this is accentuated using the LAG-UMEDA scheme which makes disappear  even the big structures. However, the PSM scheme solution shows a lot of small structures, actually at the scale of the mesh, which are suspected of coming from spurious oscillations already observed on the 1D step test case. Moreover, these oscillations often violate the theoretical maximum principle of the solution and may lead to crash the simulation.

\paragraph{}
In figure \ref{fig:4Dent}, we see the effect of the ENT limiter on the PSM scheme. It reduces a lot the oscillations with a figure showing much less small structures, but showing a more diffusive behaviour. The PSM-SLS K=1 scheme is extremely diffusive and makes disappear  even the big structures. 
\paragraph{}
In figure \ref{fig:4Doslsls}, we compare the PSM scheme result with PSM-OSL and PSM-SLS limiters results. Both limiters reduces efficiently the oscillations. The limiter PSM-SLS seems a little more diffusive than PSM-OSL.
 
\newpage
\paragraph{}
 \begin{figure}
 \centering
 \begin{tabular}{cc}
  \includegraphics[scale=0.3]{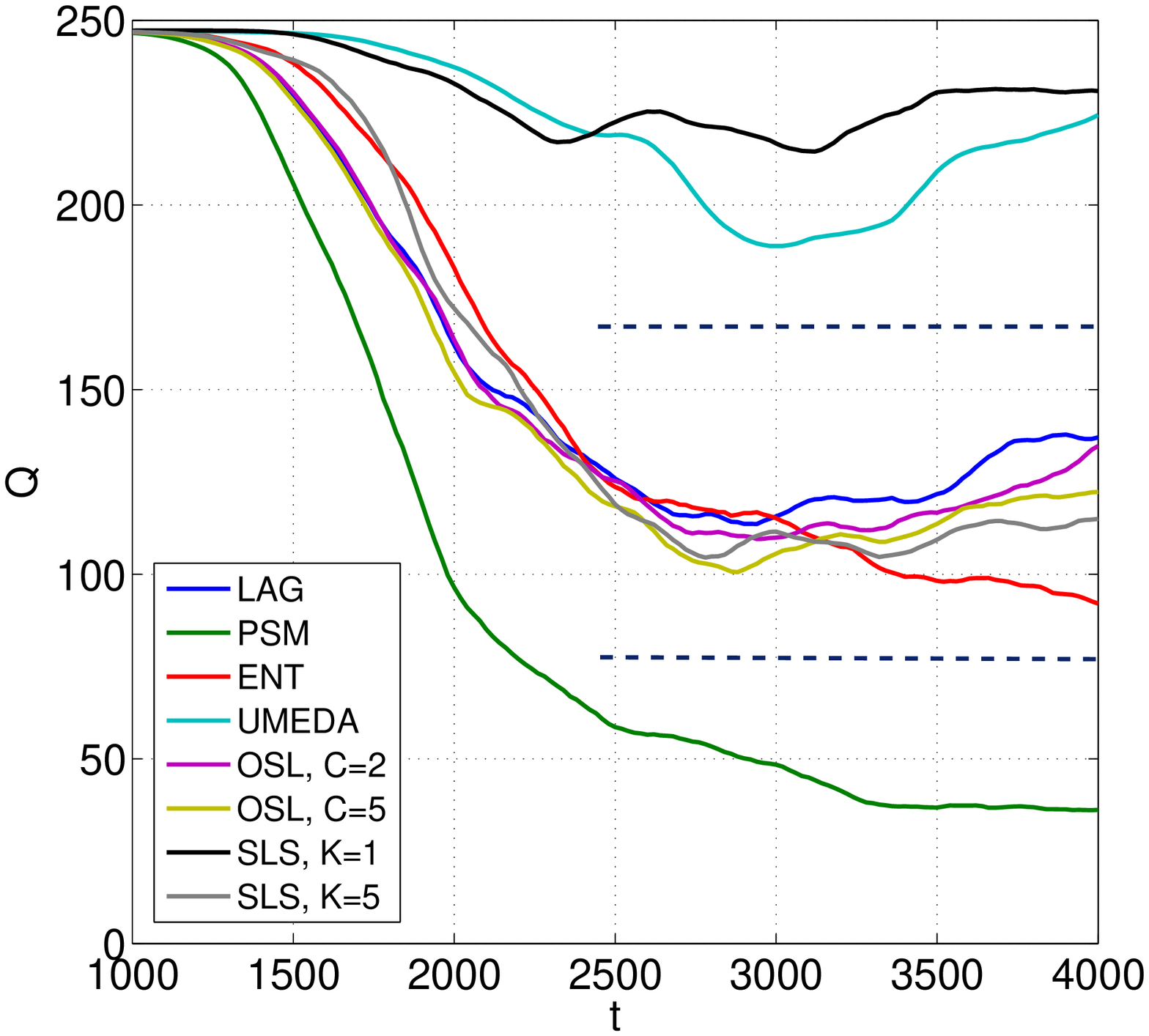}           & \includegraphics[scale=0.3]{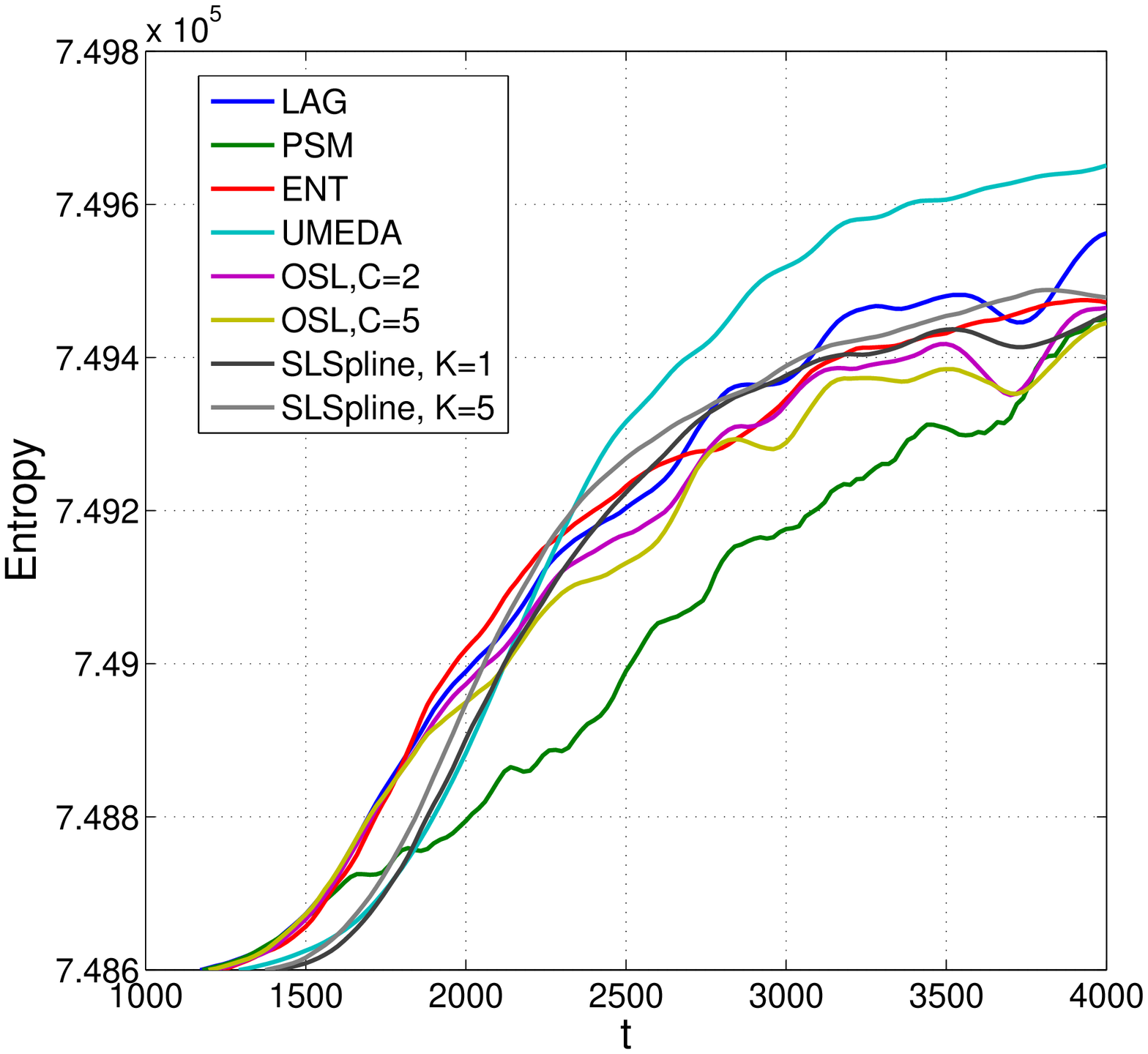}\\
  \includegraphics[scale=0.3]{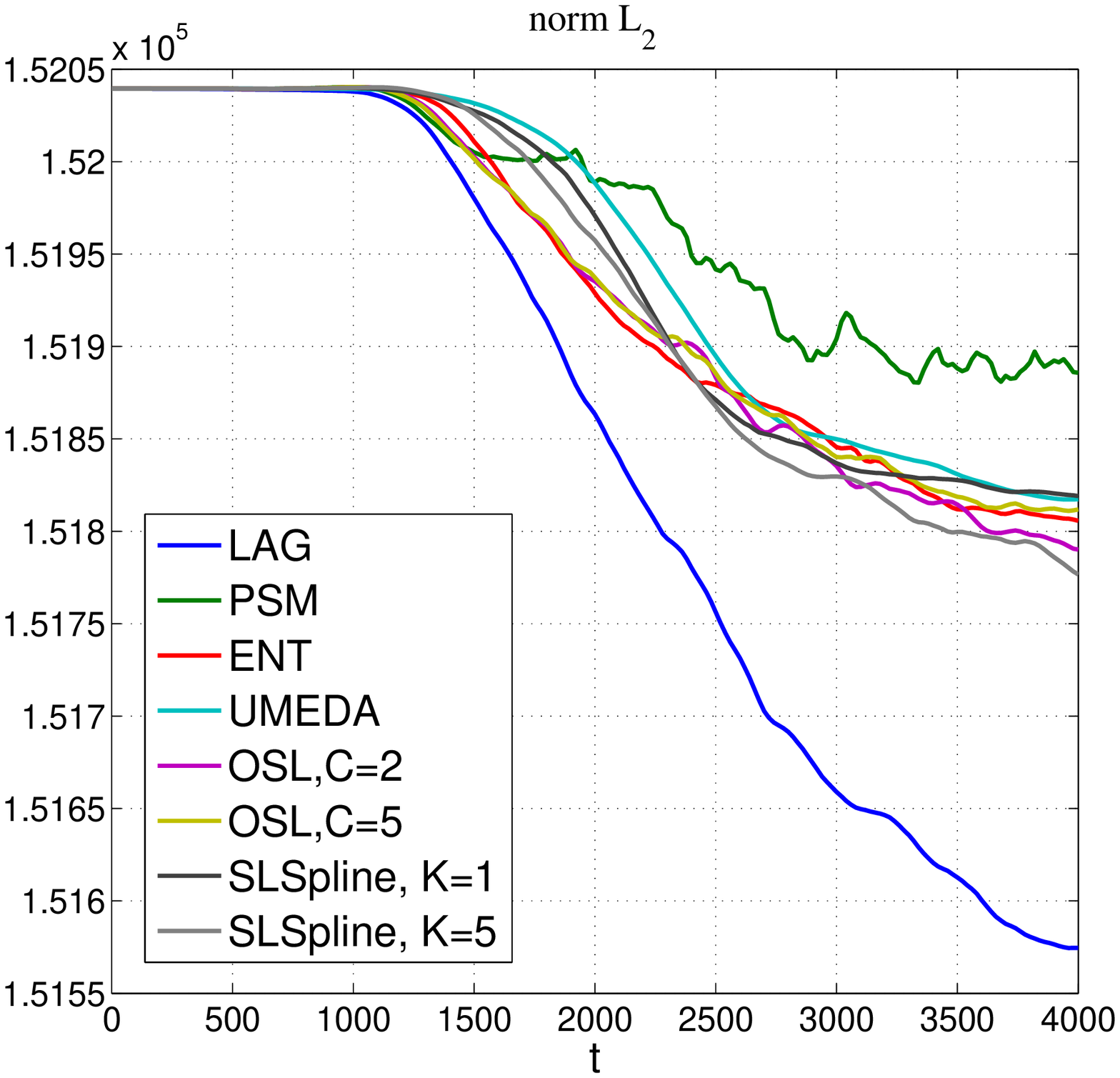}      & \includegraphics[scale=0.3]{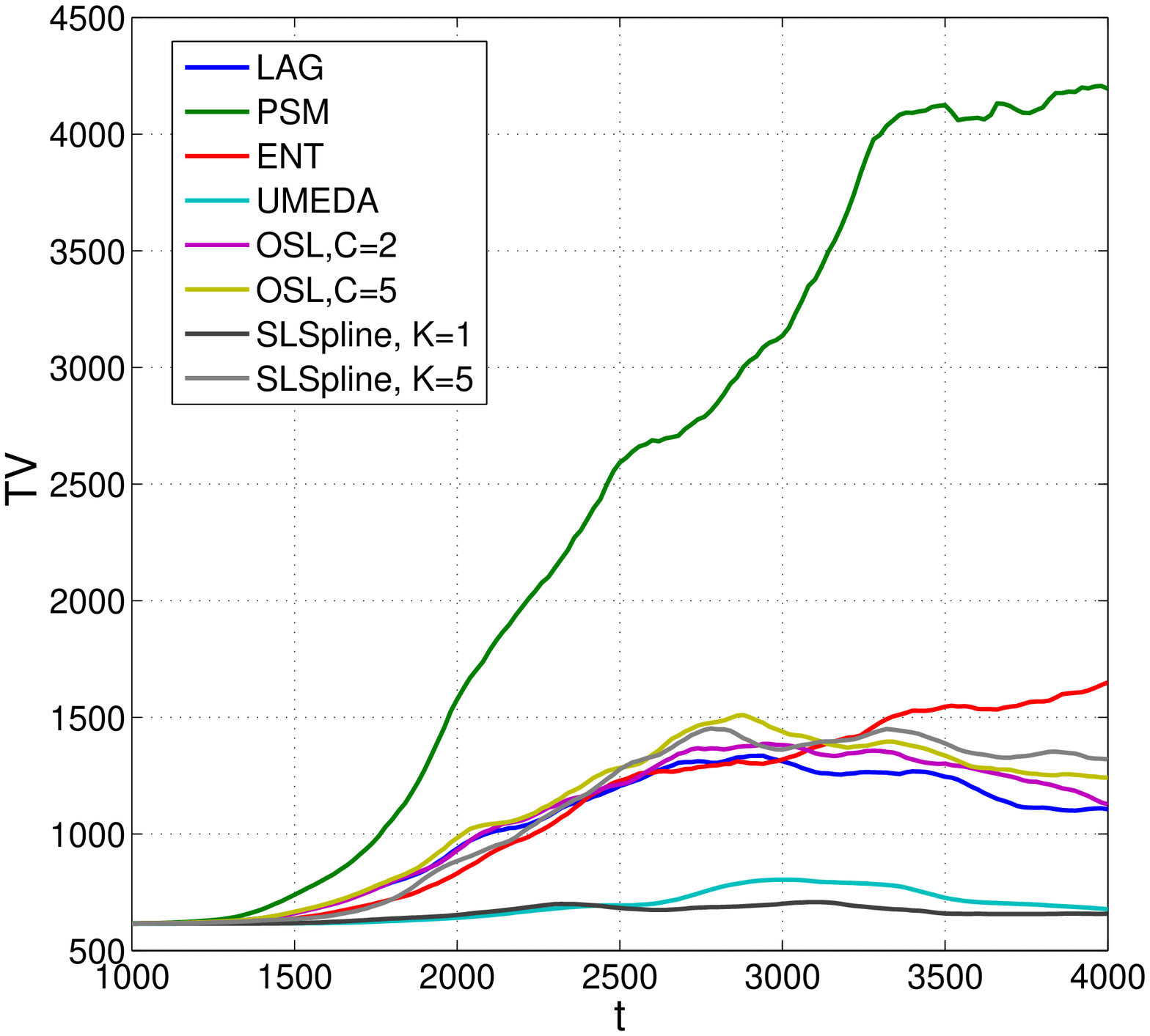}
 \end{tabular}
 \caption{Quality factor  (left top), entropy (right top), $L^2$ norm (left bottom), TV norm (right bottom) for the 4D simulation low resolution $64 \times 128  \times 16 \times 16$  with all the methods presented.}
 \label{fig:4Dresultsall}
 \end{figure}
In figure  \ref{fig:4Dresultsall} we see the effect of each limiter on the quality factor, entropy, $L^2$ norm and TV norm. 
\begin{itemize}
\item The LAG scheme gives very poor results in the conservation of the $L^2$ norm, but shows a medium quality factor.
\item The LAG-UMEDA and PSM-SLS K=1 schemes are still diminishing the TV norm at a very low rate as for the 1D step test case, because they are very diffusive. Even if the quality factor is good and the $L^2$ norm close to the other schemes, the qualitative results of figure \ref{fig:4Dlag}  and \ref{fig:4Dent} show that even big structures are not well captured with these schemes.
\item The PSM scheme without limiter has still a very different behaviour compare with the others : the $L^2$ norm and entropy plots show that  this scheme is less diffusive and the TV norm plots confirms the qualitative 2D profiles \ref{fig:4Doslsls}  that it is much more oscillating than the other schemes. Thus the quality factor is lower than the other schemes. 
\item The PSM scheme with limiters ENT, OSL, and SLS K=5 show equivalent medium quality factors and TV norms and conservation of entropy and $L^2$ norm for this coarse benchmark. 
\end{itemize}

\paragraph{Conclusion of the 4D results with low resolution section\\}
The results obtained for this 4D test cas low resolution confirms those obtained for the 1D linear step test case. The PSM scheme shows oscillations as the TV norm is high compare to the other schemes, but preserves the better $L^2$ norms and the entropy. The PSM-SLS K=1 and  LAG-UMEDA schemes are diffusive and seem experimentally of lower order of accuracy than the other schemes. The  PSM-SLS K=1 limiter, PSM-OSL limiter and ENT limiter are showing a similar level of accuracy. The high resolution simulation results following will permit to go further in these investigations.

\newpage
\subsection{High resolution simulation}
\label{sec:HR}
We propose in this section to compare PSM, OSL and SLS with an high resolution mesh relatively to the previous resolution. The mesh is constituted with $N_{r} \times N_{\theta} \times N_{\phi} \times  N_{v_{\parallel}}=256\times 512 \times 32 \times 16$. We consider that this mesh is refined enough such that conservation of the total energy becomes relevant to estimate the limiter quality, what was not the case for the low resolution mesh used in the preceding section \ref{low}. \\
We have chosen standard values for PSM limiters according to the preceding results : C=2 for the OSL limiter which depends weakly on this value and K=5 fot the SLS limiter which seems to be the minimum value to limit diffusion (see results with K=1 in section \ref{low}).  \\
The results are presented in three ways : full $(r,\theta)$ cut planes,  zoom in smaller boxes in these cut planes to better see the influence of the limiters on small structures and 1D diagnostics with the total energy, the quality factor \eqref{quality}, the entropy, the $L^2$ and TV norms.  
 \newpage
 \begin{figure}
 \centering
 \begin{tabular}{c}
  \includegraphics*[scale=0.7]{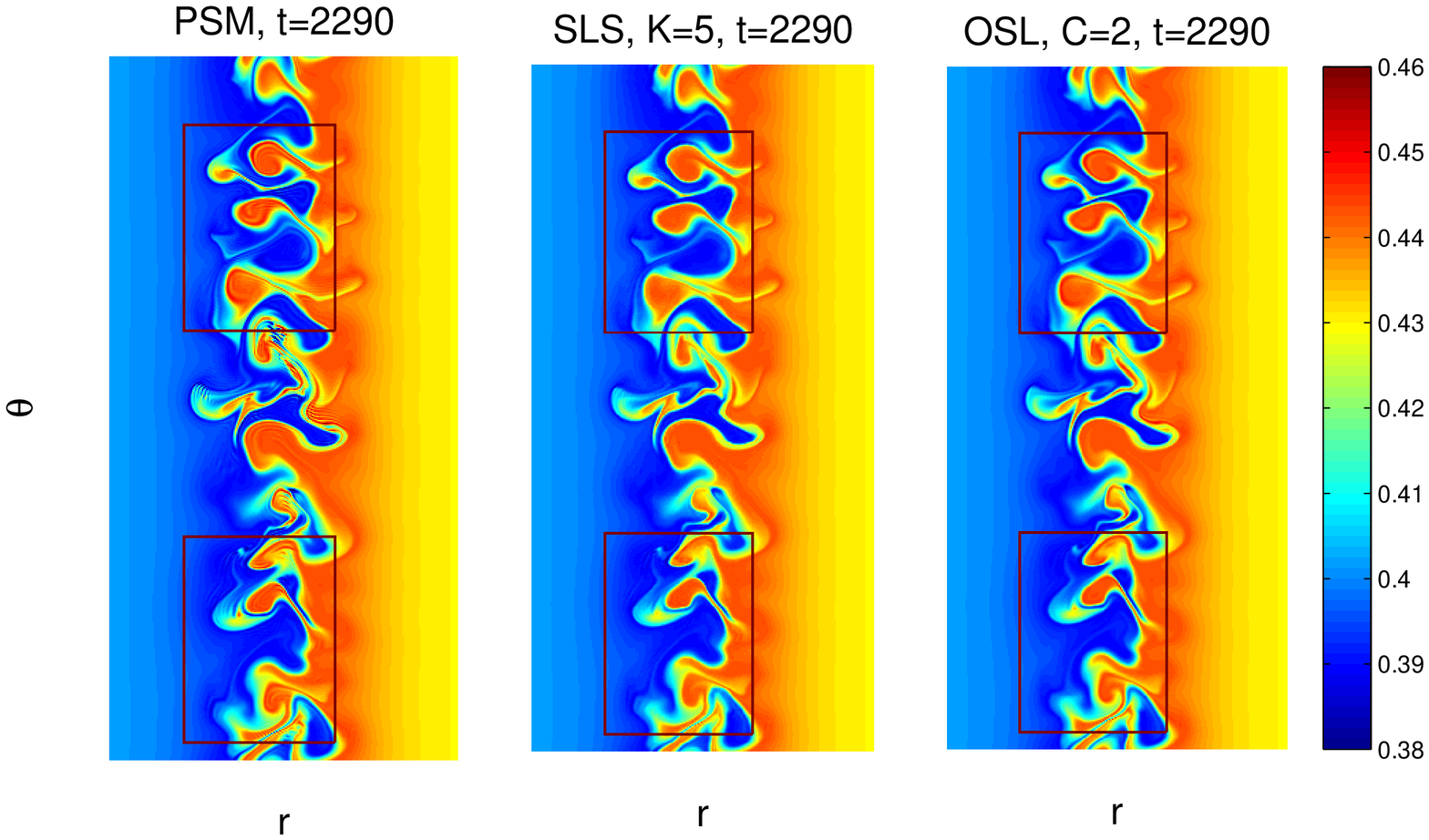} \\
  \includegraphics*[scale=0.7]{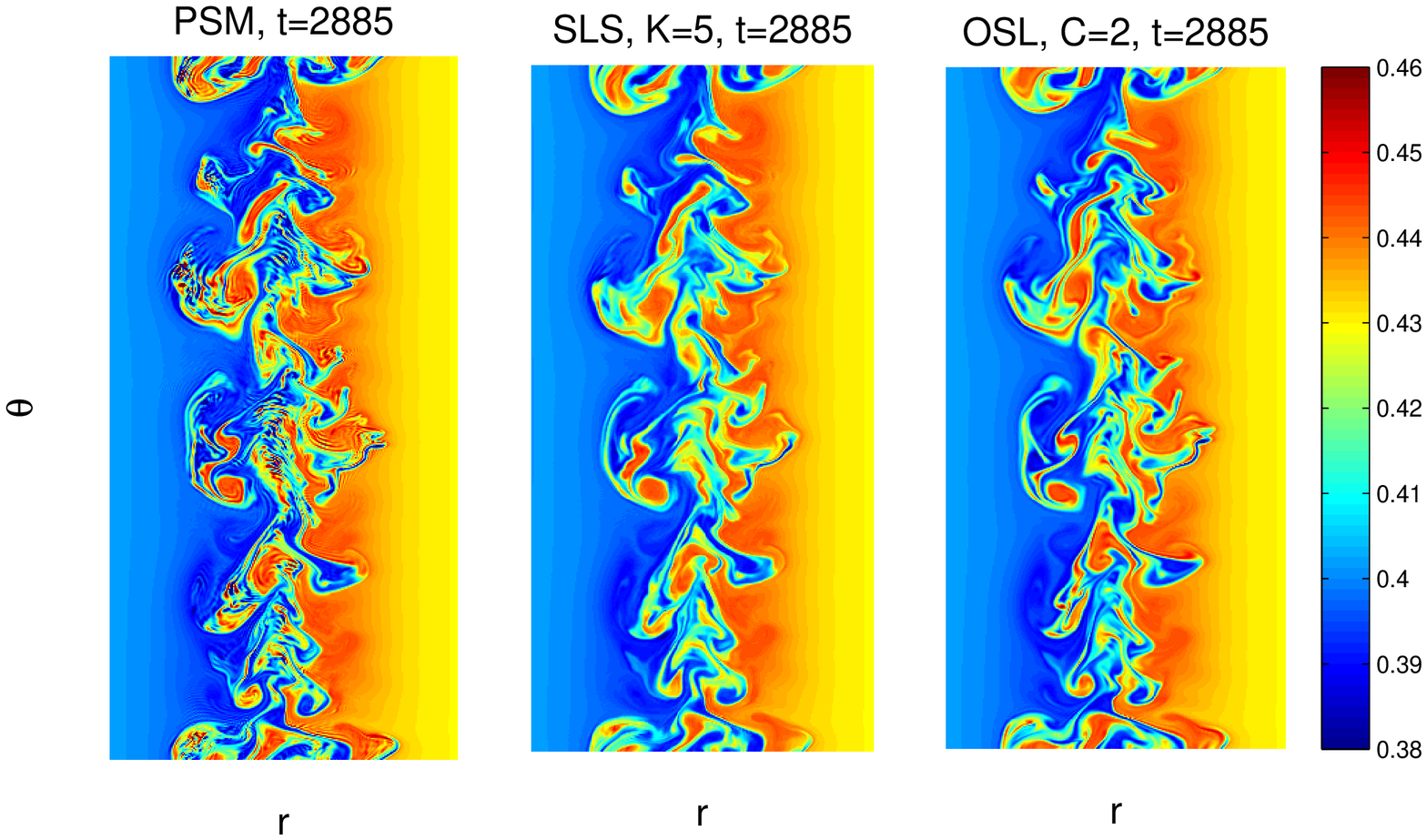} \\
 \end{tabular}
 \caption{4D simulations on $N_r \times N_\theta=256 \times 512$ cells (high resolution) along $(r,\theta)$  with PSM (left), (center) SLS, K=5 and OSL (C=2) (right) at three different times: t=2290 (top), t=2885 (bottom).}
 \label{fig:HR4D1}
 \end{figure}
 \newpage
 \begin{figure}
 \centering
 \begin{tabular}{c}
   \includegraphics*[scale=0.7]{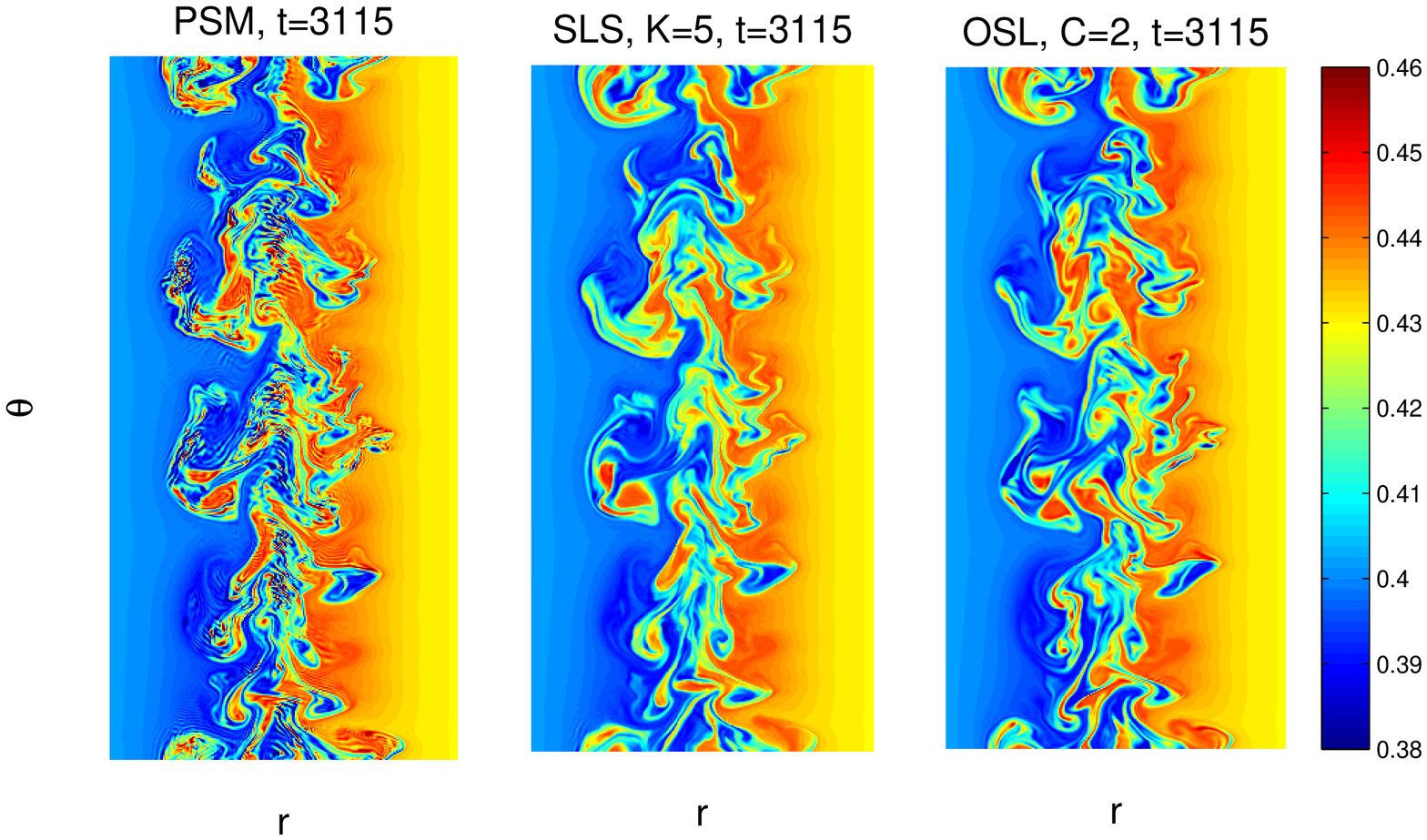} \\
   \includegraphics*[scale=0.7]{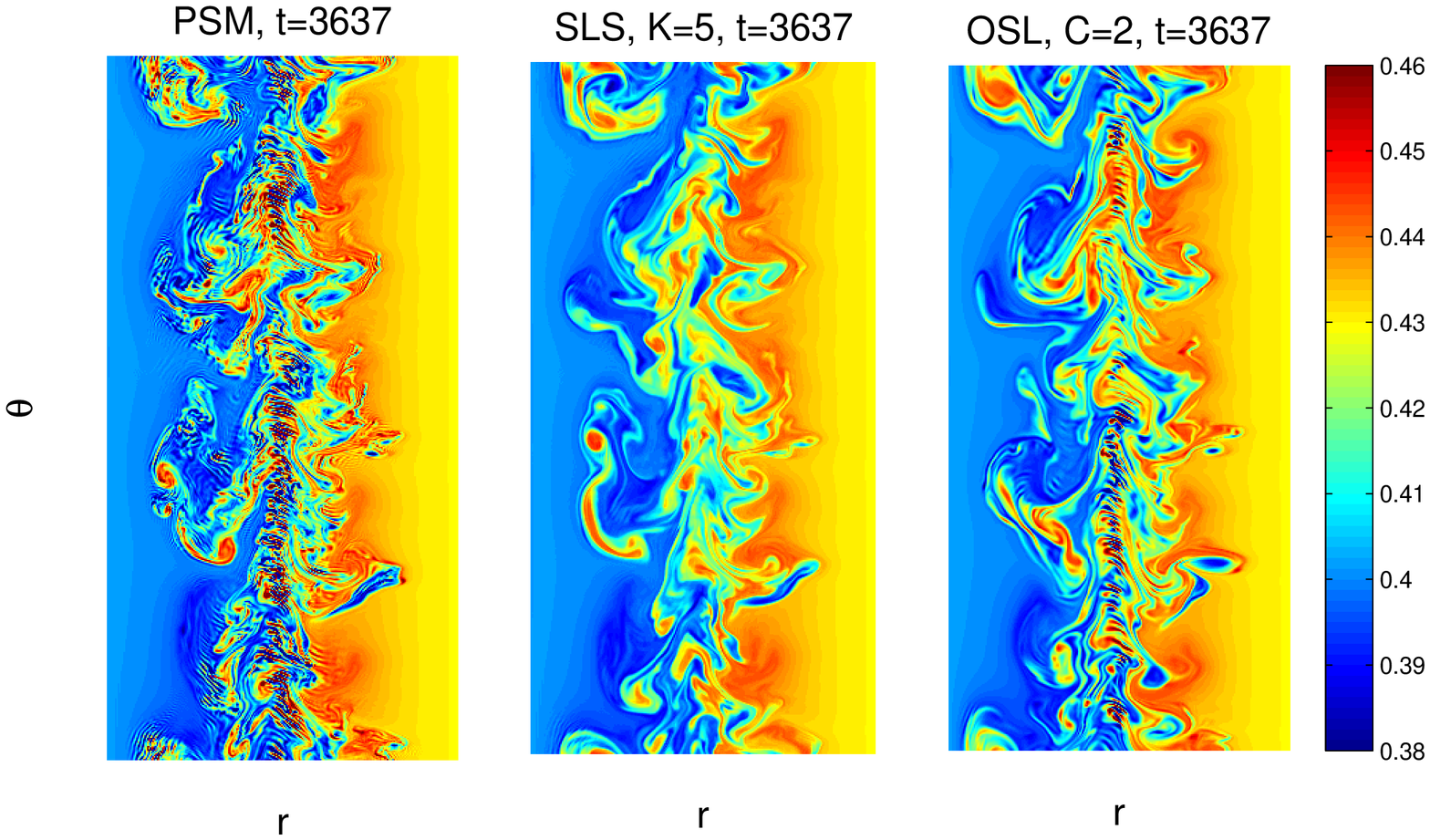} \\
 \end{tabular}
 \caption{4D simulations on $N_r \times N_\theta=256 \times 512$ cells (high resolution) along $(r,\theta)$  with PSM (left),  SLS, K=5(center) and OSL, C=2(right) at three different times: t= 3115(top),  t= 3637 (bottom). }
 \label{fig:HR4D2}
 \end{figure}
 
 \paragraph{Comparison of results in $(r,\theta)$ cut planes\\}
 In figure \ref{fig:HR4D1} and \ref{fig:HR4D2}, we present at different times $(r,\theta)$ profiles of the distribution function in the non-linear phase.  The influence of the SLS or OSL limiters is weak at the beginning of the linear phase, i.e. time $t=2290$. Afterward, a lot of small structures as filaments develop, where steep gradients exist. Therefore for later times, the differences between schemes results increase because of the different actions of SLS and OSL limiters. \\
 \begin{itemize}
\item The PSM scheme shows spurious oscillations very quickly in the non-linear phase. As seen in the 1D step test-case and in the 4D test-case with low resolution, the PSM scheme develops oscillations when transporting discontinuities or steep gradients which occurs in this flow where filaments and vortex develop.  Unfortunately, these oscillations lead to the crash of the PSM simulation at final time presented $t=3637$.
\item The PSM-SLS scheme does almost develop no oscillations even in the latest time. The price to pay is a qualitatively more diffused profile of distribution function.
\item The PSM-OSL scheme does neither develop oscillations at early times of the non-linear phase. The OSL limiter produces less diffusion than the SLS limiter. However, at late times the same kind of oscillations than the PSM scheme develops. 
 \end{itemize}
  \begin{figure}
 \centering
  \begin{tabular}{c}
   \includegraphics*[scale=0.6,viewport=40 90 600 340]{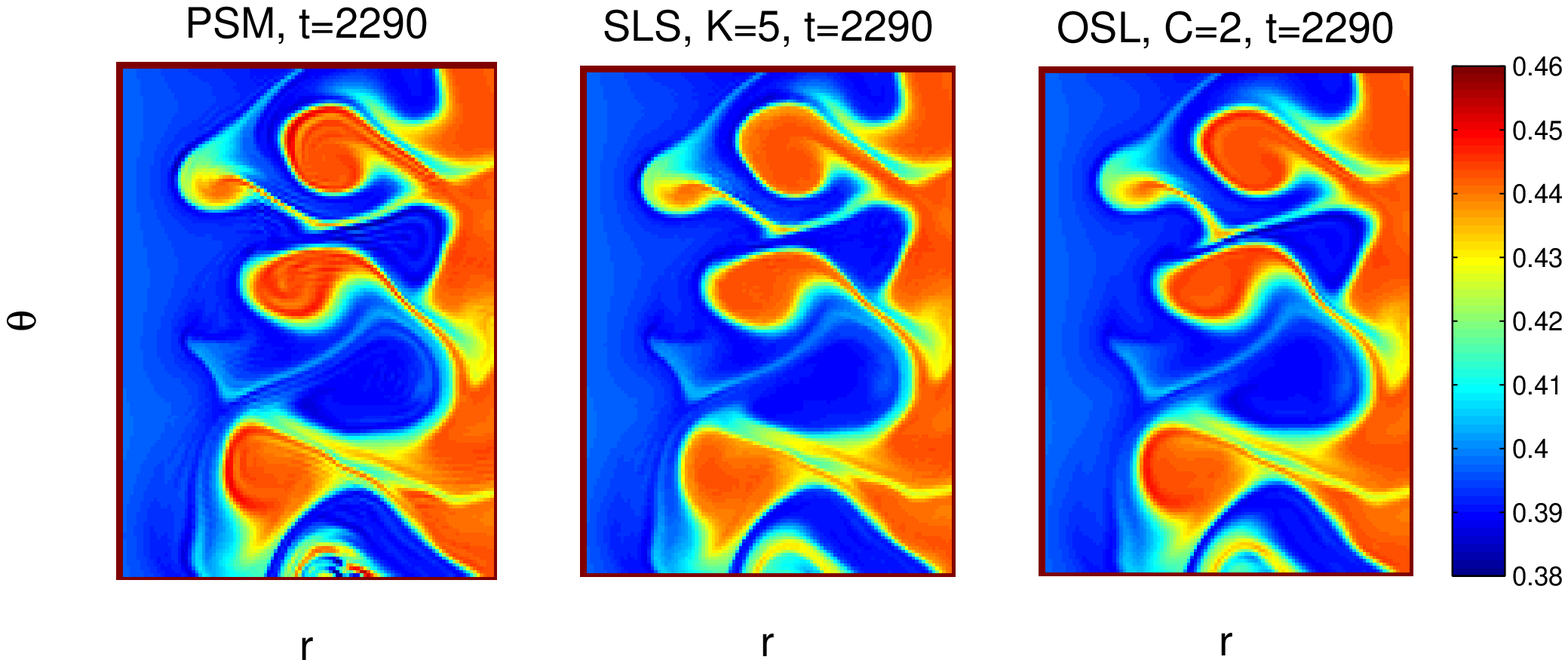}  \\
   \includegraphics*[scale=0.6,viewport=40 90 600 340]{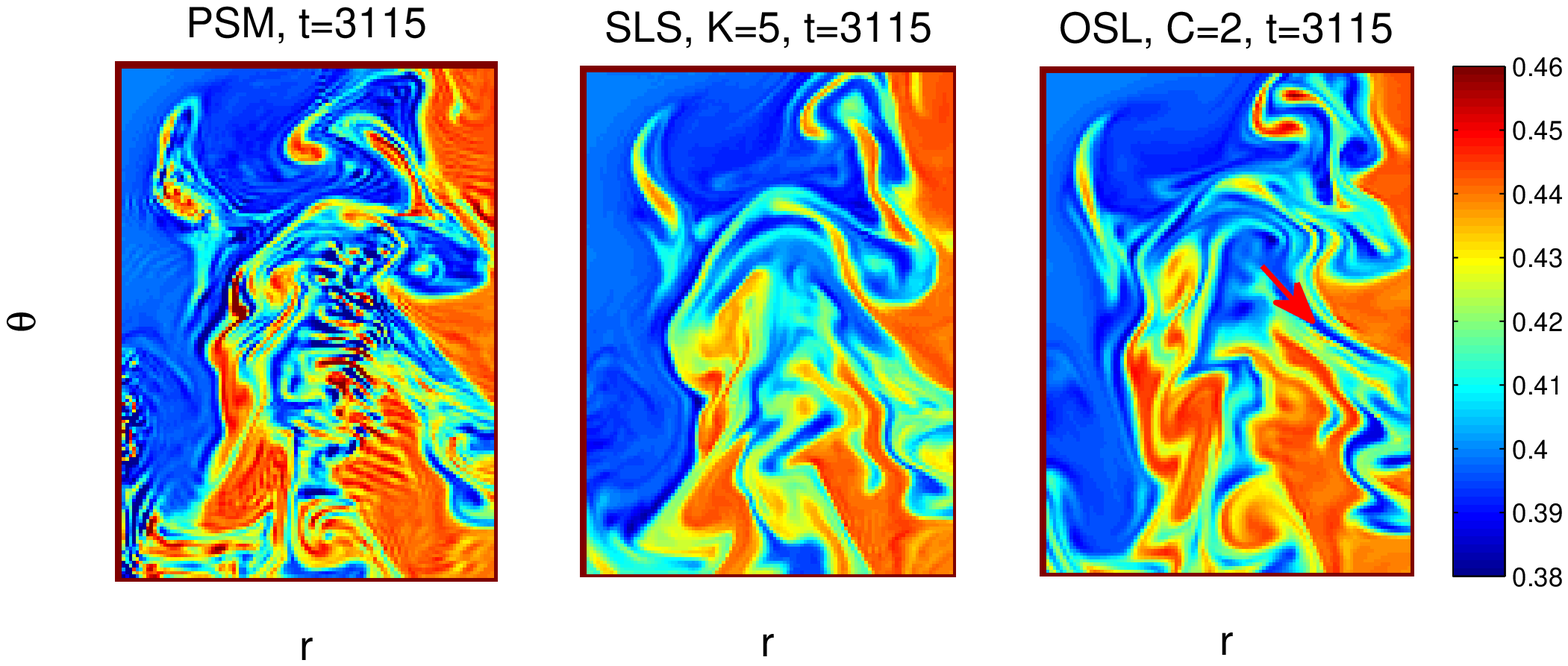} \\  
   \includegraphics*[scale=0.6,viewport=40 90 600 340]{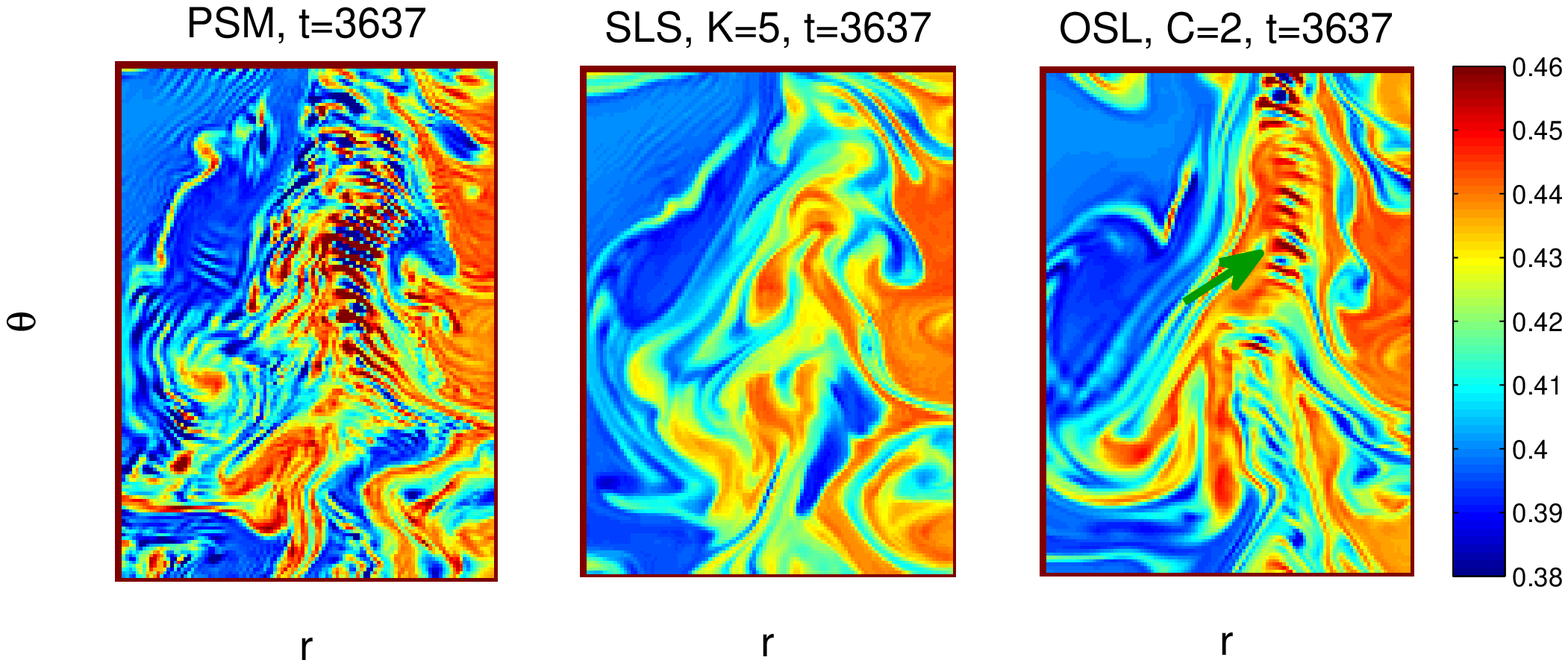} \\
 \end{tabular}
 \caption{Zoom on the zone A of the high-resolution $(r,\theta)$ profiles (\fref{HR4D1})}
 \label{fig:zoomA}
 \end{figure}
 \begin{figure}
 \centering
 \begin{tabular}{c}
   \includegraphics*[scale=0.6,viewport=40 90 600 340]{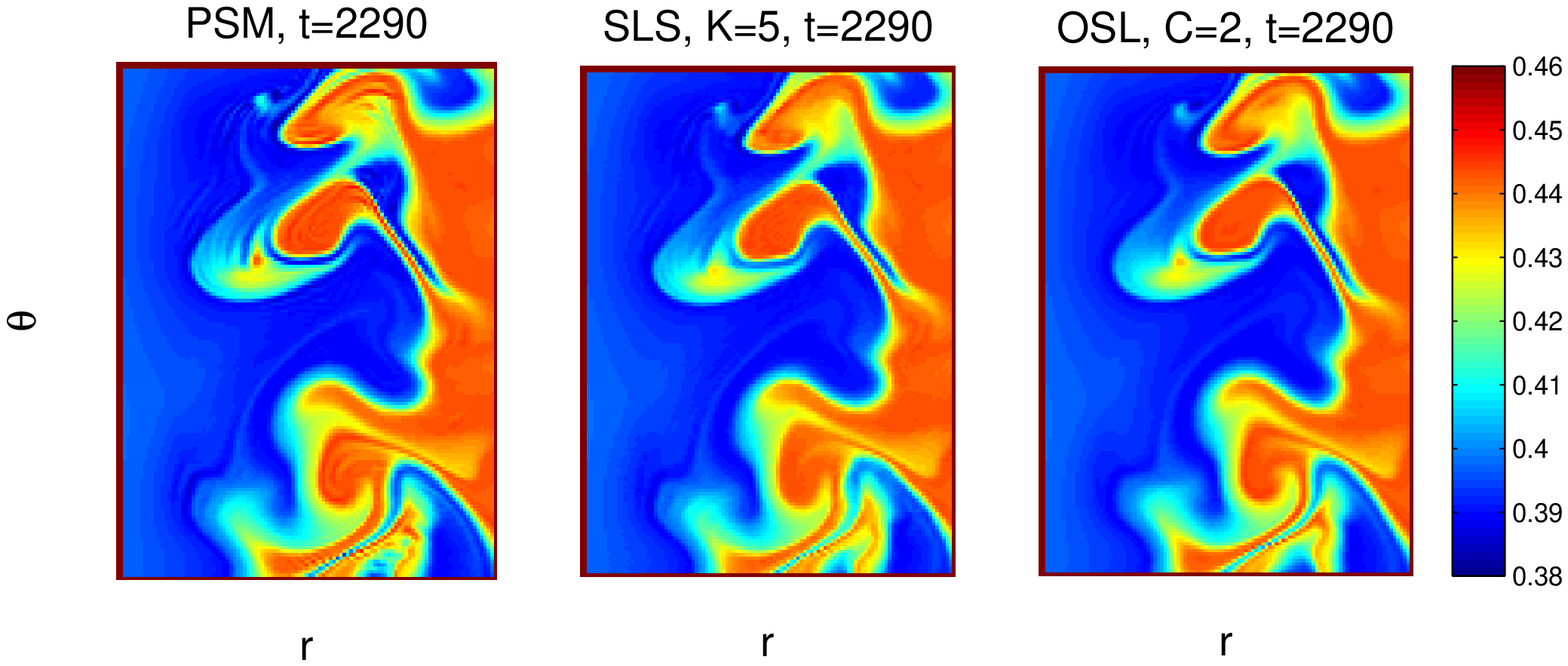}  \\
   \includegraphics*[scale=0.6,viewport=40 90 600 340]{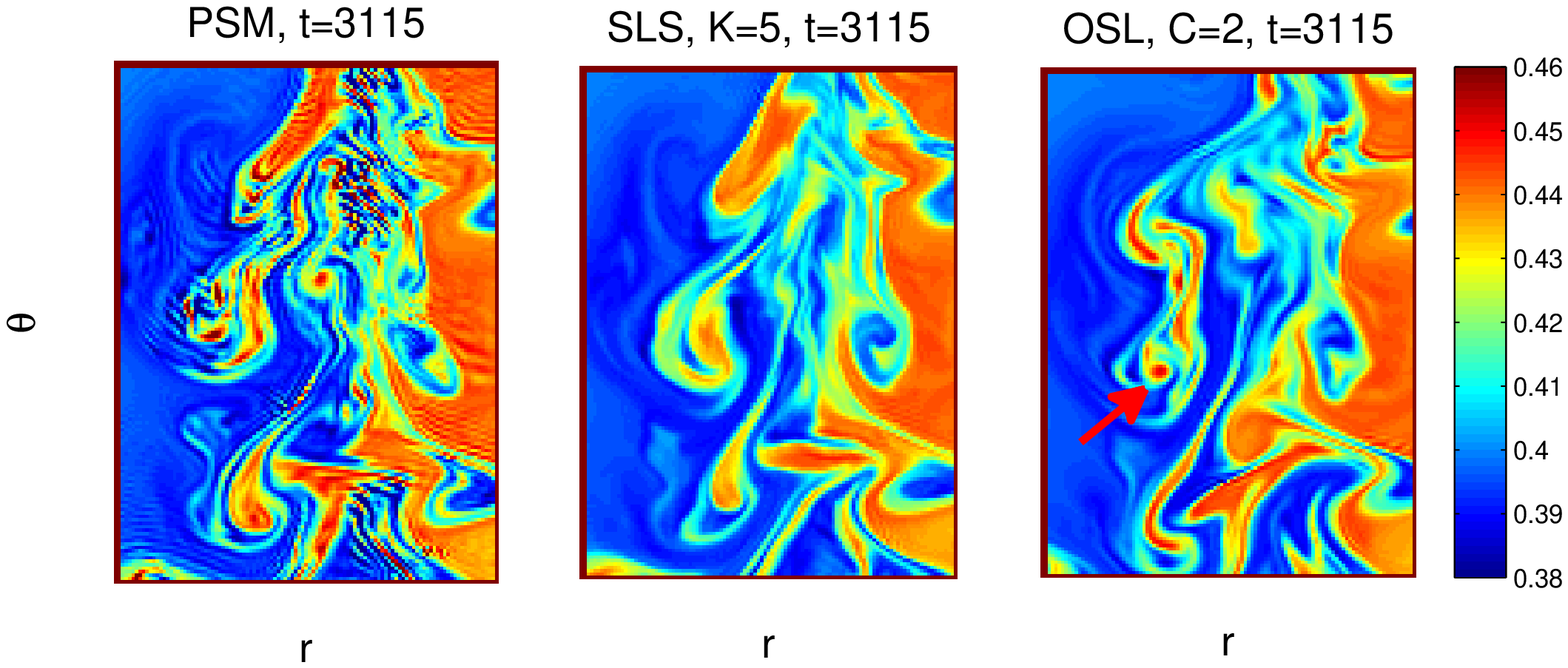} \\  
   \includegraphics*[scale=0.6,viewport=40 90 600 340]{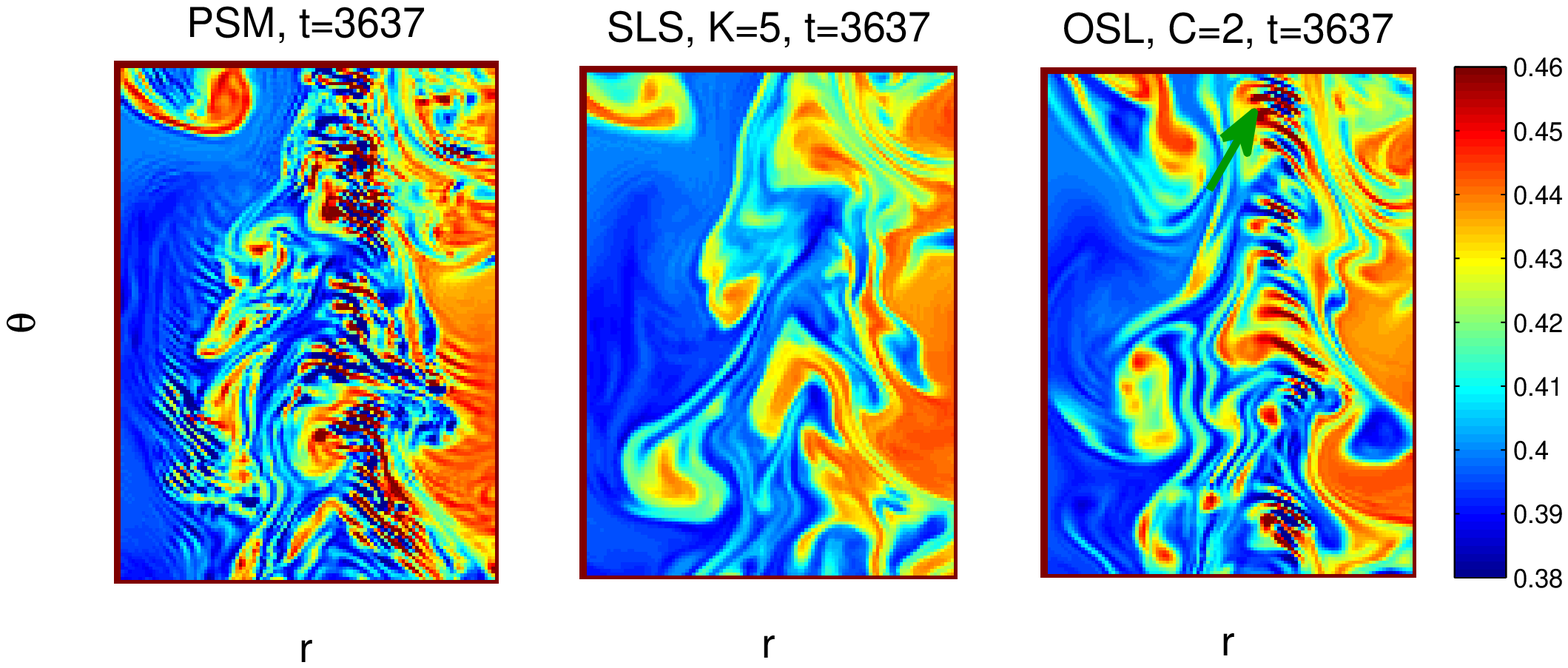} \\
 \end{tabular}
 \caption{Zoom on a the zone B of the high-resolution $(r,\theta)$ profiles (\fref{HR4D1})}
 \label{fig:zoomB}
 \end{figure}         

\begin{figure}
 \centering
 \begin{tabular}{cc}
  \includegraphics[scale=0.32]{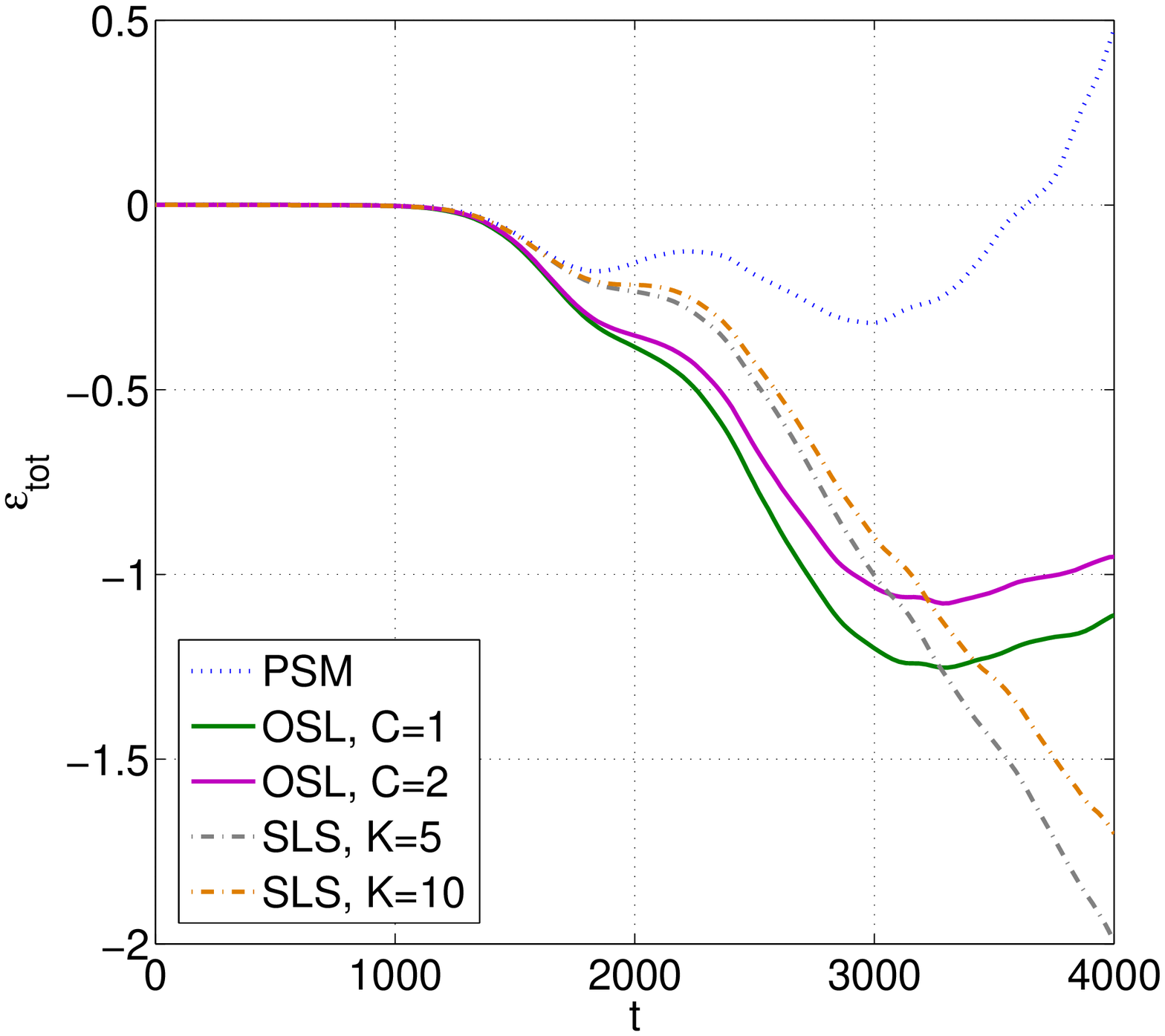} & \includegraphics[scale=0.32]{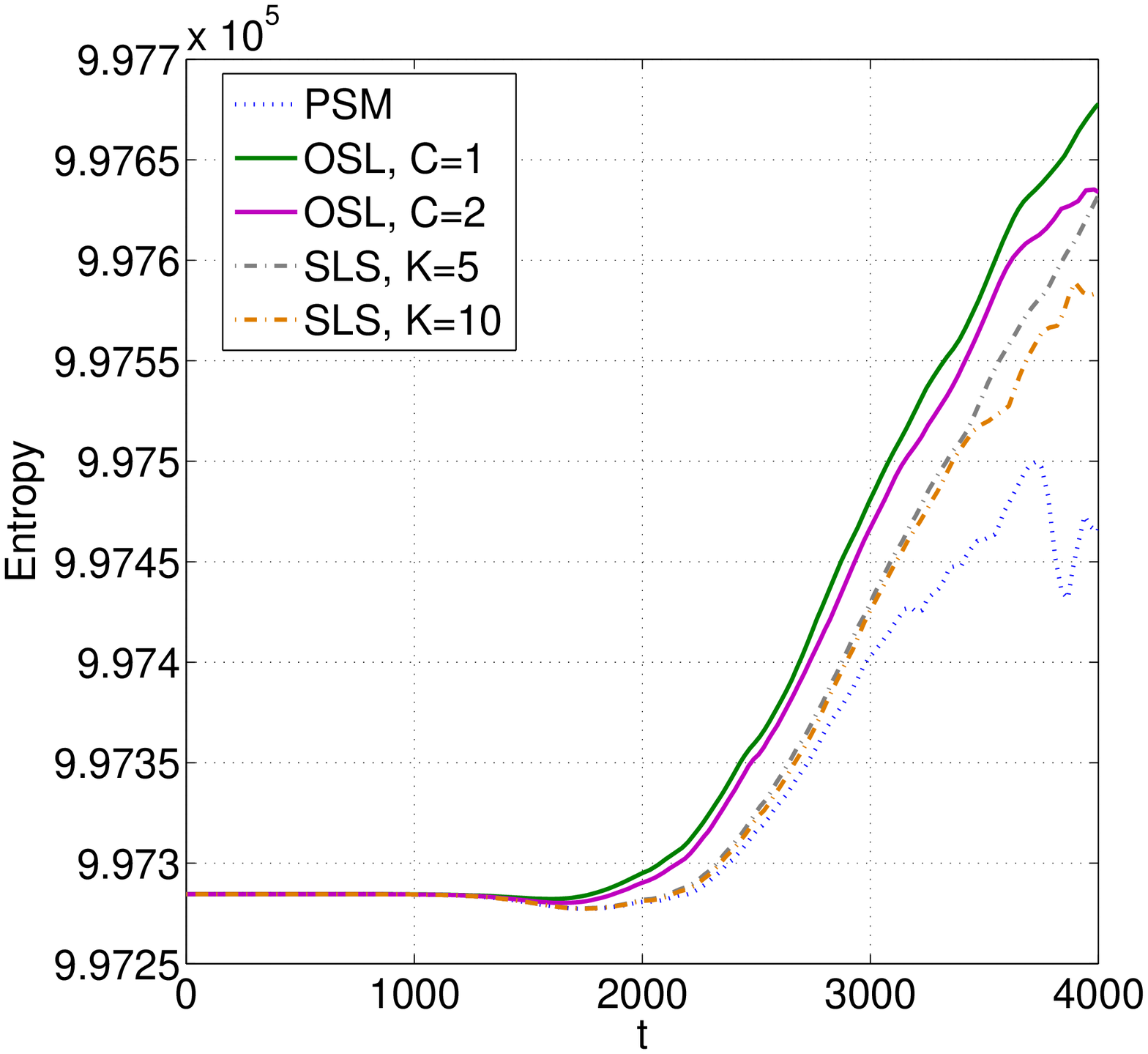}\\
  \includegraphics[scale=0.32]{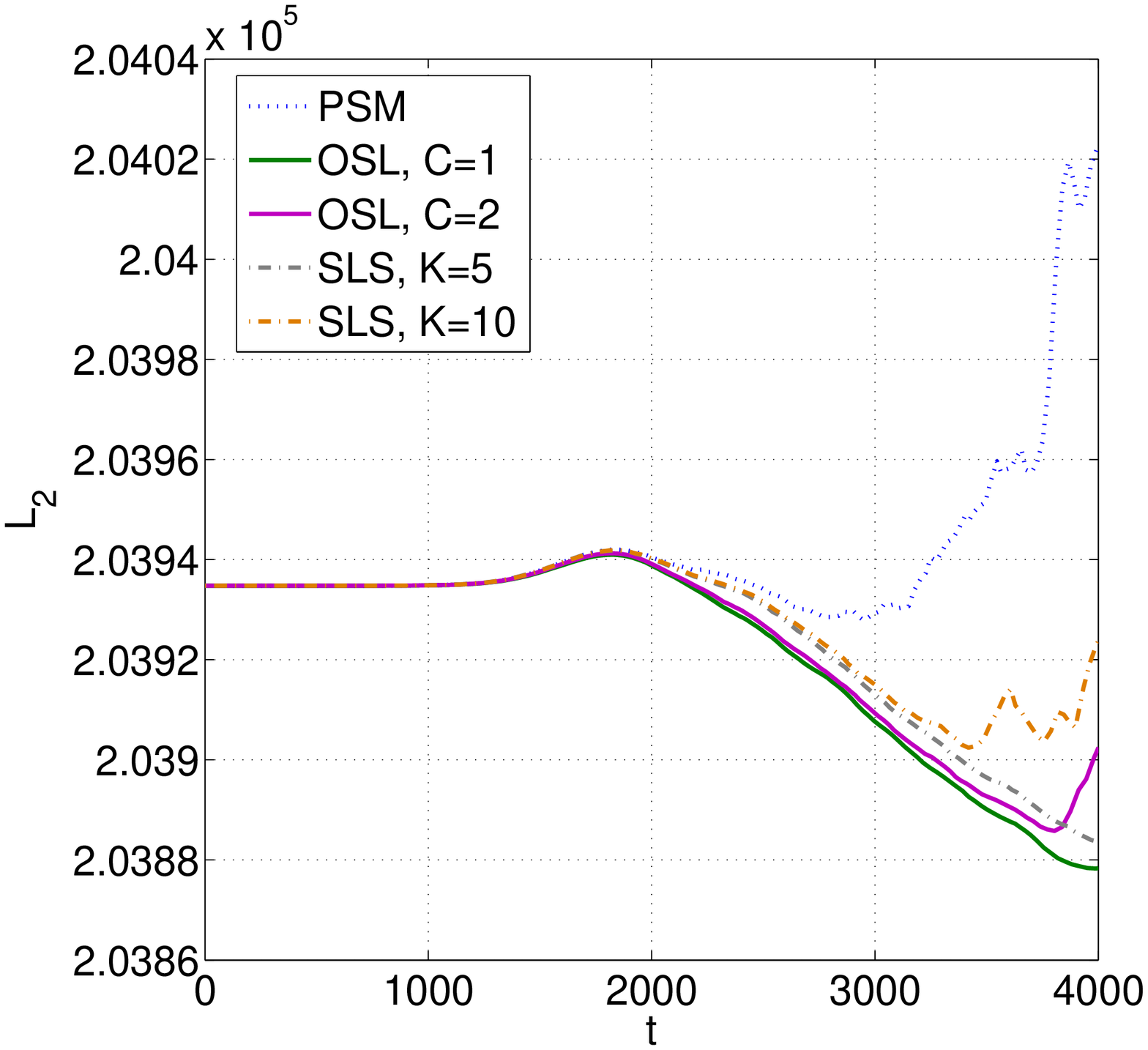}         & \includegraphics[scale=0.32]{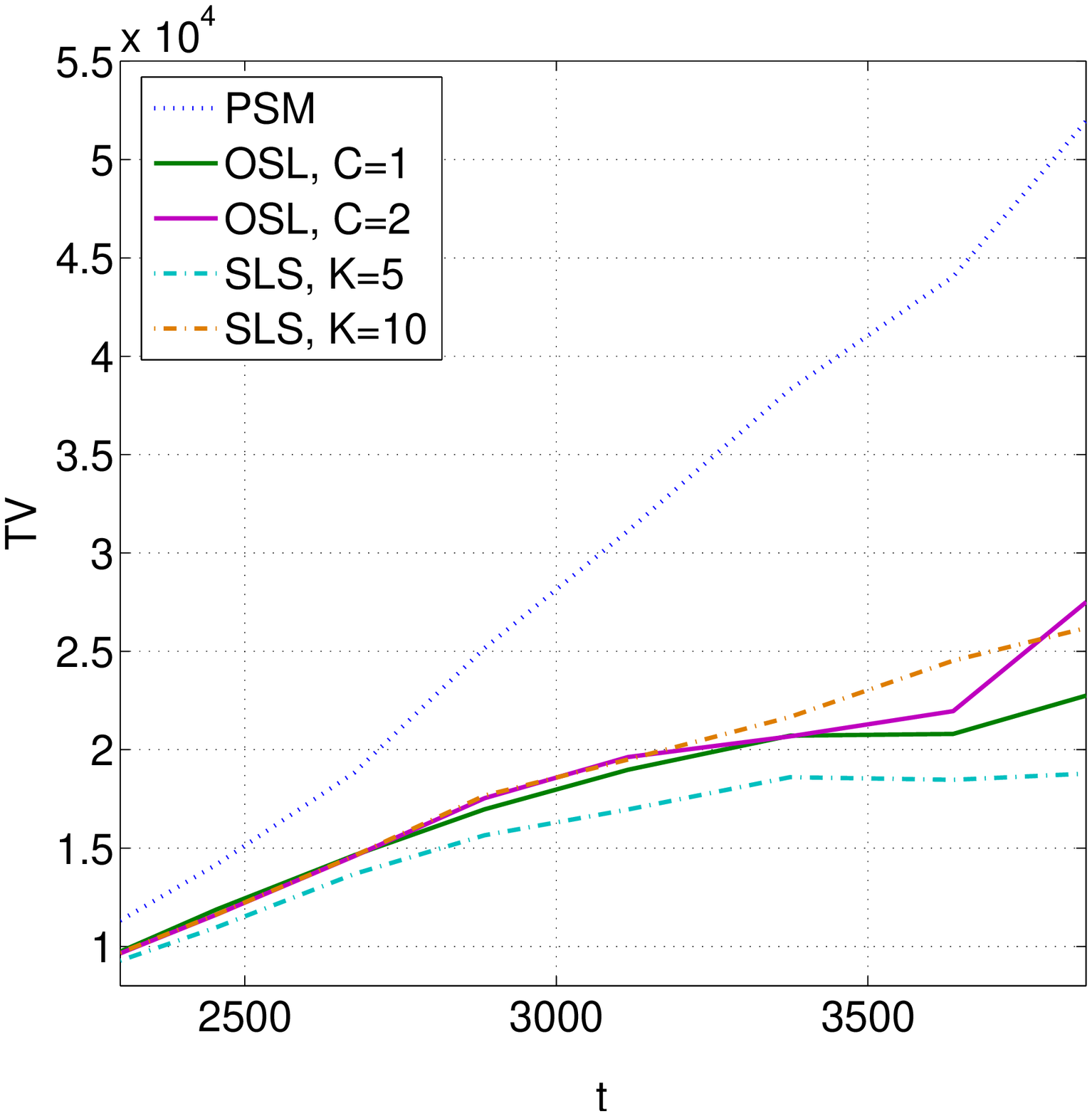}
 \end{tabular}
   \includegraphics[scale=0.4]{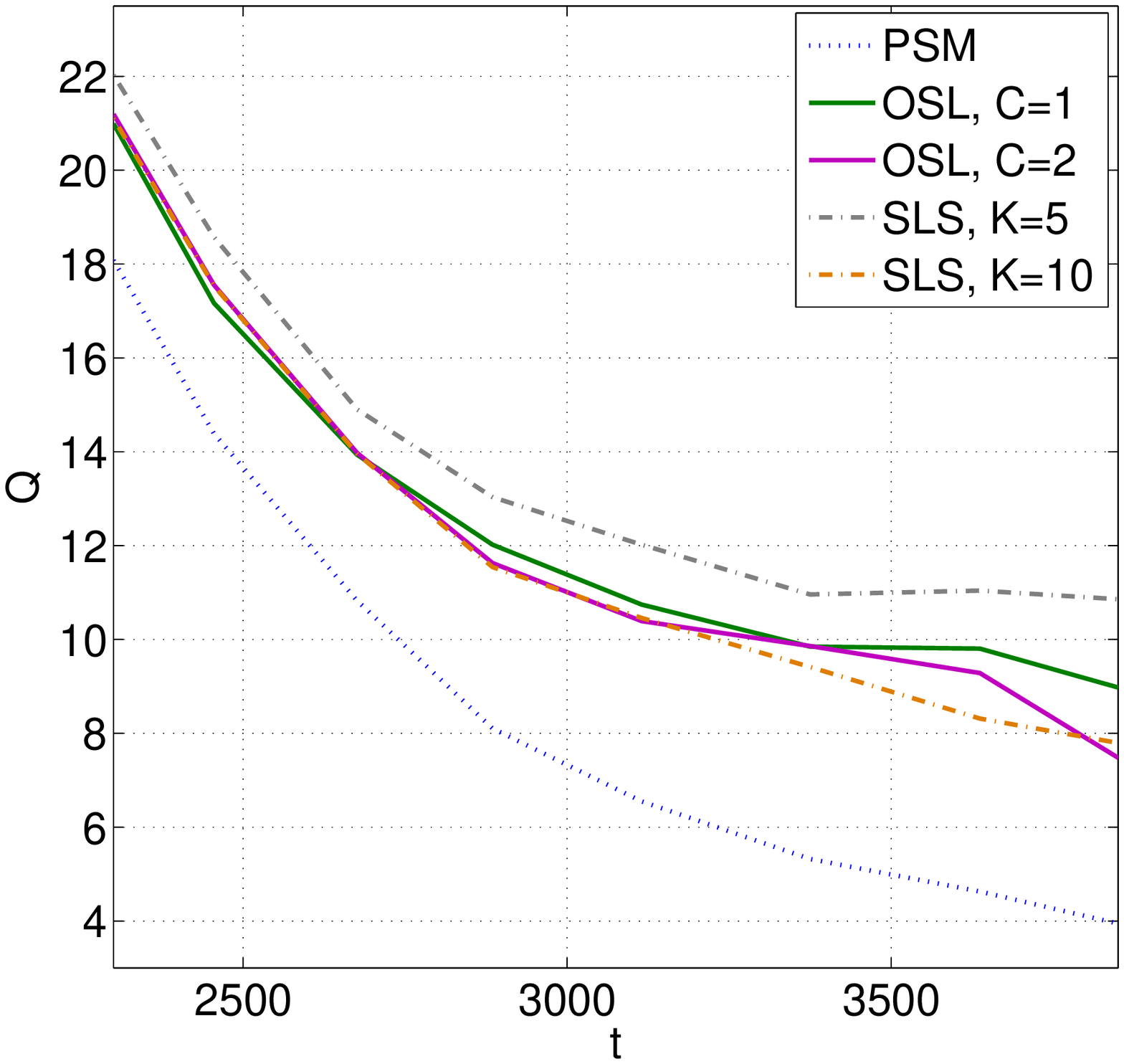}
 \caption{Total energy (left top), entropy (right top), norm $L^2$ (left bottom) and TV norm (right bottom) for the 4D simulations (high resolution), quality factor Q on the lower part.}
 \label{fig:4DHRall}
 \end{figure}

\newpage 
 \paragraph{Comparison of 1D integrated quantities}
We present in figure \ref{fig:4DHRall} the time evolution of the total energy, the entropy, the $L^2$ and TV norms and the quality factor \eqref{quality}. The results are for the PSM scheme with no limiter, the PSM-SLS scheme with K=5 and K=10 and the PSM-OSL scheme with C=1 and C=2.  \\
 \begin{itemize}
\item The PSM scheme conserves the total energy the better  until the simulation crashes down. Entropy and $L^2$ norm show that this scheme is the less diffusive but with the highest TV norm, which comes from the spurious oscillations seen on $(r,\theta)$ profiles and leads to a bad quality factor. 
\item On the contrary of what expected looking at  $(r,\theta)$ profiles, the entropy and $L^2$ norm show that OSL limiter leads to a little more dissipation than the SLS limiter. The total energy is a little better conserved with the SLS limiter than with the OSL limiter at the beginning of the simulation, but it constantly decreases while the total energy obtained with the OSL limiter stabilise at late times. 
 \item The quality factor is a little better with PSM-SLS than with PSM-OSL. However, both achieve the objective of stabilising the PSM scheme oscillations. 
\end{itemize}

\newpage
\section{Conclusion}
In this paper, we have investigated conservative schemes for a Vlasov 4D drift-kinetic model for plasma turbulence. These schemes should be accurate to bear long time simulations, so we investigate fourth order schemes. However, high order schemes usually experience difficulties when dealing with steep gradients in the flow, which happens when turbulence develops thin structures or vortices. In particular, spurious oscillations may appear which are not damped by Vlasov models.  \\
We thus have tested limiters to cut off these oscillations. It consists in adding diffusion to the high order scheme at locations where oscillations may appear, i.e. at steep gradients. Then two questions show up:  how to introduce diffusion in the scheme and how to detect the location where diffusion is needed. Of course, it should be done efficiently without loosing the accuracy of the scheme by introducing too much diffusion. We have tested the LAG and PSM schemes and the Umeda, ENT, OSL, SLS limiters, which use quite different strategies. We have tested all the schemes on a 1D step linear advection test-case and on a 4D drift-kinetic model test-case.\\
Another question is how  to quantify the action of limiters and to propose a quantitative way of comparison among them. In a classical way, we have checked the conservation of the Vlasov equation invariants as the $L^2$ norm, the entropy and the total energy for the drift-kinetic model. Even if these quantities permit to compare the behaviour of the schemes, they do not provide a quantitative way to tell one scheme is better than another. We propose an attempt  for a quality factor Q \eqref{quality}, based on the dissipation measured with the $L^2$ norm and spurious oscillations measured with the total variation $TV$ norm. It gives a clear answer the question which is the best of two schemes, but obviously it is still subjective because of the choice of quality factor itself. However, comparison of the schemes results with this quality factor is in good agreement with intuition when looking at results graphs in 1D, but in 4D it helps to evaluate the balance between diffusion and spurious oscillations. At the end of this work, we will not conclude that one scheme is better than another, we just have given these results and tools to help the reader to decide by himself.

\newpage

\section{Annexes}
\subsection{Annexe A: Detailed LAG reconstruction}
\label{sec:annexeA} 
We give hereafter detailed calculations for the LAG reconstruction.

First, we build the primitive function G at the characteristic feet:
\Eas
g^{n+1}_i=G(x^*_{i+1/2})-G(x^*_{i-1/2})
\\G(x^*_{i+1/2})=\int^{x^*_{i+1/2}}_{x_{-1/2}} g^{n}_{j(i)}(x)dx
\Ease
$G(x^*_{i+1/2})$ is given by the interpolation with Lagrange polynoms:
\[G(x)=\sum^{i+2}_{j=i-1} G_{j-1/2}L_j(x), x\in[x_{i-1/2},x_{i+1/2}]\]
where $L_j$ are the Lagrange polynoms defines as:
\[L_j(x)=\prod_{k=i-1,k\neq j}^{i+2}\frac{x-x_{k-1/2}}{x_{j-1/2}-x_{k-1/2}}\]
Hence,
\Eas
L_{i-1}(x)=\prod_{k=i-1,k\neq i-1}^{i+2}\frac{x-x_{k-1/2}}{x_{i-3/2}-x_{k-1/2}}
\\L_{i}(x)=\prod_{k=i-1,k\neq i}^{i+2}\frac{x-x_{k-1/2}}{x_{i-1/2}-x_{k-1/2}}
\\L_{i-1}(x)=\prod_{k=i-1,k\neq i+1}^{i+2}\frac{x-x_{k-1/2}}{x_{i+1/2}-x_{k-1/2}}
\\L_{i-1}(x)=\prod_{k=i-1,k\neq i+2}^{i+2}\frac{x-x_{k-1/2}}{x_{i+3/2}-x_{k-1/2}}
\Ease
We note $x=x_{i+1/2}-\alpha,\alpha>0$ then
\Eas
L_{i-1}(x)=\prod_{k=i-1,k\neq i-1}^{i+2}\frac{-\alpha+x_{i+1/2}-x_{k-1/2}}{x_{i-3/2}-x_{k-1/2}}
\\L_{i}(x)=\prod_{k=i-1,k\neq i}^{i+2}\frac{-\alpha+x_{i+1/2}-x_{k-1/2}}{x_{i-1/2}-x_{k-1/2}}
\\L_{i-1}(x)=\prod_{k=i-1,k\neq i+1}^{i+2}\frac{-\alpha+x_{i+1/2}-x_{k-1/2}}{x_{i+1/2}-x_{k-1/2}}
\\L_{i-1}(x)=\prod_{k=i-1,k\neq i+2}^{i+2}\frac{-\alpha+x_{i+1/2}-x_{k-1/2}}{x_{i+3/2}-x_{k-1/2}}
\Ease
Assuming $\Delta x=x_{i+1/2}-x_{i-1/2}, \forall i$,
\Eas
L_{i-1}(x)=\prod_{k=i-1,k\neq i-1}^{i+2}\frac{-\alpha+\Delta x(i-k+1)}{\Delta x(i-1-k)}
\\L_{i}(x)=\prod_{k=i-1,k\neq i}^{i+2}\frac{-\alpha+\Delta x(i-k+1)}{\Delta x(i-k)}
\\L_{i-1}(x)=\prod_{k=i-1,k\neq i+1}^{i+2}\frac{-\alpha+\Delta x(i-k+1)}{\Delta x(i+1-k)}
\\dL_{i-1}(x)=\prod_{k=i-1,k\neq i+2}^{i+2}\frac{-\alpha+\Delta x(i-k+1)}{\Delta x(i+2-k)}
\Ease
Writing $\beta=\alpha / \Delta x$, we deduce:
\Eas
L_{i-1}(x)=\prod_{k=i-1,k\neq i-1}^{i+2}\frac{-\beta+(i-k+1)}{(i-1-k)}
\\L_{i}(x)=\prod_{k=i-1,k\neq i}^{i+2}\frac{-\beta+ (i-k+1)}{(i-k)}
\\L_{i+1}(x)=\prod_{k=i-1,k\neq i+1}^{i+2}\frac{-\beta+ (i-k+1)}{(i+1-k)}
\\L_{i+2}(x)=\prod_{k=i-1,k\neq i+2}^{i+2}\frac{-\beta+(i-k+1)}{(i+2-k)}
\Ease
and then,
\Eas
L_{i-1}(x)=\frac{-\beta+(i-(i)+1)}{(i-1-(i))}\frac{-\beta+(i-(i+1)+1)}{i-1-(i+1)}\frac{-\beta+(i-(i+2)+1)}{(i-1-(i+2))}
\\L_{i}(x)=\frac{-\beta+ (i-(i-1)+1)}{(i-(i-1))}\frac{-\beta+ (i-(i+1)+1)}{(i-(i+1))}\frac{-\beta+ (i-(i+2)+1)}{(i-(i+2))}
\\L_{i+1}(x)=\frac{-\beta+ (i-(i-1)+1)}{(i+1-(i-1))}\frac{-\beta+ (i-(i+2)+1)}{(i+1-(i+2))}\frac{-\beta+ (i-(i)+1)}{(i+1-(i))}
\\L_{i+2}(x)=\frac{-\beta+(i-(i-1)+1)}{(i+2-(i-1))}\frac{-\beta+(i-(i)+1)}{(i+2-(i))}\frac{-\beta+(i-(i+1)+1)}{(i+2-(i+1))}
\Ease
\Eas
L_{i-1}(x)=\frac{-\beta+1}{-1}\frac{-\beta}{-2}\frac{-\beta-1}{-3}
\\L_{i}(x)=\frac{-\beta+2}{1}\frac{-\beta}{-1}\frac{-\beta-1}{-2}
\\L_{i+1}(x)=\frac{-\beta+2}{2}\frac{-\beta-1)}{-1}\frac{-\beta+ 1}{1}
\\L_{i+2}(x)=\frac{-\beta+2}{3}\frac{-\beta+1}{2}\frac{-\beta}{1}
\Ease
\Eas
L_{i-1}(x)=1/6(\beta-1)(\beta)(\beta+1)
\\L_{i}(x)=1/2(\beta- 2)(\beta)(\beta+1)
\\L_{i+1}(x)=1/2(\beta-2)(\beta+1)(\beta-1)
\\L_{i+2}(x)=-1/6(\beta-2)(\beta-1)(\beta)
\Ease

\Eas
L_{i-1}(x)=1/6(\beta^3-\beta)
\\L_{i}(x)=-1/2(\beta^3-\beta^2-2\beta)
\\L_{i+1}(x)=1/2(\beta^3-2\beta^2-\beta+2)
\\L_{i+2}(x)=-1/6(\beta^3-3\beta^2+2\beta)
\Ease
Hence
\Eas
G(x)&=&1/6G_{i-3/2}(\beta^3-\beta)
\\&-&1/2G_{i-1/2}(\beta^3-\beta^2-2\beta)
\\&+&1/2G_{i+1/2}(\beta^3-2\beta^2-\beta+2)
\\&-&1/6G_{i+3/2}(\beta^3-3\beta^2+2\beta)
\Ease

Finally, the interpolation can be written:
\Eas
G(x)&=&\beta^3(1/6 G_{i-3/2}-1/2 G_{i-1/2}+1/2G_{i+1/2}-1/6 G_{i+3/2})
\\ &+&\beta^2 (1/2 G_{i-1/2}-G_{i+1/2}+1/2 G_{i+3/2})
\\&+&\beta(-1/6G_{i-3/2}+G_{i-1/2}-1/2 G_{i+1/2}-1/3G_{i+3/2})
\\&+&G_{i+1/2}
\Ease

We note $x=x_{i+1/2}-\Delta x \beta=x_{i-1/2}-\Delta x(\beta-1)$ \newline
and by defining $\theta=1-\beta>0$, $x=x_{i-1/2}+\Delta x\theta$, then
\Eas
G(x)&=&1/6G_{i-3/2}((1-\theta)^3-(1-\theta))
\\&-&1/2G_{i-1/2}((1-\theta)^3-(1-\theta)^2-2(1-\theta))
\\&+&1/2G_{i+1/2}((1-\theta)^3-2(1-\theta)^2-(1-\theta)+2)
\\&-&1/6G_{i+3/2}((1-\theta)^3-3(1-\theta)^2+2(1-\theta)))
\Ease

\Eas
G(x)&=&1/6G_{i-3/2}(-\theta^3+3\theta^2-2\theta)
\\&-&1/2G_{i-1/2}(-\theta^3+2\theta^2+\theta-2)
\\&+&1/2G_{i+1/2}(-\theta^3+\theta^2+2\theta)
\\&-&1/6G_{i+3/2}(-\theta^3+\theta)
\Ease
 
\Eas
G(x)&=&\theta^3(-1/6 G_{i-3/2}+1/2 G_{i-1/2}-1/2G_{i+1/2}+1/6 G_{i+3/2})
\\ &+&\theta^2 (1/2 G_{i-3/2}-G_{i-1/2}+1/2 G_{i+1/2} )
\\&+&\theta(-1/3 G_{i-3/2}-1/2G_{i-1/2}+ G_{i+1/2}-1/6G_{i+3/2})
\\&+&G_{i-1/2}
\Ease

Using the previous formula with $x=x_{i+1/2}-\Delta x \beta$, $\beta<0$, $i\rightarrow i+1$, $\theta \rightarrow -\beta$, we obtain for negative displacement
\Eas
G(x)&=&\beta^3(1/6 G_{i-1/2}-1/2 G_{i+1/2}+1/2G_{i+3/2}-1/6 G_{i+5/2})
\\ &+&\beta^2 (1/2 G_{i-1/2}-G_{i+1/2}+1/2 G_{i+3/2} )
\\&+&\beta(1/3 G_{i-1/2}+1/2G_{i+1/2}-G_{i+3/2}+1/6G_{i+5/2})
\\&+&G_{i+1/2}
\Ease

We obtain the following results with $x=x_{i+1/2}-\Delta x \beta$:
\newline if $\beta>0$, 
\Eas
G(x)&=&\beta^3(1/6 G_{i-3/2}-1/2 G_{i-1/2}+1/2G_{i+1/2}-1/6 G_{i+3/2})
\\ &+&\beta^2 (1/2 G_{i-1/2}-G_{i+1/2}+1/2 G_{i+3/2})
\\&+&\beta(-1/6G_{i-3/2}+G_{i-1/2}-1/2 G_{i+1/2}-1/3G_{i+3/2})
\\&+&G_{i+1/2}
\Ease
\newline if $\beta<0$,
 \Eas
G(x)&=&\beta^3(1/6 G_{i-1/2}-1/2 G_{i+1/2}+1/2G_{i+3/2}-1/6 G_{i+5/2})
\\ &+&\beta^2 (1/2 G_{i-1/2}-G_{i+1/2}\\&+&1/2 G_{i+3/2} )
\\&+&\beta(1/3 G_{i-1/2}+1/2G_{i+1/2}-G_{i+3/2}+1/6G_{i+5/2})
\\&+&G_{i+1/2}
\Ease
We define $g^{n}_i\Delta x=G_{i+1/2}-G_{i-1/2}$.
The previous expressions become:
\newline if $\beta>0$, 
\Eas
G(x)&=&\beta^3(1/6 G_{i-3/2}-1/6 G_{i-1/2}-2/6G_{i-1/2}\\&+&2/6G_{i+1/2}+1/6G_{i+1/2}-1/6 G_{i+3/2})
\\ &+&\beta^2 (1/2 G_{i-1/2}-1/2G_{i+1/2}-1/2G_{i+1/2}+1/2 G_{i+3/2})
\\&+&\beta(-1/6G_{i-3/2}+1/6G_{i-1/2}+5/6G_{i-1/2}\\&-&5/6 G_{i+1/2}+1/3 G_{i+1/2}-1/3G_{i+3/2})
\\&+&G_{i+1/2}
\\&=&\beta^3(-1/6g^{n}_{i-1}+1/3g^{n}_{i}-1/6g^{n}_{i+1} )
\\ &+&\beta^2 (-1/2g^{n}_{i}+1/2g^{n}_{i+1})
\\&+&\beta(1/6g^{n}_{i-1}-5/6g^{n}_{i}-1/3g^{n}_{i+1})
\\&+&G_{i+1/2}
\Ease
\newline if $\beta<0$, \Eas
G(x)&=&\beta^3(1/6 G_{i-1/2}-1/6 G_{i+1/2}-2/6G_{i+1/2}\\&+&2/6G_{i+3/2}+1/6G_{i+3/2}-1/6 G_{i+5/2})
\\ &+&\beta^2 (1/2 G_{i-1/2}-1/2G_{i+1/2}-1/2G_{i+1/2}+1/2 G_{i+3/2} )
\\&+&\beta(1/3 G_{i-1/2}-1/3G_{i+1/2}+5/6G_{i+1/2}\\&-&5/6G_{i+3/2}-1/6G_{i+3/2}+1/6G_{i+5/2})
\\&+&G_{i+1/2}
\\&=&\beta^3(-1/6g^{n}_{i} +2/6g^{n}_{i+1}-1/6g^{n}_{i+2})
\\ &+&\beta^2 (-1/2g^{n}_{i}+1/2g^{n}_{i+1} )
\\&+&\beta(-1/3g^{n}_{i}-5/6g^{n}_{i+1}+1/6g^{n}_{i+2})
\\&+&G_{i+1/2}
\Ease

Let now compare theses fluxe expressions with the regular expression we found for the Hermite formulation noticed ${\phi_{k+1/2,j_{k+1/2}=k}(\beta)}$.

\subsection{Positive displacement: $\beta>0,\delta=1,j_{k+1/2}=k$}
 \Eas
\frac{\phi_{k+1/2,j_{k+1/2}=k}(\beta)}{\Delta x} &=& g^{+}_{{k}-1/2} \left(\beta^2(\beta-1) \right)\\ &+& g^{-}_{{k}+1/2}\left( \beta(1-\beta)^2\right)\\ &+& \bar{g}_{{k}}\left(\beta^2(3-2\beta)\right)
\\&=&\beta^3(g^{+}_{{k}-1/2}+g^{-}_{{k}+1/2}-2\bar{g}_{{k}})
\\&+& \beta^2(-g^{+}_{{k}-1/2}-2g^{-}_{{k}+1/2}+3\bar{g}_{{k}})
\\&+& \beta(g^{-}_{k+1/2})
\Ease
We have introduced the value of the distribution at the cell's faces that we note $g^{+}_{k-1/2}=g(x\geqq x_{k-1/2})$ and $g^{-}_{k+1/2}=g(x\leqq x_{k+1/2})$.
We set these coefficients:
\[\begin{cases}
g^{+}_{{k}-1/2}=1/3g^{n}_{k-1}+5/6 g^{n}_{k}-1/6g^{n}_{k+1}  \\ 
g^{-}_{{k}+1/2}=-1/6g^{n}_{k-1}+5/6 g^{n}_{k}+1/3g^{n}_{k+1}\\
\end{cases}
\]
that leads to:
\Eas 
\frac{\phi_{k+1/2,j_{k+1/2}=k}(\beta)}{\Delta x} &=& g_{{k}-1/2} \left(\beta^2(\beta-1) \right)\\ &+& g_{{k}+1/2}\left( \beta(1-\beta)^2\right)\\ &+& \bar{g}_{{k}}\left(\beta^2(3-2\beta)\right)
\\&=&\beta^3(1/3g^{n}_{k-1}+5/6 g^{n}_{k}-1/6g^{n}_{k+1}-1/6g^{n}_{k-1}\\&+&5/6 g^{n}_{k}+1/3g^{n}_{k+1}-2{g}^{n}_{{k}})
\\&+& \beta^2(-1/3g^{n}_{k-1}-5/6 g^{n}_{k}+1/6g^{n}_{k+1} \\&+&1/3g^{n}_{k-1}-10/6 g^{n}_{k}-2/3g^{n}_{k+1}+3g^{n}_{k+1})
\\&+& \beta(-1/6g^{n}_{k-1}+5/6 g^{n}_{k}+1/3g^{n}_{k+1})
\Ease
\Eas 
\frac{\phi_{k+1/2,j_{k+1/2}=k}(\beta)}{\Delta x} &=& g_{{k}-1/2} \left(\beta^2(\beta-1) \right)\\ &+& g_{{k}+1/2}\left( \beta(1-\beta)^2\right)\\ &+& \bar{g}_{{k}}\left(\beta^2(3-2\beta)\right)
\\&=&\beta^3(1/6g^{n}_{k-1}-1/3{g}^{n}_{{k}}+1/6g^{n}_{k+1})
\\&+& \beta^2(1/2 g^{n}_{k}-1/2g^{n}_{k+1})
\\&+& \beta(-1/6g^{n}_{k-1}+5/6 g^{n}_{k}+1/3g^{n}_{k+1})
\Ease
Remembering that  
\Eas
G(\beta)&=&\beta^3(-1/6g^{n}_{i-1}+1/3g^{n}_{i}-1/6g^{n}_{i+1} )
\\ &+&\beta^2 (-1/2g^{n}_{i}+1/2g^{n}_{i+1})
\\&+&\beta(1/6g^{n}_{i-1}-5/6g^{n}_{i}-1/3g^{n}_{i+1})
\\&+&G_{i+1/2}
\Ease
we deduce:
\[
G(x)=-\phi_{k+1/2,j_{k+1/2}=k}(\beta)+G_{i+1/2}
\]

\subsection{Negative displacement $\delta=0$}
\Eas
\frac{\phi_{k+1/2,{j_{k+1/2}}=k+1}(\beta)}{\Delta x} &=& g_{k+1/2} \left(\beta(\beta+1)^2 \right)
\\ &+&g_{k+3/2}\left( \beta^2(1+\beta)\right)
\\ &+& {g}^{n}_{k+1}\left(\beta^2(-3-2\beta)\right)
\\ &=& \beta^3(g^{-}_{k+3/2}+g^{+}_{k+1/2}-2{g}^{n}_{k+1})
\\ &+& \beta^2(2g^{+}_{k+1/2}+g^{-}_{k+3/2}-3g^{n}_{k+1})
\\ &+& \beta (g^{+}_{k+1/2})
\Ease
 \Eas
G(x)&=&\beta^3(-1/6g^{n}_{i} +2/6g^{n}_{i+1}-1/6g^{n}_{i+2})
\\ &+&\beta^2 (-1/2g^{n}_{i}+1/2g^{n}_{i+1} )
\\&+&\beta(-1/3g^{n}_{i}-5/6g^{n}_{i+1}+1/6g^{n}_{i+2})
\\&+&G_{i+1/2}
\Ease

We set these coefficients:
\[
\begin{cases}
g^{+}_{{k}+1/2}=1/3g^{n}_{k}+5/6 g^{n}_{k+1}-1/6g^{n}_{k+2}  \\ 
g^{-}_{{k}+3/2}=-1/6g^{n}_{k}+5/6 g^{n}_{k+1}+1/3g^{n}_{k+2}\\
\end{cases}
\]
\Eas
\frac{ \phi_{k+1/2,{j_{k+1/2}}=k+1}(\beta)}{\Delta x} &=& \beta^3(1/3g^{n}_{k}+5/6 g^{n}_{k+1}-1/6g^{n}_{k+2}-1/6g^{n}_{k}\\&+&5/6 g^{n}_{k+1}+1/3g^{n}_{k+2}-2{g}^{n}_{k+1})
\\ &+& \beta^2(2/3g^{n}_{k}+10/6 g^{n}_{k+1}-2/6g^{n}_{k+2}-1/6g^{n}_{k}\\&+&5/6 g^{n}_{k+1}+1/3g^{n}_{k+2}-3g^{n}_{k+1})
\\ &+& \beta (1/3g^{n}_{k}+5/6 g^{n}_{k+1}-1/6g^{n}_{k+2})
\Ease
\Eas
\frac{ \phi_{k+1/2,{j_{k+1/2}}=k+1}(\beta)}{\Delta x} &=& \beta^3(1/6g^{n}_{k}-1/3 g^{n}_{k+1}+1/6g^{n}_{k+2})
\\ &+& \beta^2(1/2g^{n}_{k}-1/2g^{n}_{k+1})
\\ &+& \beta (1/3g^{n}_{k}+5/6 g^{n}_{k+1}-1/6g^{n}_{k+2})
\Ease
Remembering that:
 \Eas
G(x)&=&\beta^3(-1/6g^{n}_{i} +1/3g^{n}_{i+1}-1/6g^{n}_{i+2})
\\ &+&\beta^2 (-1/2g^{n}_{i}+1/2g^{n}_{i+1} )
\\&+&\beta(-1/3g^{n}_{i}-5/6g^{n}_{i+1}+1/6g^{n}_{i+2})
\\&+&G_{i+1/2}
\Ease
we deduce
\[
G(x)=-\phi_{k+1/2,j_{k+1/2}=k+1}(\beta)/\Delta x+G_{i+1/2}
\]
To conclude,
\Eas
g^{new}_i&=&-\phi_{k+1/2,j_{k+1/2}}(\beta)/\Delta x+G_{i+1/2}+G_{k-1/2,j_{k-1/2}}(\beta)/\Delta x-G_{i-1/2}
\\g^{new}_i&=&g^{n}_i-(\phi_{k+1/2,j_{k+1/2}}(\beta)+\phi_{k-1/2,j_{k-1/2}}(\beta))/\Delta x
\Ease

\subsection{Annexe B: Algorithms}
\label{sec:annexeB}
We present some algorithms which are representative of the Gysela algorithms. We change the notation of the flux $\phi_{i+1/2}$ to $H_{i+1/2}$ for practical reasons. The 1D advection equation are solved with \algoref{main}.
\begin{algorithm}[H]
\caption{Advection}
\begin{algorithmic}[1]

\REQUIRE{$\left(x_k \right)_{k={-1/2}\cdots N+1/2},\left(\bar{g}^{n}_k \right)_{k=0\cdots N}$}
\STATE a=1,b=-1,c=0,p=1,q=4
\STATE L,U=MatrixCoefficient(a,b,c,p,q)
\STATE $\left(g^{n}_{k-1/2}\right)_{k=0\cdots N+1}$=CoefficientHermite($\left(\bar{g}^{n}_k \right)_{k={0}\cdots N}$,L,N,0,N+1,a,b,c)
\STATE $\beta$=FootCharacteristic(...)
\IF{$\beta=0$}
\STATE $H_{0-1/2}=0$
\ELSE
\IF{$\beta>0$}
	\STATE $H_{0-1/2}=\beta \bar{g}^{n}_0$
\ELSE  
	\STATE $H_{0-1/2}$=FluxComputation($\bar{g}^{n}_{0},g^{n}_{0-1/2},g^{n}_{0+1/2},\beta$)
\ENDIF 
\ENDIF
\FOR{i=1:N-1}
	\STATE $\beta$=FootCharacteristic()
	\IF{$\beta=0$}
	\STATE $H_{i+1/2}=0$
	\ELSE
	\IF{$\beta>0$}
		\STATE j=i
	\ELSE 
		\STATE j=i+1
	\ENDIF

\STATE $H_{i+1/2}$=FluxComputation($\bar{g}^{n}_{j},g^{n}_{j-1/2},g^{n}_{j+1/2},\beta$)
\ENDIF
\ENDFOR
\STATE $\beta$=FootCharacteristic(...)
\IF{$\beta=0$}
\STATE
$H_{N+1/2}=0$
\ELSE
\IF{$\beta<0$}
	\STATE $H_{N+1/2}=\beta \bar{g}^{n}_N$
\ELSE  
	\STATE $H_{N+1/2}$=FluxComputation($\bar{g}^{n}_{N},g^{n}_{N-1/2},g^{n}_{N+1/2},\beta$)
\ENDIF
\ENDIF
\FOR{i=1:N-1}
\STATE $g^{n+1}_i=g^{n}_i-(H_{i+1/2}-H_{(i-1)+1/2})$
\ENDFOR
\RETURN $g^{n+1}_{k=0\cdots N}$

\end{algorithmic}
\label{main}\end{algorithm}

\subsubsection{FluxComputation and FootCharacteristic}
\begin{algorithm}[H]
\caption{ FluxComputation}
\begin{algorithmic}[1]
\REQUIRE $\bar{g}^{n}_{j_k},g^{n}_{j_k-1/2},g^{n}_{j_k+1/2},\beta$


\IF {$\beta>0$}
   \STATE $\delta= 0$
   \ELSE
    \STATE $\delta= 1$
\ENDIF
\STATE $H_{k+1/2}=g_{j_k-1/2}\left( \beta(1-\delta)+\beta^2(2-3\delta)+\beta^3\right)$
\STATE$H_{k+1/2}=H_{k+1/2}+g_{j_k+1/2}\left( \beta \delta +\beta^2(1-3\delta)+\beta^3\right)$
\STATE$H_{k+1/2}=H_{k+1/2}+\bar{g}_{j_k}\left(\beta^2(-3+6\delta)+\beta^3 (-2)\right)$
\RETURN $H_{k+1/2}$
\end{algorithmic}
\end{algorithm}

\begin{algorithm}[H]
\caption{FootCharacteristic}
\begin{algorithmic}[1]
\RETURN $x^*_{i+1/2}$
\end{algorithmic}
\end{algorithm}

\subsubsection{Hermite coefficient computation for natural conditions}
For a domain with natural boundary conditions, we compute the value of the distribution at the cell's faces $g_{k-1/2}$ and $g_{k+1/2}$. We want to solve: $AX=B$  and we decompose $A$ as $A=LU$ with \algoref{matrixnat} then we solve $X=U^{-1}L^{-1}B$ with \algoref{coeffnat}.

\begin{algorithm}[H]
\caption{CoefficientHermiteNatural}
\begin{algorithmic}[1]

\REQUIRE{$\left(\bar{g}^{n}_k \right)_{k={0}\cdots N},\left( L_k\right)_{k=0 \cdots N+1},\left( U_k\right)_{k=0 \cdots N+1}$,$0$,$N+1$,a,b,c}
\STATE $X_{0}=5\bar{g}^{n}_{0}$
\FOR{$k=1:N$}
\STATE $X_{k}=3 \left( \bar{g}^{n}_{k-1}+\bar{g}^{n}_{k} \right)-L_{k-1}X_{k-1}$
\ENDFOR
\STATE $X_{N+1}=5\bar{g}^{n}_{N}-L_{N+1}X_{N-1}-L_{N}X_{N}$
\STATE $X_{N+1}=\frac{X_{N+1}}{U_{N+1}}$
\FOR{$k=N:2,-1$}
\STATE $X_{k}=\frac{X_{k}-X_{k+1}}{U_k}$
\ENDFOR
\STATE $X_{1}=\frac{ X_{1}-\left(1-c/a\right)X_{2}}{U_k}$
\STATE $X_{0}=\frac{ X_{0}-bX_{1}-cX_{2}}{U_0}$
\RETURN $\left(g_{k-1/2}\right)_{k={0}\cdots N+1}$
\end{algorithmic}
\label{coeffnat}
\end{algorithm}

We compute the matrix coefficients of the LU decomposition.

\[ \begin{bmatrix} A \end{bmatrix}=\begin{bmatrix} L \end{bmatrix}\begin{bmatrix} U \end{bmatrix}
\]
\[
 \begin{bmatrix}
a     & b      & c      & 0      & \cdots & \cdots & 0      \\
p     & q      & p      & 0      & \cdots & \cdots & \vdots      \\
0     & p      & q      & p      & 0 & \cdots & \vdots \\
\vdots & \ddots & \ddots & \ddots & \ddots & \ddots & \vdots \\ 
\vdots & \cdots & 0 & p      & q      & p      & 0      \\
\vdots & \cdots & \cdots & 0      & p      & q      & p      \\
0     & \cdots & \cdots &0       & d      & e     & f
 \end{bmatrix}=\]
 \[
 \begin{bmatrix}

1     & 0      & \cdots      & \cdots      & \cdots & \cdots & 0      \\
L_1     & 1      & 0      & \cdots      & \cdots & \cdots & \vdots      \\
0     & L_2      & 1     & 0      & \cdots & \cdots & \vdots \\
\vdots & \ddots & \ddots & \ddots & \ddots & \ddots & \vdots \\ 
\vdots & \cdots & 0 & L_{N-2}      & 1      & 0      & \vdots      \\
\vdots & \cdots & \cdots & 0      & L_{N-1}      & 1      & 0      \\
0     & \cdots & \cdots &0       & L_{N+1}      & L_N     & 1
 \end{bmatrix}
 \begin{bmatrix}
U_0     & b      & c      & 0      & \cdots & \cdots & 0      \\
0     & U_1      & 1-\frac{c}{a} & 0      & \cdots & \cdots & \vdots      \\
\vdots     & 0      & U_2      & 1      & 0 & \cdots & \vdots \\
\vdots & \ddots & \ddots & \ddots & \ddots & \ddots & \vdots \\ 
\vdots & \cdots & \cdots & 0      &  U_{N-1}      & 1     & 0      \\
\vdots & \cdots & \cdots & \cdots     & 0      &  U_{N}      & 1      \\
0     & \cdots & \cdots &  \cdots       &\cdots     & 0     & U_{N+1}
 \end{bmatrix}
 \]

\begin{algorithm}[H]
\caption{MatrixCoefficientNatural}
\begin{algorithmic}[1]
\REQUIRE a,b,c,d,e,f,p,q,0,N+1
\STATE $U_{0}=a$
\STATE $L_{0}=\frac{p}{U_{0}}$
\STATE $U_{1}=q-L_{0}\left(1-\frac{c}{a}\right)$
\STATE $L_{1}=\frac{p}{U_{1}}$
\STATE $U_{2}=4-L_{1}\left(1-\frac{c}{a}\right)$
\FOR{$k=2:N-1$}
\STATE $L_k=\frac{p}{U_k}$
\STATE $U_{k+1}=q-L_k$
\ENDFOR
\STATE $L_{N+1}=\frac{d}{U_{N-1}}$
\STATE $L_{N}=\frac{e-L_{N+1}}{U_{N}}$
\STATE $d_{N+1}=f-U_{N}$

\RETURN $L,U$
\end{algorithmic}
\label{matrixnat}
\end{algorithm}

\subsubsection{Hermite coefficient computation for periodic conditions}
We compute the value of the distribution at the cell's faces $g_{k-1/2}$ and $g_{k+1/2}$.
With periodic boundary conditions, we decompose $A$ as $A=LD L^T$ with \algoref{matrixper} then we solve $Y=D^{-1}L^{-1}B$ with \algoref{coeffper}.
\begin{algorithm}[H]
\caption{CoefficientHermitePeriodic}
\begin{algorithmic}[1]
\REQUIRE{$\left(\bar{g}^{n}_k \right)_{k={0}\cdots N-1},\left( L_k\right)_{k=0 \cdots N},\left( U
D_k\right)_{k=0 \cdots N}$,a,b,c}
\STATE $X_{0}=3(\bar{g}^{n}_{N}+\bar{g}^{n}_{0})$
\FOR{$k=1:N-1$}
\STATE $X_{k}=3 \left( \bar{g}^{n}_{k}+\bar{g}^{n}_{k+1} \right)-L_{k-1}X_{k-1}$
\ENDFOR
\STATE S=0
\FOR{$m=0:N-2$}
\STATE $S=S+\delta_m X_m$
\ENDFOR
\STATE $X_{N}=3(\bar{g}^{n}_{N}+\bar{g}^{n}_{0})-sum^{}-L_{N-1}X_{N-1}$
\FOR{$m=0:N$}
\STATE $X_m=\frac{X_m}{D_m}$
\ENDFOR
\STATE $X_{N-1}=X_{N-1}-L_(N-1)X_N$
\FOR{$k=N-2:0,-1$}
\STATE $X_{k}=X_k-L_k X_{k+1}-\delta_k X_N$
\ENDFOR
\RETURN $\left(g_{k-1/2}\right)_{k={0}\cdots N}$
\end{algorithmic}
 \label{coeffper}
\end{algorithm}
We compute the matrix coefficients of the LDL decomposition.

\[
 \begin{bmatrix} A \end{bmatrix}=\begin{bmatrix} L \end{bmatrix}\begin{bmatrix} D \end{bmatrix}\begin{bmatrix} L^{T} \end{bmatrix}
\]
\[
 \begin{bmatrix}
q     & p      & 0      & 0      & \cdots & \cdots & 1      \\
p     & q      & p      & 0      & \cdots & \cdots & \vdots      \\
0     & p      & q      & p      & 0 & \cdots & \vdots \\
\vdots & \ddots & \ddots & \ddots & \ddots & \ddots & \vdots \\ 
\vdots & \cdots & 0 & p      & q      & p      & 0      \\
\vdots & \cdots & \cdots & 0      & p      & q      & p      \\
1     & \cdots & \cdots &0       & 0      & p     & q
 \end{bmatrix}=\]
 \[
 \begin{bmatrix}

1     & 0      & \cdots      & \cdots      & \cdots & \cdots & 0      \\
L_1     & 1      & 0      & \cdots      & \cdots & \cdots & \vdots      \\
0     & L_2      & 1     & 0      & \cdots & \cdots & \vdots \\
\vdots & \ddots & \ddots & \ddots & \ddots & \ddots & \vdots \\ 
\vdots & \cdots & 0 & L_{N-2}      & 1      & 0      & \vdots      \\
\vdots & \cdots & \cdots & 0      & L_{N-1}      & 1      & 0      \\
\delta_1     & \delta_2 & \cdots &  \cdots &    \delta_{N-2}  & L_{N-1}     & 1
 \end{bmatrix}
 \begin{bmatrix}
D_1     & 0      & 0      & 0      & \cdots & \cdots & 0      \\
0     & D_1      & 0 & 0      & \cdots & \cdots & \vdots      \\
\vdots     & 0      & D_2      & 0      & 0 & \cdots & \vdots \\
\vdots & \ddots & \ddots & \ddots & \ddots & \ddots & \vdots \\ 
\vdots & \cdots & \cdots & 0      &  D_{N-2}      & 0     & 0      \\
\vdots & \cdots & \cdots & \cdots     & 0      &  D_{N-1}      & 0     \\
0     & \cdots & \cdots &  \cdots       &\cdots     & 0     & D_{N}
 \end{bmatrix}
 \]
 
\begin{algorithm}[H]
\caption{MatrixCoefficientPeriodic}
\begin{algorithmic}[1]
\REQUIRE p,q,$i_{start}$,$i_{end}$
\STATE $D_1=q$
\STATE $L_1=\frac{p}{D_1}$
\STATE $\delta_1=\frac{1}{D_1}$
\FOR{$k=2:N-2$}
\STATE $D_k=q-p L_{k-1}$
\STATE $L_k=\frac{p}{D_k}$
\STATE $\delta_k=-\delta_{k-1}\frac{p}{D_k}$
\ENDFOR
\STATE $D_{N-1}=q-pL_{N-2}$
\STATE $L_{N-1}=\frac{p-p\delta_{N-2}}{D_{N-1}}$
\STATE $D_{N}=q-\sum^{N-2}_{j=1}D_i {\delta_i}^2-D_{N-1}L_{N-1}^2$
\RETURN $L,D$
\end{algorithmic}
 \label{matrixper}
\end{algorithm}

\end{document}